\magnification=\magstep1
\baselineskip=1.3\baselineskip

\expandafter\ifx\csname bookmacros.tex\endcsname\relax \else\endinput\fi
 
\expandafter\edef\csname bookmacros.tex\endcsname{%
       \catcode`\noexpand\@=\the\catcode`\@\space}
\catcode`\@=11

\def\undefine#1{\let#1\undefined}
\def\newsymbol#1#2#3#4#5{\let\next@\relax
 \ifnum#2=\@ne\let\next@\msafam@\else
 \ifnum#2=\tw@\let\next@\msbfam@\fi\fi
 \mathchardef#1="#3\next@#4#5}
\def\mathhexbox@#1#2#3{\relax
 \ifmmode\mathpalette{}{\m@th\mathchar"#1#2#3} 
 \else\leavevmode\hbox{$\m@th\mathchar"#1#2#3$}\fi}
\def\hexnumber@#1{\ifcase#1 0\or 1\or 2\or 3\or 4\or 5\or 6\or 7\or 8\or
 9\or A\or B\or C\or D\or E\or F\fi}

\font\tenmsa=msam10
\font\sevenmsa=msam7
\font\fivemsa=msam5
\newfam\msafam
\textfont\msafam=\tenmsa
\scriptfont\msafam=\sevenmsa
\scriptscriptfont\msafam=\fivemsa
\edef\msafam@{\hexnumber@\msafam}
\mathchardef\dabar@"0\msafam@39

\font\tenmsb=msbm10
\font\sevenmsb=msbm7
\font\fivemsb=msbm5
\newfam\msbfam
\textfont\msbfam=\tenmsb
\scriptfont\msbfam=\sevenmsb
\scriptscriptfont\msbfam=\fivemsb
\edef\msbfam@{\hexnumber@\msbfam}
\def\Bbb#1{{\fam\msbfam\relax#1}}
\def\widehat#1{\setbox\z@\hbox{$\m@th#1$} 
 \ifdim\wd\z@>\tw@ em\mathaccent"0\msbfam@5B{#1} 
 \else\mathaccent"0362{#1}\fi}
\def\widetilde#1{\setbox\z@\hbox{$\m@th#1$}
 \ifdim\wd\z@>\tw@ em\mathaccent"0\msbfam@5D{#1} 
 \else\mathaccent"0365{#1}\fi}
\font\teneufm=eufm10
\font\seveneufm=eufm7
\font\fiveeufm=eufm5
\newfam\eufmfam
\textfont\eufmfam=\teneufm
\scriptfont\eufmfam=\seveneufm
\scriptscriptfont\eufmfam=\fiveeufm

\newsymbol\centerdot 1205
\newsymbol\square 1003

\csname bookmacros.tex\endcsname

\def\R{{\Bbb R}}
\def\E{{\Bbb E}}
\def\P{{\Bbb P}}
\def\Q{{\Bbb Q}}

\def\F{{\cal F}}

\def\lam{{\lambda}}

\def\al{{\alpha}}

\def\proof{{\medskip\noindent {\bf Proof. }}}
\def\longproof#1{{\medskip\noindent {\bf Proof #1.}}}
\def\qed{{\hfill $\square$ \bigskip}}

\def\section#1#2{{\bigskip\bigskip \centerline{\bf #1. #2}\bigskip}}

\def\chapter#1#2{{\bigskip\bigskip \centerline{#1. #2}\bigskip}}
\def\cite#1{{[#1]}}

\def\eps{\varepsilon}

\def\norm#1{\Vert #1 \Vert}

 \def\qq {\qquad}
\def\frac#1#2{{#1\over #2}}

\def\wt{\widetilde}
\def\ol{\overline}
\def\wh{\widehat}

\def\ni{\noindent }
\def\ms{\medskip}
\def\bs{\bigskip}

\def \half {{{1/ 2}}}

\parindent=30pt

\def\ftN{\rl\kern-0.13em\rN}

\def\Cov{{\mathop {{\rm Cov\, }}}}

\def\sqr#1#2{{\vcenter{\vbox{\hrule height.#2pt
        \hbox{\vrule width.#2pt height#1pt \kern#1pt
           \vrule width.#2pt}
        \hrule height.#2pt}}}}

\def\var{{\rm Var\,}}

\def\sign{\mathop{\rm sgn}}
\def\df{{\mathop {\ =\ }\limits^{\rm{df}}}}
\def\bone{{\mathop{\hbox{\bf 1}}}}

\def\sqr#1#2{{\vcenter{\vbox{\hrule height.#2pt
        \hbox{\vrule width.#2pt height#1pt \kern#1pt
           \vrule width.#2pt}
        \hrule height.#2pt}}}}
\def\squareinf{\quad\mathchoice\sqr56\sqr56\sqr{2.1}3\sqr{1.5}3}

\input epsf
\newdimen\epsfxsize

\centerline{\bf STOCHASTIC BIFURCATION MODELS}

\footnote{$\empty$}{\rm Research partially supported by NSF grant DMS-9700721.}

\vskip0.4truein
\centerline{{\bf Richard F.~Bass}
\qq and \qq
{\bf Krzysztof Burdzy}}

\vskip0.4truein

\noindent{\bf Abstract}.
We study an ordinary differential equation controlled by a stochastic
process.
We present results on existence
and uniqueness of solutions, on associated local
times (Trotter and Ray-Knight theorems), and on time
and direction of bifurcation. A relationship with Lipschitz
approximations to Brownian paths is also discussed.

\bigskip

\vfill\eject

\noindent{\bf 1. Introduction}.

Let $B_t$ be a continuous function of $t$, let $t_0, x_0, \beta_1,
\beta_2\in \R$,  and consider the 
ordinary differential equation

$$\frac{dX_t}{dt} = \cases{ \beta_1 & if $X_t < B_t$,
\cr \beta_2 & if $X_t > B_t$, \cr} \qquad t \in \R, \qq X(t_0) = x_0. \eqno(1.1)$$

\ni Among the results we prove are the following:

(1) Although in general there will not be a unique solution to (1.1),
there will be a unique Lipschitz solution to (1.1) 
if $B_t$ is a typical Brownian motion path.

(2) Let $B_t$ be a Brownian motion with $B_0=0$
and let $X_t^{x_0} $ denote the solution to 
(1.1) when $t_0=0$ and $X(t_0)=x_0$. The map $y\to X_t^y$ is
a one-to-one map of $\R$ onto $\R$. The smoothness of this
map is controlled by the local time at 0 of $X_t^{y}-B_t$. If
we call this local time $L_t^{y}$ and $\beta_1, \beta_2$
satisfy suitable assumptions, then $L_t^y$ is jointly
continuous in $y$ and $t$ and 
$\{L_\infty^y, y\geq 0\}$ and
$\{L_\infty^{-y}, y\geq 0\}$ are strong Markov processes.
We show that this implies that for a fixed $t>0$,
the function $y \to X^y_t$ is of class $C^{1+\gamma}$
with $\gamma < 1/2$, but it is not $C^{3/2}$.

(3) As we shall see below, (1.1) is an example of a bifurcation
model; if $B_t$ is a Brownian motion,
$\beta_1<0$ and $\beta_2>0$,
each of the events $\{\lim_{t\to \infty} X_t=+\infty\}$
and $\{\lim_{t\to \infty} X_t=-\infty\}$ has positive probability.
The bifurcation time is defined by
$= \sup\{t: X_t = B_t\}$. 
We calculate both the probability of $\{\lim_{t\to \infty} X_t=+\infty\}$ and 
the expectation of the bifurcation time using excursion theory.

(4) The equation (1.1) sheds light on the best Lipschitz
approximation to Brownian paths. In particular we obtain 
an estimate on the lower
bound on the best constant in the Koml\'os-Major-Tusn\'ady result
concerning strong approximations of Brownian motion by random walks.

\bs

Equation (1.1) is similar to 
an equation that
arose in the course of an economic study
and its accompanying probabilistic model
in Burdzy, Frankel, and Pauzner (1997, 1998).
These papers 
introduce and study an economics model whose technical side
is based on the following equation:
$$\frac{dX_t}{dt} =
\cases{ - \beta X_t & if $X_t < f(B_t)$,
\cr \beta(1-X_t) & if $X_t > f(B_t)$, \cr}\qquad t \geq 0, \qquad X(0)
=x_0\in (0,1),\eqno(1.2)$$
where $B_t$ is a Brownian motion starting from $B_0 = b_0$,
$\beta >0$ is a fixed constant, and $f$ is a
non-increasing Lipschitz function.
The case when $x_0 = f(b_0)$ is of special interest.
Results on the time and direction of the stochastic
bifurcation were crucial elements of these two  papers.

We also consider the following equation,
more general than (1.1).
$$
\frac{dX_t}{dt} =  \cases{  \beta_1 |X_t - B_t|^{\alpha_1} & if $X_t < B_t$,
\cr \beta_2  |X_t - B_t|^{\alpha_2}& if $X_t > B_t$, \cr}
\qquad t \in \R, \qq X(t_0)=x_0.\eqno(1.3)
$$

\ni (If $\al_1=\al_2=0$, then (1.3) reduces to (1.1).)
Equation (1.3)
was inspired by the following 
model. Consider a pendulum
with rigid arm which is turned upside down (see Fig.~1.1).

\bigskip
\vbox{
\epsfxsize=2.0in
  \centerline{\epsffile{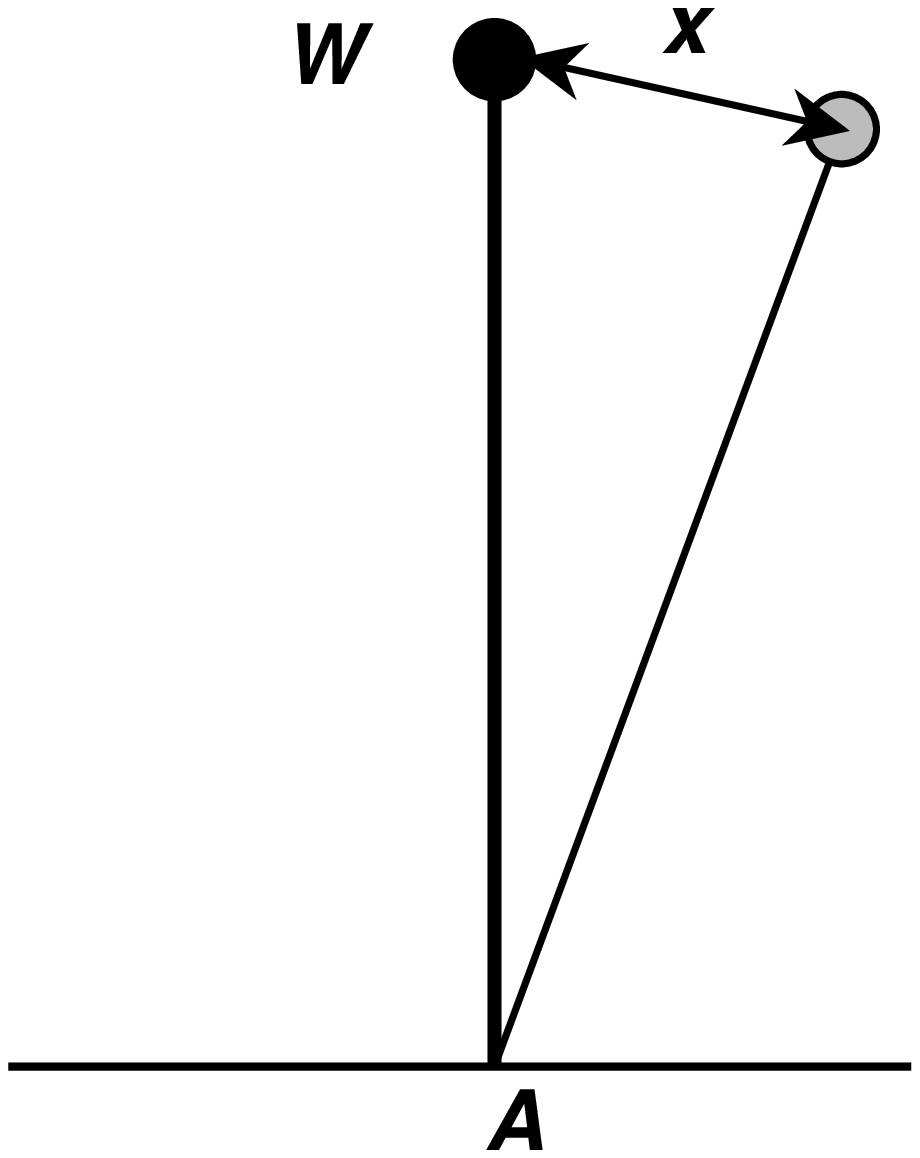}}

\centerline{Figure 1.1.}
}
\bigskip

Let $X_t$ denote the distance of the weight
$W$ from its unstable rest position at the top
of the vertical arm.
When $X_t= x$ and $x $ is small, the weight is about $c_1 x^2$
units below its rest position and, therefore
$c_2 x^2$ units of potential energy must have
been converted to kinetic energy, given by
$c_3 (dX/dt)^2$. Hence, we have the approximate
relationship $dX/dt = c_4 X_t$, assuming infinitesimally
small velocity at the rest position. Note that if the initial
velocity at the rest position is close to zero, then
the time it takes the pendulum to move any fixed
non-zero distance from the rest position is very large.
We now add stochastic oscillations to our pendulum model.
We suppose that the
base $A$ of the pendulum vibrates according to  a Brownian motion $B_t$.
Then the position $X_t$ of the weight $W$ relative
to $A$ is $X_t - B_t$ and we have 
$dX/dt = c_4 (X_t - B_t)$, which is (1.3)
with $\alpha_1=\alpha_2 =1$ and $-\beta_1=\beta_2= c_4$.
\bs

The solutions to (1.1) exhibit fast switching between two kinds
of excursions.
See Karatzas and Shreve (1988, Sect.~6.5) for a 
closely related model. Mandelbaum, Shepp, and Vanderbei (1990) 
also consider a model with fast switching between two kinds of excursions,
but we were not able to find a direct connection
with our own model.
\bs

The rest of the paper consists of five sections.
Section 2 contains results on existence and uniqueness
of solutions to (1.1), (1.3), and related equations. 
The process $B_t$ will generally be a Brownian motion,
but Theorems 2.3 and 2.4 also apply to some fractional Brownian motions
(see Examples 2.3 and 2.4).

Let $X^y_t$ denote the solution to (1.1)
with $X^y_0 =y$. For a fixed $t\geq 0$,
the function $y \to X^y_t$ is a transformation
of $\R$ onto  itself. 
How smooth is this map? How many derivatives does the function
$y \to X^y_t$ have and are they continuous? To answer these
questions, one is  led to study the local time of $X_t^y-B_t$.
Section 3 is devoted to  a number of
results about local times related to (1.1), including
analogues of the Trotter and Ray-Knight theorems.
See Knight (1981), Leuridan (1998), 
Norris, Rogers and Williams (1987),
Revuz and Yor (1991) and Yor (1997)
for old and new variants of the Ray-Knight theorem.
Our local times are defined as local times
at points, but they may 
also be  viewed as local times of Brownian motion
on a random curve---see (5.15) in
F\"ollmer, Protter, and Shiryaev (1995)
for a result on local times on non-random curves.

Section 4 gives explicit formulae for the
probability of upward bifurcation for the equation (1.3)
and the expected bifurcation time for (1.1), with some
indication how to proceed in the more general case (1.3).
This extends results from 
Burdzy, Frankel and Pauzner (1998).
Section 5 takes a look at the solutions
to (1.1) as Lipschitz approximations to the Brownian path.
As a consequence we obtain some lower bounds related to the
Koml\'os-Major-Tusn\'ady construction; see Theorem 5.7.
Finally, Section 6 is a list of open problems.

In Sections 3-5, we consider
Brownian motion defined
on the whole real line $\R$, i.e., the process 
$\{B_t, -\infty < t < \infty\}$, where 
$\{B_t, t \in (0,\infty)\}$ and $\{B_{-t}, t \in (0,\infty)\}$
are independent Brownian motions starting from $0$
with variance $\E  B^2_t = \E  B^2_{-t} = \sigma^2 t$.
Unless stated otherwise, we will assume that all Brownian motions
(including those with drift and/or reflection)
have infinitesimal variance $\sigma^2$,  and that all constants
are strictly positive and finite.

Section 3 of the paper was inspired by unpublished heuristic
calculations involving local times which were a part
of an earlier project of David Frankel, Ady Pauzner, and the second author.
We would like to thank the many colleagues who kindly gave
us advice on various aspects of the model: Robert Adler, Ludwig Arnold,
Jean Bertoin,
Miklos Cs\"org\"o, Burgess Davis, Laurent Decreusefond, David Frankel, 
Mike Harrison, Haya Kaspi, Frank Knight, Jim Kuelbs, Avi Mandelbaum, Ady Pauzner, 
Jim Pitman, Philip Protter, Emmanuel Rio,
Ruth Williams, Marc Yor, and Ofer Zeitouni.

\bigskip

\bigskip
\noindent{\bf 2. Existence and uniqueness of solutions}.
In this section we present several theorems on the existence and uniqueness
of solutions to differential equations similar to (1.1). There is considerable
overlap among the theorems, but each contains cases not covered
by the other.
We first present our main results. They are followed by
some remarks and examples. The proofs are relegated to the end
of the section.

We start with the equation

$$
\frac{dX}{dt} =  \cases{  \beta_1 |X_t - B_t|^{\alpha_1} & if $X_t < B_t$,
\cr \beta_2  |X_t - B_t|^{\alpha_2}& if $X_t > B_t$, \cr}
\qquad t \in \R, \qq X(t_0)=x_0.\eqno(2.1)
$$

\ni where $B_t$ is a Brownian motion, $\al_1, \al_2>-1$, and $\beta_1, \beta_2\in \R$.

First note that the function $X_t =B_t$ is a 
solution to (2.1) with $t_0=0$ and $x_0 = 0$,
because neither of the conditions on the right
hand side of (2.1) is ever satisfied. 
We would like to disregard such a solution for two
reasons. First, the economics model behind (1.2)
required that the solutions to (1.2) be Lipschitz.
Second, the example $X_t = B_t$ is rather artificial.
For $\alpha_1, \alpha_2\geq 0$ it is natural to require
that $X_t$ is a Lipschitz function.
We generalize this to all $\alpha_1, \alpha_2 > -1$
by writing an integrated version of (2.1), namely,
$$X_t  = x_0 +
\int_{t_0}^t \left[
\beta_1  |X_s-B_s|^{\alpha_1} \bone_{\{X_s - B_s \leq 0\}} 
+ \beta_2  |X_s-B_s|^{\alpha_2}\bone_{\{X_s - B_s > 0\}}
\right] ds. \eqno(2.2)$$
It is easy to see that solutions to (2.2) satisfy (2.1),
but the example $X_t = B_t$ shows that the opposite
statement is not true.

\bigskip
\proclaim Theorem 2.1.  
For fixed
$t_0,x_0,\beta_1,\beta_2\in \R$, $\sigma^2>0$, and
$ \alpha_1,\alpha_2>-1$, there exist a Brownian motion
$B_t$ and a process $X_t$ which satisfy (2.2) with
the initial condition as in (2.1). The solution
$X_t$ is unique in law. We may construct $X_t$
in such a way that $(X_t,B_t)$ is a strong Markov
process relative to the appropriate filtration.
If we assume in addition that 
$ \alpha_1,\alpha_2\geq 0$, then for a given
Brownian motion $B_t$ there exists a unique
solution to (2.2), a.s.

\bigskip

Our next theorem is a result on existence.
We will state the result for the following generalization
of the equation (1.1),
$$
\frac{dX_t}{dt} =  \cases{  F_1(X_t)& if $X_t > B_t$,
\cr F_2(X_t)& if $X_t < B_t$, \cr}
\qquad t \in \R,\qq X(t_0) = x_0.\eqno(2.3)
$$

\bigskip
\noindent{\bf Theorem 2.2}. {\sl
Assume that 
$F_1$ and $ F_2$ are continuous functions and that
$|F_1|$ and $| F_2|$ are bounded by $\beta< \infty$.
If $B_t$ is a continuous process, then (2.3) has a Lipschitz solution, a.s. 
There exists a maximal Lipschitz solution $\{X_t^+, t\geq t_0\}$
to (2.3); it is adapted to the
filtration ${\cal F}_t = \sigma(B_s, s\in [t_0, t])$.
}
\bigskip

Haya Kaspi pointed out to us that measurability of a solution
to (2.3) is the most delicate point of Theorem 2.2.

\bigskip

We will say that $L^x_t$ is a local time for a process $B_t$
if it is the occupation time density:
$$\int_{-\infty}^\infty h(x)L^x_t\, dx=\int_0^t h(B_t)\, dt, \qq {\rm a.s.},$$
for all $h$ bounded and measurable.
Note that if $B_t$ is continuous
and the local time $L^x_t$ is jointly continuous, then
$\sup_x L^x_t<\infty$, a.s. for each $t$.

We will use the traditional
Markovian notation $\P^x$ to denote the distribution
of $\{B_t, t\geq t_0\}$ conditioned by $\{B_{t_0}=x\}$,
even though we do not assume the Markov property
for $B_t$ in Theorems 2.3 and 2.4 below.

\bigskip
\noindent{\bf Theorem 2.3}. {\sl Let $t_0>0$,
$ x_0, \beta_1, \beta_2\in \R$.
Assume that 
\item{(i)} the process $B_t$ is continuous
and has a jointly continuous local time  $L^x_t$, and
\item{(ii)} if $A_t$ is an adapted process with $A_{t_0}=x_0$
whose paths are Lipschitz
continuous with Lipschitz constant $M$, then for each $x$ the law
of $\{B_t+A_t, t_0\leq t \leq t_0 + s\}$ under $\P^x$ is mutually absolutely continuous
with respect to the law of $\{B_t, t_0\leq t \leq t_0 + s\}$ under $\P^{x+x_0}$,
for every $s>0$.
\par
\ni Then with probability one there exists a random $s_0>0$
and a unique Lipschitz solution to (1.1) on $[t_0, t_0 + s_0]$.

If in addition we assume that $B_t$ is strong Markov 
then there is a unique Lipschitz solution to (1.1)
for all $t\geq t_0$.}

\bigskip
\ni{\bf Remark 2.4.} If $W_t$ is a Brownian motion and $f$ is a strictly increasing
function such that both $f$ and $f^{-1}$ are Lipschitz continuous,
it is easy to check that $B_t=f(W_t)$ is a strong Markov process that
satisfies the other assumptions of Theorem 2.3.

\bigskip
\noindent{\bf Theorem 2.5}. {\sl Let $t_0,x_0\in \R$. Assume that 
$F_1$ and $ F_2$ are bounded,  Lipschitz functions. Suppose that
both are bounded by $M$ and that both have Lipschitz constant less
than or equal to $M$.
Let $B_t$ be a continuous process such that
\item{(i)} there exist $c_1>0$ and
$\gamma\in (0,1)$ such that  whenever $s<t$,
$$\P(B_t\in dy\mid\F_s)\leq\frac{c_1}{(t-s)^\gamma}\, dy,\qq y\in \R,
\eqno(2.4)$$
\item{(ii)} if $A_t$ is an adapted process with $A_{t_0}=x_0$
whose paths are Lipschitz
continuous with Lipschitz constant $M$, then for each $x$ the law
of $\{B_t+A_t, t_0\leq t \leq t_0 + s\}$ under $\P^x$ is mutually absolutely continuous
with respect to the law of $\{B_t, t_0\leq t \leq t_0 + s\}$ under $\P^{x+x_0}$,
for every $s>0$.

Then with probability one, there exists a unique solution to 
(2.3) for all $t\geq t_0$.
}

\bigskip

We will show in Example 2.10 below that Theorem 2.5 applies
to some fractional Brownian motions.
As in Remark 2.4,
some functions of fractional Brownian motions also satisfy the hypotheses of
Theorem 2.5.

\bigskip
Let $f(x,b) = \beta_1 \bone_{\{x \leq b\}} + \beta_2 \bone_{\{x>b\}}$
and suppose that $\alpha_1=\alpha_2=0$.
Then (2.2) may be written as
$$X_t  = x_0 +
\int_{t_0}^t f(X_t, B_t) ds. \eqno(2.5)$$
The function $(x,b)\to f(x,b)$ is discontinuous.
In applications,
such as that in Burdzy, Frankel, and Pauzner (1997),
it may be argued that a model with continuous
$dX/dt$ might be more realistic. 
Let us replace $f$ with a continuous approximation,
$$f_\eps(x,b) = \beta_1 \bone_{\{x < b-\eps\}} 
+ \beta_2 \bone_{\{x>b+\eps\}} +
\left[ {\beta_2 - \beta_1 \over 2 \eps}
(x - b + \eps) + \beta_1\right]
\bone_{\{b-\eps \leq x \leq b+\eps\}},$$
and consider the corresponding equation
$$X^\eps_t  = x_0 +
\int_{t_0}^t f_\eps(X^\eps_t, B_t) ds. \eqno(2.6)$$
We will show that the solutions to (2.6) converge
to those of (2.5), and thus many results about solutions
to (2.5) proved later in this article may be applied to 
give asymptotic results for the solutions to (2.6). 

\bigskip
\noindent{\bf Theorem 2.6}. {\sl
Assume that the equations (2.5) and (2.6)
are defined relative to the same Brownian motion $B_t$.
The equation (2.6) has a unique Lipschitz solution.
As $\eps \to 0$, the functions $X^\eps_t$
converge to the unique solution $X_t$ of (2.5), a.s.
}

\bigskip
Note that the convergence in Theorem 2.6 is uniform
on compact sets as all functions $X^\eps_t$ are Lipschitz
with constant $\max\{|\beta_1|, |\beta_2|\}$.

\bigskip
\noindent{\bf Remark 2.7}. For the economics
model behind (1.2), one does not necessarily want to require the Markov property to 
hold. The proof of Theorem 2.3 uses the strong
Markov property to do an induction argument. For Theorem 2.5
we  have in mind examples
where $B_t$ is a Gaussian process; see Example 2.10 below.
In general, $B_{T+t}-B_T$ will not be Gaussian
when $T$ is a stopping time. 

\bigskip
\noindent{\bf Example 2.8}.
We present an elementary example of
a continuous deterministic function $t \to B_t$
for which there are multiple solutions to (1.1).
Let $\beta_1 <0$, $\beta_2 >0$,
$$B_t= \cases{(1 + \beta_2) t & for $t \in[0,1]$, \cr
1+\beta_2 & for $ t >1$, \cr
0 & for $ t < 0$. } $$
There are uncountably many solutions to (1.1) with this
choice of $B_t$ and the initial condition
$X_0 = 1$. Here are two of them:
$$X^1_t = \cases{ 0 & for $t \leq - 1/\beta_2$, \cr 
1 +  \beta_2 t & for $ t > - 1/\beta_2$ ; }$$
$$X^2_t = \cases { 0 & for $t \leq - 1/\beta_2$, \cr 
1 +  \beta_2 t & for $ t \in( - 1/\beta_2, 1]$ \cr
1+ \beta_2 & for $t \in (1, 5]$, \cr
1 + \beta_2 + 5 \beta_1 + \beta_1 t & for $t >5$. }
$$
\bigskip

\noindent{\bf Example 2.9}. As we noted earlier in this section,
$X_t=B_t$ is a solution to (1.1) but a rather trivial one.
In this example, we will show a less trivial and perhaps
more interesting non-Lipschitz solution to (1.1).
Take $\beta_1=\beta_2=0$ in (1.1); in other words, consider
the equation
$$\frac{dX_t}{dt} = 0 \qquad \hbox{  if  }\qquad X_t \ne B_t
 \qquad t \in \R, \qq X(t_0) = x_0. $$
The function $X_t =0$ is a solution to this equation
and, moreover, it is the only Lipschitz solution,
by Theorem 2.1.
Let $Y_t$ be a skew Brownian motion, i.e., a process
which may be constructed by flipping positive
excursions of a standard Brownian motion $\wt B_t$ to the negative side
with probability $p_1$ and negative excursions
to the positive side with probability $p_2$, independently
of each other. Suppose that $p_1\ne p_2$ so that
the process $Y_t$ is not a standard Brownian motion.
Let $L_t$ be the local time of $Y_t$ at $0$.
By a result of Harrison and Shepp (1981) (see
also Exercise X (2.24) in Revuz and Yor (1991)), for a suitable constant
$c_1 \ne 0$, the process $Y_t - c_1 L_t$ is a standard
Brownian motion. If we take $B_t = Y_t - c_1 L_t$ then
$X_t = c_1 L_t$ is a non-Lipschitz solution to our equation.

\bigskip
\noindent{\bf Example 2.10}.
We provide an example of a process satisfying
the assumptions of Theorem 2.5 that is not strong Markov. Let $B_t$ be fractional Brownian motion
of index 
$H\in (0,1/2]$. This means that $B_t$ is a mean zero Gaussian
process with
$$\Cov(B_s,B_t)=c_1(s^{2H}+t^{2H}-|t-s|^{2H}).$$
$B_t$ has a stochastic integral representation
$$B_t=\int_{-\infty}^t R(t,u)\, dZ_u,$$
where $Z_u$ is a standard Brownian motion and
$$R(t,u)=c_2[((t-u)^+)^{H-1/2}-(u^-)^{H-1/2}];$$ 
see, e.g., Rogers (1997). Conditioning
on $\F_s$ with $s>0$, the law of $B_t$ given $\F_s$ is that of a Gaussian
process with variance
$$c_2^2\E\Big[\Big(\int_s^t (t-u)^{H-1/2}\, dZ_u\Big)^2\mid \F_s\Big]
=c_2^2\int_s^t (t-u)^{2H-1}\, du=c_3(t-s)^{2H}.$$
Assumption (i) of Theorem 2.5 is immediate from this. 

We now show (ii). We give the argument for the case $t_0=x_0=0, s=1$;
the extension to the general case is routine. 

If $H=1/2$, then $B_t$ is standard Brownian
motion, and (ii) follows from the Girsanov theorem; so we suppose
$H<1/2$. Let $\al=H+1/2$. 
See Decreusefond and \"Ust\"unel (1997) for more details of
some of the steps in the following  argument.
Let
$F(a,b,c,z)$ be the standard Gauss hypergeometric function
and
define an operator $K_H$ on functions on $[0,1]$ by
$$\eqalign{(K_H&f)(t)\cr
&=\frac{1}{\Gamma(H+1/2)}\int_0^t (t-x)^{H-1/2}F(H-1/2,1/2-H, H+1/2,
1-t/x)f(x)dx.\cr}$$
Let ${\cal H}_H=\{K_Hh: h\in L^2([0,1])\}$ and define
$$\norm{f}_{{\cal H}_H}=\norm{K_H^{-1}f}_{L^2}.$$
For $\beta\in(0,1)$ define
$$(I^\beta f)(x) =\frac{1}{\Gamma(\beta)}\int_0^x f(t)(x-t)^{\beta-1}dt$$
and 
$$(D^\beta f)(x)=\frac{d}{dx}\Big(I^{1-\beta}f\big)(x).$$
By 
Decreusefond and \"Ust\"unel (1997) (Theorem 2.1, Theorem 3.3, and the proof of Theorem 3.3),
we have that ${\cal H}_H$ is
dense in the set of continuous functions on $[0,1]$ that are null at 0
and that $K_H$ is an isomorphism from $L^2([0,1])$ onto $I^{H+\half}(L^2([0,1]))$.
By Proposition 2.1 of that paper, $D^\beta$ is the inverse to  $I^\beta$.

Since $K_H^{-1}$ is continuous from $I^{H+\half}(L^2)$ into $L^2$,
then $K_H^{-1}\circ I^{H+\half}$ is continuous from 
$L^2$ into itself, and so there exists a constant $c_4$
such that
$$\norm{K_H^{-1}I^{H+\half}g}_{L^2}\leq c_4 \norm{g}_{L^2}.$$
Thus if 
$f\in {\cal H}_H$, then
$$\norm{K_H^{-1}f}_{L^2}\leq c_4\norm{D^{H+\half}f}_{L^2},$$
or
$$\norm{f}_{{\cal H}_H}\leq c_4 \norm{D^\al f}_{L^2}.$$

Let $A_t$ be  a uniformly Lipschitz process  as in  the statement of
Theorem 2.5. By Theorem 4.9 of
Decreusefond and \"Ust\"unel (1997) and the Novikov
condition discussed just after that theorem, (ii) will hold if for
each $T\in(0,1)$ we have
$$\E\exp[\norm{A(\cdot)}_{{\cal H}_H}^2/2]<\infty.$$
By the above paragraph, it is enough to show
$$\E \exp\Big(\int_0^T |D^\al A_t|^2\, dt/2\Big)<\infty. \eqno (2.7)$$
To show (2.7), by an approximation argument it
suffices to show that for each fixed $T> 0$ there exists
$c_5$ (depending on $T$) such that  if $f$ is a $C^\infty$ function on
$[0,\infty)$ with $f(0)=0$, then
$$\sup_{0\leq t\leq T} |D^\al f(x)|\leq c_5 \norm{f'}_\infty; \eqno (2.8)$$
(2.7) will then follow easily from (2.8) and our assumptions on $A_t$.

Note that by a change of variables,
$$I^{1-\al}f(x)=c_6\int_0^x f(x-t)t^{-\al}dt,$$
and by the Leibniz formula and the fact that $f(0)=0$,
$$\frac{d}{dx} I^{1-\al}f(x)=c_6\int_0^x f'(x-t)t^{-\al}dt
=c_6\int_0^x f'(t)(x-t)^{-\al}dt=I^{1-\al}f'(x).$$
Since $\al=H+1/2<1$, then  $|x-t|^{-\al}$ is integrable on $[0,x]$.
So, for $u=f'$,
$$|D^\al f(x)| =
|I^{1-\al}u(x)|\leq \norm{u}_\infty\int_0^x |x-t|^{-\al}\, dt
\leq c_7\norm{u}_\infty$$
for $x\leq T$. This gives (2.8), and thus
a fractional Brownian motion with parameter $H\in(0,1/2]$
satisfies the assumptions of Theorem 2.5.

\bs

\bigskip
\noindent{\bf Example 2.11}. The weaker version of Theorem 2.3,
i.e., the one
without the assumption on the Markov character of $B_t$,
applies to fractional Brownian motions with 
parameter $H\in(0,1/2]$.
Assumption (ii) of Theorem 2.3 is the same as (ii)
of Theorem 2.5; we have verified that assumption in the previous
example. As for assumption (i) of Theorem 2.3,
the joint continuity of the local time for the
fractional Brownian motion follows from Lemma 8.8.1,
Theorem 8.8.2 and the proof of Theorem 8.8.4
in Adler (1981).

\bs
\ni{\bf Example 2.12.} Fabes and Kenig (1981) gave an example of a process
$B_t$ satisfying 
$$dB_t= \sigma(B_t,t)\, dW_t,$$
where $W_t$ is a standard Brownian motion, $\sigma$ is H\"older continuous in 
the first variable, $\sigma$ is bounded above and below by positive constants,
and the distribution of $B_1$ does not have a density with respect to Lebesgue measure.
$B_t$ is a space-time strong Markov process. Because  $\sigma$ is bounded below,
it is not hard to see that $B_t$ has a jointly continuous local time (cf.\ Revuz and Yor (1991), Ch.~6) and that  hypothesis (ii) of Theorem 2.3
holds. Thus this process $B_t$ is an example where
the assumptions of Theorem 2.3 hold, but those of Theorem 2.5 do not.

\bigskip
The rest of the section contains proofs of our main
results. The following lemma is immediate.

\bigskip
\proclaim Lemma 2.13. 
Let $\wt B_t = B_{-t}$ and $\wt X_t = X_{-t}$.
If $X_t$ is a solution to (2.1) then $\wt X_t$
is a solution to
$$
\frac{d\wt X}{dt} =  \cases{- \beta_1 |\wt X_t - \wt B_t|^{\alpha_1} &
if $\wt X_t < \wt B_t$,
\cr -\beta_2  |\wt X_t - \wt B_t|^{\alpha_2}& if $\wt X_t > \wt B_t$, \cr}
\qquad t \in \R, \qq \wt X(-t_0) = x_0.
$$

\longproof{of Theorem 2.1}
For simplicity, assume that $t_0=0$.
The equation
$$Y_t  = x_0 + \int_0^t \left[
\beta_1 |Y_s|^{\alpha_1}\bone_{\{Y_s \leq 0\}} 
+ \beta_2  |Y_s|^{\alpha_2}\bone_{\{Y_s > 0\}}\right] ds 
- \int_0^t dB_s ,\qquad t \geq 0,$$
has a weak solution which is unique in law by
Theorem 5.15 in Karatzas and Shreve (1988). 
For $X_t = Y_t + B_t$, the last equation
is equivalent to (2.2) for $t\geq 0$. This proves the first assertion
of the theorem. The strong uniqueness in the case
$\alpha_1,\alpha_2\geq0$ follows from
Proposition 5.17 of Karatzas and Shreve (1988).
We note that although the function $y\to y^\alpha$
is not bounded, that proposition clearly applies
by using a truncation argument. The part of the solution to (2.1)
for $t < t_0=0$ can be obtained in a similar way
using Lemma 2.13.
That $X_t$ may be constructed so that $(X_t, B_t)$ is a strong
Markov process follows from the weak uniqueness in a standard
manner; see Bass (1997), Section I.5, or Stroock and Varadhan
(1979), Chapter 6.
\qed

\bigskip
\longproof{of Theorem 2.2} We start by showing that 
for each $\omega$ and
for any $u_1$ and $z_1$ there exists
a maximal solution $\wt X^{u_1, z_1}_t$ to the equation
$$dX_t/dt = F_1(X_t), \qquad t \in \R , \qquad 
X(u_1) = z_1.$$
First of all, it is well known that there exists
at least one solution to the equation since $F_1$
is continuous. Since $|F_1|$ is bounded by $\beta$,
all solutions are Lipschitz with constant $\beta$
and so their supremum $\wt X^{u_1, z_1}_t$ is also a Lipschitz function
with constant $\beta$. Next note that the maximum
of any two solutions is also a solution to the equation.
This and the Lipschitz property of solutions
easily imply that there exists a sequence
of solutions converging to $\wt X^{u_1, z_1}_t$,
uniformly on compact intervals. Now a standard argument
can be used to show that $\wt X^{u_1, z_1}_t$ is a solution
to the equation.

The analogous maximal solution to $dX_t/dt = F_2(X_t)$
with the initial condition $X(u_1) = z_1$ will be denoted
$\wh X^{u_1, z_1}_t$.

We start by proving the existence of a solution to (2.3)
for $t \geq t_0$. 
Consider a small $\delta>0$. We proceed to  define a $\delta$-approximate
solution $X_t^\delta$ to (2.3).
First suppose that $B_{t_0} < x_0$. 
By the continuity of the paths of $B_t$, for almost every path
of $B_t$, there exist a unique time $t_1 \in ( t_0,\infty]$
and a function $X^\delta_t$ defined for $t \in (t_0,t_1)$,
such that $X^\delta_{t_0}=x_0$,
$X^\delta_{t_1} = B_{t_1}$ if $t_1 < \infty$,
and $X^\delta_t = \wt X^{t_0, x_0}_t$ 
for all $t \in (t_0,t_1)$.
We then let
$X^\delta_t = X^\delta_{t_1} + \beta(t-t_1)$ 
for all 
$t \in [t_1, t_1 +\delta]$, if $t_1 < \infty$.
If $B_{t_0} > x_0$ we use the same procedure to define
$X^\delta_t$ for $t \in[t_0, t_1+\delta]$
except that we use the function 
$\wh X^{t_0, x_0}_t$ in place of $\wt X^{t_0, x_0}_t$.
If $B_{t_0} = x_0$, we let $t_1 = t_0$ and 
$X^\delta_t = X^\delta_{t_1} + \beta(t-t_1)$ for 
$t \in [t_1, t_1 +\delta]$.

We have defined $X^\delta_t$ on an interval
$[t_0, t_1 + \delta]$. Let $x_1 = X^\delta_{t_1 + \delta}$.
Let us replace the initial condition in (2.3)
by $X(t_1 + \delta) = x_1$ and define an approximate
solution $X^\delta_t$ to (2.3) on an interval $[t_1 + \delta, t_2+\delta]$
using the same method as above. By induction,
we can construct a (possibly infinite) sequence
of times $\{t_k\}$ and a continuous function $X^\delta_t$
which satisfies (2.3) on every interval
$(t_k+\delta, t_{k+1} )$ and which is linear
on every interval $[t_k, t_k +\delta]$, for $k\geq 1$.
Note that the function $X^\delta_t$ is defined
for all $t\geq t_0$ because $t_{k+1} \geq t_k + \delta$
for every $k$. 

By construction, the $\delta$-approximate
solution $X^\delta_t$ is a Lipschitz function with 
Lipschitz constant $\beta$.  

For every integer $m\geq 1$
consider a $1/m$-approximate solution $X^{1/m}_t$. All of these functions
are Lipschitz with the same constant $\beta$,
and they all satisfy $X^{1/m}_{t_0} = x_0$.
Let $X_t$ be defined by
$$X_t = \limsup_{m\to \infty} X^{1/m}_t =
\lim_{n\to \infty}  \sup_{m > n} X^{1/m}_t.$$

The supremum of an arbitrary family of Lipschitz
functions with constant $\beta$ is a Lipschitz function
with the same constant, and the same remark applies
to the limit of a sequence of such functions.
Hence, for every $n$, the function
$Y^n_t= \sup_{m > n} X^{1/m}_t$ is Lipschitz with constant $\beta$,
and the same is true of $X_t$.
Note that $Y^n_t$ converge in a monotone way to $X_t$,
uniformly on compact intervals, because all these
functions are Lipschitz with the same constant $\beta$.

We will show that $X_t$ is a solution to (2.3). 
Let 
$$W(\delta) = \bigcup_{\{(s,x): s \geq t_0, B_s = x\}}
\{ (t,y): y = x + (t- s) \beta, t\in[s,s+\delta]\}.$$ 
For $\delta \leq \delta_1$,
the portion of the graph of $X^\delta_t$ which lies
outside $W(\delta_1)$ satisfies (2.3), by construction.

The set of $t$ such that
$B_t = X_t$ is closed because both functions
$B_t$ and $X_t$ are continuous. Consider any interval
$(s_1,s_2)$ such that $B_t \ne X_t$ for all
$t \in (s_1,s_2)$. 
Suppose without loss of generality that
$B_t < X_t$ for all $t \in (s_1,s_2)$
Choose an arbitrarily small $\delta_1 >0$.
Note that as $\delta \to 0$, the open sets
$W^c(\delta)$ converge to the complement of
$\{(s,x): s \geq t_0, B_s = x\}$.
Let 
$\delta_2>0$ be so small that the (closed) portion
of the graph of $X_t$ between 
$s_1 + \delta_1$ and $ s_2 - \delta_1$
does not intersect $W(\delta_2)$. 
Let $s_0 = s_1 + \delta_1$. Since the $Y^n_t$
converge to $X_t$, there exists a sequence
$m_j$ such that $X^{1/m_j}_{s_0} \to X_{s_0}$.
For sufficiently large $j$, the point
$(s_0, X^{1/m_j}_{s_0})$ lies outside
$W(\delta_2)$ and we also have $1/m_j < \delta_2$. Then,
for $t$ in a neighborhood of $s_0$, the function $X^{1/m_j}_t$
must be given by
$X^{1/m_j}_t  =  \wt X_t^{s_0, X^{1/m_j}_{s_0}}$.
We will show that 
$X_t  = \wt X_t^{s_0, X_{s_0}}$ for 
$t \in (s_0, s_2 - \delta_1)$.

Suppose that this is not true and let
$s_3 = \inf\{t\in[ s_0, s_2 -\delta_1] : X_t  \ne \wt X_t^{s_0, X_{s_0}} \}$.
Since $(s_3, X_{s_3})$ lies outside $W(\delta_2)$,
an argument similar to the one given above shows that
for some $\delta_3,\delta_4>0$, and all $m > 1/\delta_2$, the 
functions $X^{1/m}_t$ must satisfy
$X^{1/m}_t  =  \wt X_t^{s_3, X^{1/m}_{s_3}}$
for $t \in[s_3, s_3 + \delta_3]$, 
if $| X^{1/m}_{s_3}- X_{s_3}| \leq \delta_4$.
A straightforward argument now implies that for large $n$,
$Y^n_t  = \wt X_t^{s_3, Y^n_{s_3}}$
for $t \in[s_3, s_3 + \delta_3]$, and this in turn
proves that
$X_t  = \wt X_t^{s_3, X_{s_3}}$
for $t \in[s_3, s_3 + \delta_3]$. This contradicts
the definition of $s_3$ and proves our claim.

Thus $X_t$ satisfies (2.3) on
$(s_1 + \delta_1, s_2 - \delta_1)$ and, in view
of arbitrary nature of $\delta_1$, the same claim extends
to the whole interval $(s_1,s_2)$. The argument
applies to all intervals 
$(s_1,s_2)$ such that $B_t \ne f(X_t)$ for all
$t \in (s_1,s_2)$.
This implies that $X_t$ is
a Lipschitz solution to (2.3). The proof
of the existence of a Lipschitz solution is complete.

The existence of the solution to (2.3)
for $t< t_0$ may be proved in a completely
analogous way. The two solutions can
be combined into one function $X_t$ in an obvious
way. It remains to check if the differential equation
(2.3) is satisfied at $t = t_0$. It is easy to see
that if $B_{t_0} < x_0$ then 
$d X_t /dt = F_1(X_t)$ for all $t$ in some 
intervals $(t_0 - \delta, t_0)$ and 
$(t_0, t_0 + \delta)$ with $\delta >0$. This 
and the continuity of $X_t$ at $t=t_0$ evidently imply
that $d X_t /dt = F_1(X_t)$ for $t = t_0$
and so (2.3) is satisfied for $t=t_0$.
The case when $B_{t_0} > x_0$ is analogous.
When $B_{t_0} = x_0$ then (2.3) is trivially satisfied by $X_t$
for $t=t_0$.

Since the functions $\{X^{1/m}_t, t\geq t_0\}$ 
are adapted to the Brownian filtration
${\cal F}^B_t= \sigma(B_s, s\in[t_0, t])$,
so is their lim sup, $X_t$.
It follows that the process
$\{(B_t, X_t), t\geq t_0\}$ is strong Markov
with respect to the filtration
$\{ {\cal F}^B_t, t\geq t_0 \}$. 

We will show that the function
$\{X_t, t\geq t_0\}$ constructed above
is the largest of all Lipschitz solutions to (2.3),
that is, if $X^*_t$ is another Lipschitz solution, then
$X_t \geq X^*_t$ for all $t\geq t_0$.
Consider any Lipschitz solution $X^*_t$ to (2.3)
and suppose that $X^*_t > X_t$ for some $t\geq t_0$.
Then there must exist $\delta = 1/{m_j}$ such that
$X^*_t > X^\delta_t$ for some $t\geq t_0$. Fix such $\delta$
and let $S$ be the infimum of those $t$ such that
$X^*_t > X^\delta_t$. If $S \in [t_j+\delta, t_{j+1})$
for some $j$, then $X^*_S = X^\delta_S \ne B_S$ a.s.,
and, by continuity, we must have
$X^*_s \ne B_s$ and $ X^\delta_s \ne B_s$
for all $s$ in some non-degenerate interval $[S, S+\delta_1)$.
On this interval one of the conditions in (2.3)
is satisfied by  both $X^*_t$ and $X^\delta_t$,
so $X^*_s = X^\delta_s = \wt X^{S, X^*_S}$ for all $s \in [S, S+\delta_1)$
or $X^*_s = X^\delta_s = \wh X^{S, X^*_S}$ for all $s \in [S, S+\delta_1)$.
This contradicts 
the definition of $S$. Next suppose that
$S \in [t_j, t_j +\delta)$ for some $j$.
On this interval, the derivative of $X^\delta_t$
is equal to $\beta$. It is 
is easy to see that a Lipschitz solution $X^*_t$ to (2.3)
cannot grow faster than that on this interval,
and so $S \geq t_j +\delta$, a contradiction
which completes the proof of our claim.

A similar construction gives
a solution $\{X_t, t\leq t_0\}$ to (2.3) which is 
maximal among all Lipschitz solutions 
on the interval $(-\infty, t_0]$ with constant $\beta$.
Note that $X_t$ is measurable with respect to the $\sigma$-field
$\sigma(B_s, s\in[t,t_0])$ for $t < t_0$.

The maximal solution $X_t$ of (2.3) is consistent in the following sense.
Consider a fixed path $\{B_t, t\in \R\}$ and the corresponding 
maximal solution
$X_t$. Now choose any $s >0$ and suppose that
$X_s = z$. Let $\{X^*_u, u \geq s\}$ be the largest
Lipschitz solution with constant $\beta$
for the equation (2.3) on the interval
$[s , \infty)$ with the initial condition $X^*_s = z$ and the path
$\{B_t, t\in \R\} $ truncated to $\{B_t, t\geq s\} $. 
Then it is easy to see
that $X^*_u = X_u$ for all $u \geq s$.
It follows that for $s \geq 0$,
the portion $\{X_t, t \in [s,u]\}$ of the solution to
(2.3) may be defined only in terms of $X_s$
and $\{B_t, t\in[s,u]\}$. 

In a similar fashion we can construct a minimal solution  to (2.3);
this minimal solution is also adapted to the filtration of $B_t$.
Uniqueness would follow once we prove the maximal and minimal 
solutions are equal for all $s$ a.s.
\qed

\longproof{of Theorem 2.3}
Let $X^+$ and $X_-$ be the maximal and minimal solutions to (1.1). By (2.3)
the
$\P^x$ law of $B_t-X^-(t)$ is mutually absolutely continuous with
respect to the $\P^x$ law of $B_t$, so under $\P^x$, $B_t-X^-(t)$
has a jointly continuous local time $\wt L^x_t$
such that $\sup_x \wt L^x_t<\infty$, a.s. for each $t$.

Let
$$U(1)=\inf\{t>0: \sup_x \wt L^x_t\geq 1/(4\beta)\}.$$
If $t\leq U(1)$ and $a>0$, then
$$\int_0^t 1_{(B_s-X^-(s)\in [0,a])}ds=\int_0^a \wt L^x_t\, dx\leq a/(4\beta).$$
Let $a>0$ and
$$S=\inf\{t>0: X^+(t)-X^-(t)\geq a\}.$$
Since both $X^+$ and $X^-$ satisfy (1.1),
if $V=U(1)\land S$,
$$\eqalign{X^+(V)-X^-(V)&\leq 2\beta\int_0^V
1_{(X^-(u)\leq B_u\leq X^+(u))}du\cr
&\leq 2\beta \int_0^t 1_{(0\leq B_u-X^-(u)\leq X^+(u)-X^-(u))}du\cr
&\leq 2\beta \int_0^t 1_{(0\leq B_u-X^-(u)\leq a)}du\cr
&\leq 2a\beta/(4\beta)=a/2.\cr}$$
Since $X^+(V)-X^-(V)=a$ if $U(1)>S$, we must have $V=U(1)$.
This is true for all $a>0$, so $X^+(t)=X^-(t)$ for $t\leq U(1)$.

Now assume that $B_t$ is strong Markov 
and let $U(j+1)=U(j)+U(1)\circ \theta_{U(j)}$, $j=1,2, \ldots$,
where $\theta$ is the shift operator associated with the process
$B_t$. 
An induction argument using the strong Markov property at $U(j)$
shows that $X^+(t)=X^-(t)$ for $t\leq U(j+1)$
for $j=1,2,\ldots$. The continuity of $B_t$ and $ \wt L^x_t$ easily
implies $U(j)\to \infty$, a.s., so $X^+(t)=X^-(t)$ for all $t\geq t_0$.
\qed

\bs
The proof of Theorem 2.5 will be split into several
lemmas.

For the remainder of the section, 
let $\delta=(1-\gamma)/4$. Note that $\delta\in(0,1/4) $
since $\gamma \in (0,1)$. The constants $c_1, c_2, \dots$,
in the proofs in this section may depend on $\gamma$ and $\delta$.

\proclaim Lemma 2.14. Let $\al\geq 1$, $t\leq 1$, 
$A>0$, 
$C_t=\int_0^t 1_{(0<B_s<As^\al)}\, ds$. 
Assume that condition (i) of Theorem 2.5 holds.
There exist $c_1$ and $c_2$ independent
of $\al$ and $A$ such that for $\lambda >0$,
$$\P(C_t>\lam)\leq c_1\exp(-c_2\lam\al^\delta/(At^{\al+2\delta})).$$

\proof First let us compute $\E(C_t-C_u\mid \F_u)$ for $u\in [0,t]$.
Let $R=R(\al)=\al^{1/\al}$. 
Note that
$R\geq 1$, $R=\exp(\al^{-1}\log \al)\leq c_3$, 
and
$$1-R^{-1}=1-\exp(-\log \al/\al)\leq \log \al/\al\leq c_4\al^{-\half},$$
where $c_3$ and $c_4$  do not depend on $\al$ as long as $\al\geq 1$.

By condition (i) of Theorem 2.5,
$$\eqalign{\E(C_t-C_u\mid \F_u)&=\int_u^t \P(B_s\in (0,As^\al)\mid \F_u)\, ds\cr
&\leq \int_u^t \frac{c_5 As^\al}{(s-u)^{\gamma}}\, ds.\cr}$$

Let us examine
$$I=\int_u^t \frac{s^\al}{(s-u)^{\gamma}}\, ds.$$ Suppose first that
$u<t/R$. We observe, using the fact that $R\geq 1$, 
$$\eqalignno{\int_u^{t/R} \frac{s^\al}{(s-u)^{\gamma}}\, ds
&\leq \Big(\frac{t}{R}\Big)^\al\int_u^{t/R}
\frac{ds}{(s-u)^{\gamma}}\leq \frac{t^\al}{\al}\int_u^t \frac{ds}
{(s-u)^{\gamma}}\cr
&=\frac{t^\al}{\al}\int_0^{t-u}\frac{ds}{s^{\gamma}}\leq
c_6\frac{t^\al}{\al}t^{1-\gamma}. &(2.9)\cr}$$
On the other hand, in view of the inequality $1- R^{-1} \leq c_4 \alpha^{-1/2}$,
$$\eqalignno{\int_{t/R}^t \frac{s^\al}{(s-u)^{\gamma}}\, ds&\leq t^\al\int_{t/R}^t
\frac{ds}{(s-u)^{\gamma}}\cr
&\leq t^\al\int_{t/R}^t \frac{ds}{(s-t/R)^{\gamma}}\cr
&=t^\al \int_0^{t(1-1/R)}\frac{ds}{s^{\gamma}}\cr
&=c_7t^\al t^{1-\gamma}(1-R^{-1})^{1-\gamma}\cr
&\leq c_8 t^{\al+1-\gamma}/\al^{(1-\gamma)/2}. \cr}$$
Recalling that $\al\geq 1$ and combining with (2.9),
$$I\leq \frac{c_9t^{\al-1-\gamma}}{\al^{(1-\gamma)/2}}.$$

Now suppose $u\geq t/R$. Then
$$\eqalign{\int_u^t \frac{s^\al}{(s-u)^{\gamma}}\, ds&\leq t^\al\int_u^t \frac{ds}
{(s-u)^{\gamma}}=t^\al\int_0^{t-u}\frac{ds}{s^{\gamma}}\cr
&=c_9t^\al (t-u)^{1-\gamma}\leq c_9 t^\al(t-t/R)^{1-\gamma}\cr
&=c_9t^{\al+1-\gamma}(1-R^{-1})^{1-\gamma}.\cr}$$
As before, this is less than or equal to
$c_{10}t^{\al+1-\gamma}/\al^{(1-\gamma)/2}$.

Since $\delta = (1-\gamma)/4$, $t\leq 1$ and $\al\geq 1$,
$$\E(C_t-C_u\mid \F_u)\leq c_{11}At^{\al+4\delta}/\al^{2\delta}
\leq c_{11}At^{\al+2\delta}/\al^\delta.$$
This says that almost surely the process $\E(C_t\mid \F_u)$
does not exceed $C_u$ by more than $c_{11}At^{\al+2\delta}/\al^\delta$
for any $u\leq t$. In particular,
$$\E(C_t-C_T\mid \F_T)
\leq c_{11}At^{\al+2\delta}/\al^\delta$$
for every stopping time $T$ bounded by $t$.
We apply Theorem I.6.11 of Bass (1995)
to deduce that there exist $c_{12}$ and $c_{13}$ such that
$$\E\exp(c_{12}C_t\al^\delta/(At^{\al+2\delta}))\leq c_{13}.$$
Our result easily follows from this estimate.
\qed

\proclaim Lemma 2.15. 
Given $\xi >0$, 
there exist $c_1, c_2$ such that if $\al\geq 1$,
$A,B>0$, $B/A> \xi$, and $\beta=\al+\delta$, then
$$\P(C_t\geq Bt^\beta\hbox{ for some }t\leq \half)\leq c_1\exp(-c_2B\al^\delta/A).$$
 
\proof Let $t_k=2^{-1-k/\beta}$, $k=0, 1, \ldots$.
The process $C_t$ is increasing. 
So if $C_t\geq Bt^\beta$ for some $t\leq \half$,
then for some $k\geq 1$ we must have $C_{t_{k-1}}\geq B(t_{k})^\beta$.
Hence
$$\eqalignno{\P(C_t\geq Bt^\beta\hbox{ for some }t\leq \half)&\leq \P(C_{t_{k-1}}
\geq B(t_{k})^\beta\hbox{ for some }k\geq 1)\cr
&\leq \sum_{k=1}^\infty \P(C_{t_{k-1}}\geq B(t_k)^\beta). &(2.10)\cr}$$
Using Lemma 2.14, this is bounded by
$$\eqalignno{
\sum_{k=1}^\infty c_3&\exp\Big(-c_4 B(t_k)^\beta 
\al^\delta/ (A t_{k-1}^{\al +2 \delta}) \Big) \cr
&=\sum_{k=1}^\infty c_3\exp\Big(-c_4\frac{B\al^\delta}{A}
2^{-\beta - k - (-1-(k-1)/\beta) (\al +2 \delta)}\Big)\cr
&=\sum_{k=1}^\infty c_3\exp\Big(-c_4\frac{B\al^\delta}{A}
2^{ k \delta/\beta + \delta - (\al +2 \delta)/(\al + \delta)}\Big).}
$$
Since $\al \geq 1$ and $\delta\in(0,1/4)$, the quantity
$2^{\delta - (\al +2 \delta)/(\al + \delta)}$ is bounded below
and above by absolute constants, so the last displayed
formula admits a bound
$$\eqalignno{
\sum_{k=1}^\infty c_3&\exp\Big(-c_5\frac{B\al^\delta}{A}
2^{k\delta/\beta}\Big)&(2.11)\cr
&=c_3\exp\Big(-c_5\frac{B\al^\delta}{A}\Big)\sum_{k=1}^\infty
\exp\Big(-c_5\frac{B}{A}\al^\delta (2^{k\delta/\beta}-1)\Big).\cr}$$
The infinite sum in the last expression
is bounded by
$$\eqalign{\sum_{k=1}^\infty \exp(-c_5\frac{B}{A}(2^{k\delta/\beta}-1))&
\leq \sum_{k=1}^\infty \exp\Big(-\frac{c_5 B k\delta \log 2}{A\beta}\Big)\cr
&\leq \frac{1}{1-\exp(-c_5 B \delta \log 2/(A\beta))}\cr
&\leq c_6 A\beta/B .}$$
Combining this with (2.10) and (2.11) we obtain
$$\eqalign
{\P&(C_t\geq Bt^\beta\hbox{ for some }t\leq \half)
\leq c_3\exp\Big(-c_5\frac{B\al^\delta}{A}\Big)
c_6 A\beta/B \cr
&= c_3\exp\Big(-c_5\frac{B\al^\delta}{A}
+ \log c_6 + \log (A/B) + \log (\al + \delta) \Big)\cr
&\leq c_3\exp\Big(-c_5\frac{B\al^\delta}{A}
+ \log c_6 - \log \xi + \log 2 +\log \al  \Big).
}$$
The last expression is less than
$$c_7\exp\Big(-c_8\frac{B\al^\delta}{A}\Big)$$
for suitable $c_7$ and $c_8$ (depending on $\xi$ and $\delta$) and all $\al\geq 1$.\qed

Let $X_t^+$ and $X_t^-$ be the maximal and minimal solutions to  (2.3)
constructed in in the proof of Theorem 2.2. Let $Y_t=X_t^+-X_t^-$. We will show
$Y_t=0$, a.s. for $t\leq 1/2$. 

\proclaim Lemma 2.16. For each $s$,
$$\P(X_s^+=B_s)=0, \qq a.s.$$
and similarly with $X_s^+$ replaced by $X_s^-$.

\proof We know $X_s^+$ is a process whose paths are Lipschitz
continuous. 
By assumption (ii) of Theorem 2.5, there exists
a probability measure $\Q$ which is equivalent to $\P$ and
such that the $\Q$ law of $B_s-X_s^+$ is the same as the
$\P$ law of $B_s$. Then
$$\Q(X_s^+=B_s)=\Q(B_s-X_s^+=0)=\P(B_s=0).$$
This is equal to zero by (2.4). Since $\P$ and $\Q$
are equivalent, the lemma is proved. \qed

\proclaim Lemma 2.17. $Y_t=0$, a.s. if $t\leq 1/2$.

\proof The process $X_t^+$ satisfies the equation
$$X_t^+=x+\int_0^t [F_1(X_s^+)1_{(X_s^+>B_s)}
+F_2(X_s^+)1_{(X_s^+<B_s)}]\, ds.$$
$X_t^-$ satisfies a similar equation. Then, noting Lemma 2.16,
$$\eqalign{Y_t&=\int_0^t[F_1(X_s^+)-F_1(X_s^-)]1_{(B_s<X_s^-\leq X_s^+)}
\,ds\cr
&\qq +\int_0^t [F_2(X_s^+)-F_2(X_s^-)]1_{(X_s^-\leq X_s^+<B_s)}\,ds\cr
&\qq +\int_0^t [F_1(X_s^+)-F_2(X_s^-)]1_{(X_s^-<B_s<X_s^+)}\, ds.\cr}$$
Therefore
$$\eqalignno{Y_t&\leq M\int_0^t (X_s^+-X_s^-)\, ds+2M\int_0^t 1_{(X_s^-
<B_s<X_s^+)}\, ds&(2.12)\cr
&=M\int_0^t Y_s\, ds+2M\int_0^t 1_{(0<B_s-X_s^-<X_s^+-X_s^-)}\,ds\cr
&= M\int_0^t Y_s\, ds+2M\int_0^t 1_{(0<B_s-X_s^-<Y_s)}\, ds.\cr}$$

Recall that we have assumed that $F_j$ is bounded by $M$. Hence,
the process $Y_t$ is Lipschitz with constant $2M$. 
Since $X_s^-$ has Lipschitz paths, there exists, by assumption 
(ii) of Theorem 2.5,
a probability measure $\Q$ equivalent to $\P$ such that under $\Q$, 
$\{B_s-X_s^-, 0 \leq s \leq 1/2\}$ has the same law as 
$\{B_s, 0\leq s \leq 1/2\}$ does under $\P$. 
So it suffices to show that for any Lipschitz process $Y_s$ with constant $M$
satisfying
$$Y_t\leq M\int_0^t Y_s\, ds+2M\int_0^t 1_{(0<B_s<Y_s)}\, ds,\eqno (2.13)$$
we have 
$$ \P(Y_t\ne 0 \hbox{ for some } t\leq 1/2)=0.$$

Let
$$D(A,\al)=\{Y_s\geq As^\al \hbox{ for some } s\leq 1/2\}.$$
As $Y$ is Lipschitz with $|Y_t|\leq 2Mt$, then
$D(3M,1)=\emptyset$. Let $\eps>0$ and let $\eta=1/4$. 
We will choose $N \geq 1$ and $j_0\geq 0$
in a moment. Let $A_j=N^j$ if $j\leq j_0$ and $A_j=(1+\eta)^j N^{j_0}$
for $j>j_0$. Let $\al_j=1+j\delta$. We want an estimate on the probability of
$D(A_{j+1},\al_{j+1})-D(A_j,\al_j)$. If $\omega\notin D(A_j, \al_j)$, then
$Y_s\leq A_js^{\al_j}$ for all $s\leq 1/2$, and so from (2.13), for $t\leq 1/2$,
$$\eqalignno{Y_t&\leq M\int_0^t A_js^{\al_j}\,ds+
2M\int_0^t 1_{(0<B_s<A_js^{\al_j})}\, ds&(2.14)\cr
&=\frac{MA_jt^{\al_j+1}}{\al_j+1}+2M\int_0^t 1_{(0<B_s<A_js^{\al_j})}\, ds.\cr}$$

Let $\xi = (1-\eta)/2 M$ and let $c_1$ and $c_2$ be constants
chosen as in Lemma 2.15 (depending on  $\xi$).
Find large $j_0$ so that
$$(1+j_0\delta)^{\delta/2}/2M\geq 1,\eqno(2.15)$$
$$\frac{M}{1+j_0\delta}\leq \eta(1+\eta),\eqno(2.16)$$
and
$$c_1 \sum_{j=j_0}^\infty \exp(-c_2(1-\eta^2)(1+j\delta)^{\delta/2})<\eps/2.\eqno(2.17)$$
Next choose $N$ large so that
$$N\geq M/\eta\eqno(2.18)$$ 
and
$$2j_0 c_1 \exp(-c_2(1-\eta)N/2M)<\eps/2.\eqno(2.19)$$

For $j\geq j_0$, we have
$$\frac{MA_j}{(\al_j+1)}\leq \eta A_{j+1},\eqno (2.20)$$
using (2.16). The same inequality holds for $j< j_0$
in view of (2.18).

In view of (2.14) and (2.20),
for $\omega$ to be in $D(A_{j+1},\al_{j+1})-D(A_j,\al_j)$, we must have,
$$\eqalignno{
\int_0^t 1_{(0<B_s<A_js^{\al_j})}\, ds&\geq 
Y_t/(2M) - \frac{A_jt^{\al_j+1}}{2(\al_j+1)} &(2.21)\cr
& \geq A_{j+1}t^{\al_{j+1}}/(2M) - \frac{A_jt^{\al_j+1}}{2(\al_j+1)} \cr
&\geq (1-\eta)A_{j+1}t^{\al_{j+1}}/(2M)
}$$
for some $t<1/2$. 
Recall that we set $\xi = (1-\eta)/2 M$
and note that for all $j$ we have $(1-\eta)A_{j+1}/(2MA_j) \geq \xi$.
By Lemma 2.15, the probability that the inequality (2.21) 
holds is less than or equal to
$$c_1\exp\Big(-c_2\frac{(1-\eta)A_{j+1}}{2MA_j}\al_j^\delta\Big).$$
Using (2.15) and (2.17) for $j\geq j_0$, we obtain
$$c_1\sum_{j=j_0}^\infty \exp\Big(-c_2\frac{(1-\eta)A_{j+1}}{2MA_j}\al_j^\delta
\Big)\leq 
c_1\sum_{j=j_0}^\infty \exp(-c_2\frac{1-\eta^2}{2M}(1+j_0\delta)^{\delta/2}
(1+j\delta)^{\delta/2}\Big)<\eps/2.$$
{}From (2.19), 
$$\eqalign{
c_1\sum_{j=0}^{j_0-1} \exp\Big(-c_2\frac{(1-\eta)A_{j+1}}{2MA_j}\al_j^\delta
\Big)&\leq 
c_1\sum_{j=0}^{j_0-1} \exp\Big(-c_2\frac{(1-\eta)N}{2M}
\Big)\cr
&\leq 2j_0 c_1 \exp(-c_2(1-\eta)N/2M)<\eps/2.}$$
Hence,
$$c_1\sum_{j=0}^\infty \exp\Big(-c_2\frac{(1-\eta)A_{j+1}}{2MA_j}\al_j^\delta
\Big)\leq \eps,$$
and so
$$\P\Big(\bigcup_{j=0}^\infty D(A_j,\al_j)\Big)\leq \eps.$$

If $\omega\notin \bigcup_{j=0}^\infty D(A_j,\al_j)$, then
$Y_t(\omega)\leq A_jt^{\al_j}\leq (1+\eta)^jN^{j_0}(1/2)^{1+j\delta}$
for all $j\geq j_0$ and all $t\leq 1/2$. Since $(1+\eta)(1/2)<1$, letting $j\to \infty$ shows
$Y_t(\omega)=0$.
Therefore
$$\P(Y_t\ne 0\hbox{ for some } t\leq 1/2)\leq \eps.$$
Since $\eps$ is arbitrary, this proves the lemma. \qed

\longproof{of Theorem 2.5} By Lemma 2.17 we have
$Y_t=0$ a.s. for $t\leq 1/2$.
If we consider the law of $B_{t+1/2}$ given $\F_{1/2}$, 
it is not hard to see that assumptions (i) and (ii) of Theorem 2.5
apply to 
this process as well. So we apply the same argument to
$X_{t+1/2}^+$ and $X_{t+1/2}^-$, and we obtain 
$Y_{t+1/2}=0$ for $t\leq 1/2$, or $Y_t=0$ for $t\leq 2(1/2)$.
By an induction argument, we then have $Y_t=0$ for all $t$,
which proves uniqueness.
\qed

\bigskip
\noindent{\bf Proof of Theorem 2.6}.
The existence and strong uniqueness of solutions
$X^\eps_t$ to (2.6) can be proved
in the same way as in Theorem 2.1.

Consider any sequence $\eps_n \downarrow 0$
and with a slight abuse of notation
let $X^n_t = X^{\eps_n}_t$. 
Since all functions $t \to X^n_t$ are Lipschitz
with constant $\beta$, we may suppose,
passing to a subsequence, if necessary, that $X^n_t $
converge to a function $X^\infty_t$.
In order to finish the proof, it will suffice to show
that $X^\infty_t = X_t$. Since the equation (2.5)
has a unique solution a.s.,  it will be enough to show that if
$\omega$ is not in the null set where uniqueness does not hold, 
then $X^\infty_t(\omega)$ is a solution to (2.5).
The functions $X^n_t$ are Lipschitz with constant $\beta$,
so the same is true of $X^\infty_t$. Let $A$ be the set of times $t$
such that $X^\infty_t = B_t$.
The complement of the set $A$ consists of a countable
number of open intervals. Let $I = (t_1,t_2)$
be one of the intervals in the complement of $A$.
Fix any $t_3 \in I$ and suppose without loss of
generality that $X^\infty_{t_3} > B_{t_3}$.
Choose some $t_4 \in (t_1, t_3)$ and
$t_5 \in (t_3, t_2)$ and let $a$ 
be the infimum of $X^\infty_t - B_t$
over $t \in (t_4,t_5)$. For sufficiently
large $n$, we have $\eps_n < a/3$ and
$|X^n_t - X^\infty_t| < a/3$ for all $t \in (t_4,t_5)$.
It follows that for large $n$ and $t \in (t_4,t_5)$, we have
$X^n_t - B_t  > a/3 > \eps_n$. Hence, for
such $n$ and $t$, $dX^n_t/dt = \beta_2$. This
shows that $d X^\infty_t /dt = \beta_2$
for all $t \in I$. The same argument works for
all other intervals in the complement of $A$.
There is nothing to check for $t\in A$, so $X_t^\infty$ is  a solution to (2.5).
\qed

\bigskip
\noindent{\bf 3. Local time}. In the remaining part
of the article we assume that $B_t$ is a Brownian motion.
In this section we will exclusively deal with
solutions to (1.1).
We will find several explicit 
formulae for the local time spent by $B_t$
on the paths of the process $X_t$.
Moreover, we will prove analogues of
the Trotter and Ray-Knight theorems.
The results on local times provide information
about the behavior of the function $y\to X^y_t$,
for  fixed $t$; see Remark 3.9.

The first part of the section deals with exit
systems. 
Some of our results on exit systems may be of independent
interest.
We refer the  reader to Blumenthal (1992),
Burdzy (1987), Maisonneuve (1975) or Sharpe (1989)
concerning the fundamentals of excursion theory.
\ms

In this section, we will assume that $t_0=0$
and study the portion of the solution $X_t$ to (1.1)
for $t\geq 0$ only.

Let $D= \{(b,x) \in \R^2: b = x \}$.
We will construct
an exit system $(H^x, dL)$ for the process of
excursions of $(B_t,X_t)$ from the set $D$.
The first element of an exit system is a family of
excursion laws $H^x$. An excursion law
$H^x$ is an infinite $\sigma$-finite
measure on the space $C^*$ of functions 
$(e^1_t, e^2_t)$ defined on $(0,\infty)$ (note that $0$ is excluded)
which take values in $\R^2 \cup \{\Delta\}$. Here 
$\Delta$ is the coffin (absorbing) state.
Let $\nu$ be the lifetime of an excursion, i.e., 
$\nu=\inf \{t > 0 : (e^1_t, e^2_t) = \Delta\}$. 
Then $H^x$-a.e., we have $(e^1_t, e^2_t) \in \R^2$ for
$t\in(0,\nu)$ and $(e^1_t, e^2_t)=\Delta$
for $t \in [\nu,\infty)$.
The measure $H^x$ is strong Markov with respect
to the transition probabilities of the process
$\{(B_t, X_t), t\geq 0\}$ killed at the hitting
time of $D$. Moreover, the $H^x$-measure
of the set of paths for which
$\lim _{t\downarrow 0} (e^1_t, e^2_t) \ne (x,x)$
is equal to $0$.
The second element of the exit system, $dL$, denotes
the measure defined by a non-decreasing process $L_t$.
The process $L_t$ is a continuous additive functional,
also known as a local time, for
$(B_t,X_t)$ on $D$. The process $L_t$ does not
increase on any interval $(s,u)$ such that
$(B_t,X_t) \notin D$ for $t\in (s,u)$; that is,
$L_s = L_u$ for such intervals.
Consider a maximal interval $(s,u)$ such that
$B_t \ne X_t$ for $t \in (s,u)$. Suppose
$L_s = r$. Let
$(e^1_t, e^2_t)_r = (B_{s+t}, X_{s+t})$ for
$t \in (0, u-s)$ and
$(e^1_t, e^2_t)_r = \Delta$ for $t \geq u-s$.
Let $\mu(r) = \inf \{ t>0: L_t = r\}$.
The collection of all ``excursions''
$\{(r, (e^1_\cdot, e^2_\cdot)_r)\}$ 
is a Poisson point process which, roughly speaking,
has  random mean measure
$(r_2-r_1)\int _{r_1}^{r_2} H^{\mu(r)}(A) dr$
on the set $(r_1,r_2) \times A$.

Next we apply some transformations to
the excursions and excursion laws in order
to simplify our description of the exit system.
First, we note that by the translation invariance
of the Brownian motion $B_t$ and equation (1.1),
the distribution of $(e^1_t - x, e^2_t -x)$ under $H^x$
is the same for every $x \in \R$. Let this distribution
be called $H_1$. For $H_1$-almost all excursions,
the second component $e^2_t$ is a linear function
of $t$ until the excursion lifetime $\nu$, with the slope
equal to $\beta_1$ or $\beta_2$. In the first case,
$e^1_t > e^2_t$ for $t \in (0,\nu)$, while the inequality
goes the other way in the second case. Let $H_{1+}$
denote the part of the measure $H_1$ which is supported
on excursions with $e^1_t > e^2_t$ and let $H_{1-}$
be the part supported on the set where $e^1_t < e^2_t$.
Let $H_{2+}$ be the distribution of 
$\{e^1_t - e^2_t, t \in (0,\nu)\}$ under $H_{1+}$
and let $H_{2-}$ have the same definition relative
to $H_{1-}$. Note that, by definition,
the excursion laws $H_{2+}$ and $H_{2-}$ are
supported on paths in $\R \cup \{\Delta\}$ 
rather than $\R^2 \cup \{\Delta\}$, since
the second component becomes irrelevant after our last
transformation.

Our transformations preserve the strong Markov
property, but the last transformation creates
a drift so that the measure $H_{2+}$ has the transition
probabilities of Brownian motion with drift $-\beta_1$,
killed upon hitting $0$. It is standard to show
(see, e.g., Theorem 4.1 of Burdzy (1987))
that for any event $A$ defined in terms of
the process after some fixed time $s_0>0$, we have,
up to a multiplicative constant,
$$H_{2+}(A) = \lim_{x \downarrow 0}
{1 \over |x|} Q^x_{-\beta_1}(A),\eqno(3.1)$$
where $Q^x_{-\beta_1}$ stands for the distribution
of Brownian motion with drift $-\beta_1$, killed
at the hitting time of $0$. The normalization
of the excursion laws is arbitrary as long as
it matches the normalization of the local time,
so we can use the normalization
in (3.1). 
We next choose the normalization of the local time
so that it matches that of $H_{2+}$.
Given the normalization for $H_{2+}$, the
normalization for $H_{2-}$ is no longer arbitrary
and we will have to prove that
$$H_{2-}(A) = \lim_{x \uparrow 0}
{1 \over |x|} Q^x_{-\beta_2}(A).\eqno(3.2)$$

Unless specified otherwise, all excursion laws
in this paper will be normalized
as in (3.1) or (3.2).

Let $H_3$ denote the excursion law
for excursions of Brownian motion without
drift away from $0$. Let us split $H_3$ into
positive and negative parts $H_{3+}$ and $H_{3-}$,
as in the case of $H_2$. We normalize $H_3$
using a formula analogous to (3.1).
Recall that $\nu$ denotes the lifetime of an excursion $e$,
and that (3.1) defines the normalization of $H_{2+}$.

\bigskip
\noindent{\bf Lemma 3.1}. {\sl
(i) On the set where $\nu < \infty$,
$${dH_{2+} \over d H_{3+}} (e) 
= \exp(  - \beta_1^2 \nu / (2 \sigma^2)). $$

(ii) For a fixed time $s\in(0,\infty)$,
the conditional distributions of
$H_{2+}$ and $H_{3+}$ given $\{\nu=s\}$
are identical.

(iii) If $\beta_1 <0$ then
$H_{2+} ( \nu = \infty) = 2 |\beta_1| / \sigma^2$.

(iv)  Formula (3.2) is the correct normalization
for $H_{2-}$.
}

\bigskip
Parts (i)-(iii) of Lemma 3.1 have obvious analogues
for $H_{2-}$.

\bigskip
\noindent{\bf Proof}.
Fix arbitrary $0 < s_0 < s_1 < \infty$ and let $A$ be an event
measurable with respect to $\sigma\{e_t, t \in (s_0, s_1)\}$.
Since $H_{3+}$ is assumed to be normalized using a formula
analogous to (3.1), we have
$${H_{2+} (A \cap \{\nu  = s_1\})
\over  H_{3+}(A \cap \{\nu  = s_1\})} 
= \lim_{x \downarrow 0}
{Q^x_{-\beta_1}(A \cap \{\nu  = s_1\}) 
\over Q^x_0(A \cap \{\nu  = s_1\}) } .
$$
An application of Girsanov's Theorem, as in 
Karatzas and Shreve (1988) ((5.11), p.~196),
shows that 
$$
{Q^x_{-\beta_1}(A \cap \{\nu  = s_1\}) 
\over Q^x_0(A \cap \{\nu  = s_1\}) }
= \exp(  x \beta_1 /\sigma^2 - \beta_1^2 s_1 / (2 \sigma^2)).
$$
This and the previous formula imply
$${H_{2+} (A \cap \{\nu  = s_1\})
\over  H_{3+}(A \cap \{\nu  = s_1\})} 
= \exp(  - \beta_1^2 s_1 / (2 \sigma^2)),
$$
which then easily implies (i) and (ii).

As for (iii), we start with the formula
$$Q^x_{-\beta_1} ( \nu = \infty)
= 1 - \exp( 2 x \beta_1 / \sigma^2),$$
with $\beta_1 <0$ (Karlin and Taylor (1975), p.~362).
Then (3.1) yields
$$H_{2+} ( \nu = \infty) 
= \lim_{x \downarrow 0}
{1 \over |x|} Q^x_{-\beta_1}( \nu = \infty) 
= 2 |\beta_1| / \sigma^2,$$
as desired.

It remains to prove (iv). 
Fix arbitrarily small $\gamma>0$ and let
$$A_1=A_1(t) = \{\max_{s\leq t} |X_s| > t ^{1/2+\gamma}\}.$$
Note that $|X_t| \leq \beta t < t ^{1/2+\gamma}$ for small $t>0$
so we have $\P (A_1(t)) = 0$ if $t$ is small. However, we will
prove the result using only the property that
$\lim_{t\to 0} \P (A_1(t)) = 0$ because we will need this
version of the proof later in the paper. 
Let us take $t_0=0$ and $x_0=0$
so that $X_0=0$. 
Note that the excursion law normalization does not
depend on $t_0$ and $x_0$.
Let $A_+= A_+(s)$ be the event that
the first excursion $(e^1_t, e^2_t)$ of
$(B_t, X_t)$ from $D$ with the property that
$|e^1_t - e^2_t|> s^{1/2+\gamma/2}$
for some $t \in (0,\nu)$, also has  the property
that $e^1_t > e^2_t$ for $t \in (0,\nu)$.
Let $A_- $ be the analogous event with $e^1_t < e^2_t$.
Let $T(a)$ be the hitting time of $a$ by $B_t$.
For small $s>0$,
$$\{ T( s^{1/2+\gamma/2} + s^{1/2+\gamma})
< T(- s^{1/2+\gamma/2} + s^{1/2+\gamma}) < s \} \subset A_+(s) \cup A_1(s),$$
and
$$\{ T( -s^{1/2+\gamma/2} - s^{1/2+\gamma})
< T( s^{1/2+\gamma/2} - s^{1/2+\gamma}) < s \} \subset A_-(s) \cup A_1(s).$$
It is elementary to check that
$$\eqalign{ \lim_{s \to 0} 
\P &( T( s^{1/2+\gamma/2} + s^{1/2+\gamma})
< T(- s^{1/2+\gamma/2} + s^{1/2+\gamma}) < s ) \cr
&= \lim_{s \to 0} 
\P ( T( -s^{1/2+\gamma/2} - s^{1/2+\gamma})
< T( s^{1/2+\gamma/2} - s^{1/2+\gamma}) < s ) = 1/2.}$$
This and the fact that $\lim_{t\to 0} \P (A_1(t)) = 0$ imply that
$$\lim_{s\to 0} \P (A_+(s)) =
\lim_{s\to 0} \P (A_-(s)) = 1/2.\eqno(3.3)$$

The scale function $S(y)$ for Brownian motion 
with drift $-\beta_1$ is given by \break
$S(y) = \exp( 2 \beta_1 y/\sigma^2)$
(Karlin and Taylor (1981) Chapter 15.4).
Let $F_h$ be the event that the difference between the maximum 
and the minimum
of an excursion exceeds $h$. Then, by (3.1),
$$\eqalign{H_{2+}(F_h) &= 
\lim_{x \downarrow 0}
{1\over x} Q^x_{-\beta_1} (T_h < T_0)
= \lim_{x \downarrow 0}
{1\over x} \cdot{ S(x) - S(0) \over S(h) - S(0) } \cr
&= \lim_{x \downarrow 0}
{1\over x} \cdot{  \exp( 2 \beta_1 x /\sigma^2) -1 \over
\exp( 2 \beta_1 h/\sigma^2) -1 } 
= {2 \beta_1 \over \sigma^2} \cdot
{ 1 \over \exp( 2 \beta_1 h/\sigma^2) -1 }.
}$$
An analogous formula holds for $H_{2-}(F_h)$,
but we will write it with an additional multiplicative
constant $c_1$, since we have not proved that (3.2)
is the right normalization yet:
$$H_{2-}(F_h) = c_1
{2 \beta_2 \over \sigma^2} \cdot
{ 1 \over \exp( 2 \beta_2 h/\sigma^2) -1 }.$$
Our goal is to show that $c_1 =1$ is the correct choice
for the constant.

Excursion theory tells us that the 
arrival times for excursions
$(e^1_t, e^2_t)$ of
$(B_t, X_t)$ from $D$ with the property that
$|e^1_t - e^2_t|> s^{1/2+\gamma/2}$
for some $t \in (0,\nu)$, and with 
$e^1_t > e^2_t$ for $t \in (0,\nu)$, form
a Poisson point process on the local time scale
with intensity $H_{2+}(F_{s^{1/2+\gamma/2}})$.
This process is independent from the analogous
process of excursions with $e^1_t < e^2_t$.
Formula (3.3) tells us that for small $s$,
the probability that the first arrival for the
first process is earlier than the first arrival for the
second process is close to $1/2$.
Hence, the ratio of the
intensities for the two Poisson point processes must
converge to 1 as $s \to 0$. Therefore, we must have
$$\lim_{s \to 0}
{2 \beta_1 \over \sigma^2} \cdot
{ 1 \over \exp( 2 \beta_1 s^{1/2+\gamma/2}/\sigma^2) -1 }
\cdot \left(
c_1 {2 \beta_2 \over \sigma^2} \cdot
{ 1 \over \exp( 2 \beta_2 s^{1/2+\gamma/2}/\sigma^2) -1 } \right)^{-1}
 =1.$$
However, this is possible only if $c_1 =1$. This completes the proof of (iv).
\qed

\bigskip
\noindent{\bf Remark 3.2}. (i) Lemma 3.1 (iv) 
can be used to prove uniqueness
for (1.1). In order to do so, one would have to 
consider an exit system for the maximal Lipschitz solution $X_t$
to (1.1), constructed as in the proof of Theorem 2.2,
and the analogous exit system for the minimal
Lipschitz solution.
Lemma 3.1 (iv) shows that both 
exit systems are identical but this can be true only if the
maximal and minimal solutions are the same.
We will not formalize this argument as it cannot be easily
generalized to non-Markov processes. The delicate
part of the argument would be to show that the maximal solution $X_t$ is the sum
of the excursions and it contains no component corresponding
to  a ``push'' proportional to local time.

(ii) According to Lemma 3.1 (i)-(iii), if $\beta_1 <0$,
the excursion laws $H_{2+}^{\beta_1}$ and $H_{2+}^{-\beta_1}$
agree on the set of excursions with finite lifetime
and the only difference is that $H_{2+}^{\beta_1}$ gives some
mass to excursions with infinite lifetime, while
$H_{2+}^{-\beta_1}$ does not.

\bigskip

For every $x \in \R$ consider the solution
$X^x_t$ to (1.1) with $X^x_0 =x$.
Let $L^x_t$ denote the local time of 
$Y_t^x \df B_t - X^x_t$ at $0$, defined
earlier in this section as the local
time of $(B_t, X^x_t)$ on the diagonal,
accumulated between times $0$ and $t$.
Note that this is not the local time of a one-dimensional
diffusion at level $x$.

\bigskip
\noindent{\bf Proposition 3.3} {\sl
If $\beta_1 > 0$ and $\beta_2 < 0$ then
for every $x \in \R$,
$$\lim_{t\to \infty} L_t^x/t = 
\left({1\over |\beta_1|} + {1\over |\beta_2|}\right)^{-1},\qquad 
\hbox{a.s.} $$
}

\bigskip
\noindent{\bf Proof}.
Fix some $x\in\R$.
Our assumptions that $\beta_1 > 0$ and $\beta_2 < 0$
imply that there will never be an excursion
of $B_t$ from $X^x_t$ with infinite
lifetime, since the drift will always push the excursions
of $Y^x_t$ towards $0$. This in turn implies
that $L^x_t$ will grow to infinity a.s.

Recall that we used $Q^y_{-\beta_1}$ to denote
the distribution of Brownian motion with drift
$-\beta_1$, killed at the hitting time of $0$.
By Theorem 7.5.3 of Karlin and Taylor (1975), we have
for $y>0$,
$$Q^y_{-\beta_1}(\nu\in dt) = {|y| \over \sigma t^{3/2} \sqrt{2 \pi} }
\exp\left( - { (|y| - |\beta_1| t)^2 \over 2 \sigma^2 t} \right) dt,$$
where $\nu$ denotes the lifetime of the process.
The same formula holds for $y<0$, with $\beta_1$
replaced by $\beta_2$.
Using (3.1)-(3.2) we obtain
$$\eqalign{
H_2&(\nu\in dt)  =
\lim_{y \downarrow 0} {1 \over |y|}Q^y_{-\beta_1}(\nu\in dt)
+ \lim_{y \uparrow 0} {1 \over |y|}Q^y_{-\beta_2}(\nu\in dt) \cr
&= \lim_{y \downarrow 0}
{1 \over \sigma t^{3/2} \sqrt{2 \pi} }
\exp\left( - { (|y| - |\beta_1| t)^2 \over 2 \sigma^2 t} \right)dt
+ \lim_{y \uparrow 0}
{1 \over \sigma t^{3/2} \sqrt{2 \pi} }
\exp\left( - { (|y| - |\beta_2| t)^2 \over 2 \sigma^2 t} \right)dt \cr
&= {1 \over \sigma t^{3/2} \sqrt{2 \pi} }
\exp[- (\beta_1 ^2 / 2 \sigma^2) t ]dt
+{1 \over \sigma t^{3/2} \sqrt{2 \pi} }
\exp[- (\beta_2 ^2 / 2 \sigma^2) t ]dt.}
$$
Let $V^x_s$ be the inverse local time,
i.e., $V^x_s = \inf\{t>0: L^x_t > s\}$.
The process $V^x_s$ is the sum of lifetimes
of excursions which start before the local time
reaches the level $s$. The Poisson character of the excursion
process easily implies that
$$\eqalignno{
\E  V^x_s &= s \int_0^\infty t H_2(\nu\in dt) \cr
&= s \int_0^\infty t {1 \over \sigma t^{3/2} \sqrt{2 \pi} }
\exp[- (\beta_1 ^2 / 2 \sigma^2) t ] dt 
+ s \int_0^\infty t {1 \over \sigma t^{3/2} \sqrt{2 \pi} }
\exp[- (\beta_2 ^2 / 2 \sigma^2) t ] dt \cr
&= \left( {1 \over |\beta_1|} + {1 \over |\beta_2|}\right)s.& (3.4)}$$
This and the strong law of large numbers for the
L\'evy process $s\to V^x_s$
(see p.~92 of Bertoin (1996)) imply that 
$$V^x_s/s \to {1 \over |\beta_1|} + {1 \over |\beta_2|},$$
a.s., as $s\to \infty$. This can be easily translated
to the statement of the proposition.
\qed
\bigskip

We note that if $\beta_1,\beta_2>0$,
then we will eventually have
$X^x_t > B_t$, for every $x$. Hence, in this case,
$L^x_\infty < \infty$ for every $x\in \R$, a.s. 
We will prove the next lemma under the assumption
that $\beta_1 -\beta_2 >0$. We believe that
similar statements hold when $\beta_1 -\beta_2 <0$
but technical difficulties prevent us from giving a formal
proof in that case.

The following lemma contains the most complicated
and technical argument in the whole article. 

\bigskip
\noindent{\bf Lemma 3.4}. {\sl (i)
Fix $x,a,\beta_1,\beta_2\geq 0$ and assume that 
$\beta_1 -\beta_2 >0$. Then
$$\E (L^{x+\delta}_\infty \mid L^x_\infty = a)
= a - \delta{\beta_1 \over \beta_1 -\beta_2}
(1- \exp (-2 a (\beta_1 -\beta_2)  /\sigma^2) ) + o(\delta),$$
for $\delta \downarrow 0$. 

\noindent (ii) 
If $a,\beta_1,\beta_2\geq 0$, $x\leq 0$, and 
$\beta_1 -\beta_2 >0$ then
$$\E (L^{x+\delta}_\infty \mid L^x_\infty = a)
= a + \delta\left[
{\beta_2 \over \beta_1-\beta_2}
- {\beta_1 \over \beta_1-\beta_2}
\exp (-2 a (\beta_1 -\beta_2)  /\sigma^2) \right] + o(\delta),
$$
for $\delta \uparrow 0$.
}

\bigskip
\noindent{\bf Proof}.
(i) Recall that $B_0=0$, that we have fixed
$x,a,\beta_1,\beta_2\geq 0$ and assumed that 
$\beta_1 -\beta_2 >0$.
Since the proof of the lemma is quite long, we will
split it into several steps.

\medskip
\noindent{\it Step 1}. 
We start with some transformations of the
processes $X^x_t$ and $B_t$ which will enable us to 
look at $L^x_\infty$ from a slightly different perspective. It is perhaps not
necessary to make these transformations, but we find
the transformed problem much easier to 
comprehend than the original one
from an intuitive point of view.

We first offer a rough guide to our notation
(whose validity is limited to this proof).
Different Brownian motions with different drifts
and reflected barriers will be denoted $B^j_t$, for $j=1,2, \dots$.
The notation $L^{j-}_t$ and $L^{j+}_t$ will refer
to the local time of $B^j_t$ on the lower and upper reflected
barriers (if any). We will write $v_j(a) = \inf\{t: L^{j-}_t= a\}$.

It is well known that the set
$\{t: B_t = 0\}$ has zero Lebesgue measure, so the
Girsanov theorem implies that the same is true of
the set $\{t: X^x_t=B_t\}$. We will excise all
intervals where $X^x_t > B_t$. First, we define
a clock 
${\cal C}_1(t) = \int_0^t \bone_{\{X^x_s \leq B_s\}} ds$
and its inverse $b_1(t) = \inf\{s: {\cal C}_1(s)\geq t\}$.
Since $\beta_1,\beta_2>0$, we will eventually have
$X^x_t > B_t$, so we let $u_1 = \sup\{ {\cal C}_1(t): t\geq 0\}$.
Then we define new processes on the random interval
$[0,u_1]$ by
$$\eqalign{
B^1_t &= B_{b_1(t)} - \beta_2 (b_1(t) - t), \cr
X^{1,x}_t &= X^x_{b_1(t)} - \beta_2 (b_1(t) - t), \cr
X^{1,x+\delta}_t &= X^{x+\delta}_{b_1(t)} - \beta_2 (b_1(t) - t).
}$$
For $t \in [0,u_1]$,
we have $X^{1,x}_t = x + \beta_1 t$,
the process $B^1_t$ is a Brownian motion 
staying above and reflected
on the line $t \to x + \beta_1 t$,
and the process $X^{1,x+\delta}_t$ is a solution to (1.1)
with $B_t$ replaced by $B^1_t$.

Next we similarly excise the intervals where
$X^{1,x+\delta}_t < B^1_t$. Let us define a new clock
${\cal C}_2(t) = \int_0^t \bone_{\{X^{1,x+\delta}_s \geq B^1_s\}} ds$,
its inverse $b_2(t) = \inf\{s: {\cal C}_2(s)>t\}$,
a random time
$u_2 = \sup\{ {\cal C}_2(t): t\geq 0\}$, and processes
$$\eqalign{
B^2_t &= B^1_{b_2(t)} - \beta_1 (b_2(t) - t), \cr
X^{2,x}_t &= X^{1,x}_{b_2(t)} - \beta_1 (b_2(t) - t), \cr
X^{2,x+\delta}_t &= X^{1,x+\delta}_{b_2(t)} - \beta_1 (b_2(t) - t).
}$$
For $t \in [0,u_2]$,
we have $X^{2,x}_t = x + \beta_1 t$ and
$X^{2,x+\delta}_t = x + \delta + \beta_2 t$.
The process $\{B^2_t, t\in[0,u_2]\}$ is a Brownian motion reflected
on the lines $t \to x + \beta_1 t$
and $t \to x + \delta + \beta_2 t$
and confined to the region between them.
Note that $B^2_0 = x$ a.s. and that the lines
$t \to x + \beta_1 t$
and $t \to x + \delta + \beta_2 t$
intersect at $t = \delta/(\beta_1-\beta_2)$
so necessarily $u_2 \leq \delta/(\beta_1-\beta_2)$.

The time $u_2$ corresponds to the start of the infinite
excursion of $B_t$ below the graph of $X^x_t$.
By excursion theory and Lemma 3.1 (iii), the distribution of
$L^x_{u_2}$ is exponential with mean
$\sigma^2/(2 \beta_2)$. Hence, we may assume that
the process $B^2_t$ is generated in the following way.
Suppose that $B^3_t$ is a Brownian motion starting
from $B^3_0=x$, reflected
on the lines $t \to x + \beta_1 t$
and $t \to x + \delta + \beta_2 t$
and confined to the region between them,
but defined for all $t\in[ 0,\delta/(\beta_1-\beta_2))$ rather than
confined to some random time interval.
Let $L^{3-}_t $ be the local time of $B^3_t$
on the line $t \to x + \beta_1 t$ and let
$Z$ be an exponential random variable
with mean $\sigma^2/(2 \beta_2)$, independent
of $B^3_t$. If $v_3(s) = \inf\{t: L^{3-}_t = s\}$,
then the distributions of the processes
$\{B^2_t, t\in[0,u_2]\}$ and
$\{B^3_t, t\in[0,v_3(Z)]\}$ are the same.

Let $L^{3+}_t$ be the local time of $B^3_t$ accumulated
on the line $t \to x + \delta + \beta_2 t$.
The distribution of $L_\infty^{x+\delta}$
given $\{L^x_\infty = a\}$
is the same as the distribution of 
$L^{3+}_{v_3(a)}$, so we will try to
find an approximate formula for
$\E  L^{3+}_{v_3(a)} $.

We continue our transformations.
Let $B^4_t = B^3_t - x - \beta_1 t$. The process
$B^4_t$ is a Brownian motion starting from $0$,
with drift $-\beta_1$,
reflected on the horizontal axis and the line
$t\to \delta - (\beta_1 - \beta_2)t$.
The processes $L^{3-}_t$ and $L^{3+}_t$ can
be identified with the local
times $L^{4-}_t$ and $L^{4+}_t$
of $B^4_t$ on the horizontal axis and the line
$t\to \delta - (\beta_1 - \beta_2)t$, resp.
Hence, it will suffice to show that the estimate
given in part (i) of the lemma holds for $\E  L^{4+}_{v_4(a)} $.

\medskip
\noindent{\it Step 2}. In this step we will obtain
some estimates for reflected Brownian motions using
excursion theory.
Let $B^5_t$ be a Brownian motion with drift
$-\beta_1$, confined to positive values by
reflection on the horizontal axis.
The Green function $G(z,y)$ for Brownian motion
with drift $-\beta_1$, killed upon hitting
$0$ is given by
$$G(z,y) =  {1 \over \beta_1}
\left[ \exp\left ( {2 \beta_1 z \over \sigma^2}\right)
- 1\right] \exp\left( -{2 \beta_1 y \over \sigma^2} \right),
$$
for $0< z < y < \infty$, by (3.15) in Section 15.3
and Section 15.4.B of Karlin and Taylor (1981).
Let $G^5_H(y)$ denote the Green function for the
excursion law $H_5$ of $B^5_t$ from $0$,
i.e., the function defined by
$$H_5\left(\int_0^\infty 
\bone_{\{e(t) \in [z_1,z_2]\}} dt \right) 
= \int_{z_1}^{z_2} G^5_H(y) dy.$$
A formula analogous to (3.1) yields
$$G^5_H(y) = \lim_{z \downarrow 0}
{1 \over z} G(z,y)
= {2 \over \sigma^2} 
\exp\left( -{2 \beta_1 y \over \sigma^2} \right),
$$
for $y>0$.

Consider some $\delta_1 >0$ and excise
excursions of $B^5_t$ above the level $\delta_1$,
just as we did with the excursions of $B_t$
and $B^1_t$. Let
${\cal C}_3(t) = \int_0^t \bone_{\{ B^5_s \leq \delta_1\}} ds$,
$b_3(t) = \inf\{s: {\cal C}_3(s)>t\}$, and $B^6_t = B^5_{b_3(t)}$.
The process $B^6_t$ is a reflected Brownian motion 
in $[0,\delta_1]$.
Let $G^6_H(y)$ be the Green function for the
excursion law $H_6$ of $B^6_t$ from $0$. It is clear
from the nature of the transformation which generates
$B^6_t$ from the paths of $B^5_t$ that
$G^6_H(y) = G^5_H(y)$ for $y\in(0,\delta_1)$.
Hence,
$$H_6(\nu) = \int_0^{\delta_1} G^6_H(y) dy
= \int_0^{\delta_1} {2 \over \sigma^2} 
\exp\left( -{2 \beta_1 y \over \sigma^2} \right) dy
= {1 \over \beta_1} \left 
[ 1 - \exp\left( -{2 \beta_1 \delta_1 \over \sigma^2} \right)
\right].
$$
Let $L^{6-}_t$ and $L^{6+}_t$ denote the local time
of $B^6_t$ at $0$ and $\delta_1$, resp.
Let $v_6(s) = \inf\{t: L^{6-}_t = s\}$. The random
variable $v_6(s)$ is the sum of the lifetimes
of excursions of $B^6_t$ from $0$ which occur 
before $L^{6-}_t$ reaches the level $s$.
The last formula and excursion theory give
$$\E  v_6(s) = s{1 \over \beta_1} \left 
[ 1 - \exp\left( -{2 \beta_1 \delta_1 \over \sigma^2} \right)
\right]
\df s\eta(\delta_1). \eqno(3.5)
$$

Next we will derive an estimate for $H_6(\nu>t)$.
Recall that $Q^z_{-\beta_1}$ denotes
the distribution of Brownian motion with drift $-\beta_1$,
killed upon hitting $0$. Let $\wh Q^z_{-\beta_1}$
denote the distribution of Brownian motion
starting from $z\in(0,\delta_1)$, with drift $-\beta_1$,
reflected at $\delta_1$, and killed upon hitting $0$.
It is easy to see that
$$\wh Q^z_{-\beta_1} (\nu > t) \leq 
Q^z_{-\beta_1}(\nu > t),$$
for all $t>0$ and $z \in(0,\delta_1)$.
By Lemma 3.1 (i), Theorem 5.1 (iii)
of Burdzy (1987), and scaling,
$$\eqalignno{H_6(\nu>t) &\leq \lim_{z \downarrow 0}
{1 \over z}
\wh Q^z_{-\beta_1} (\nu > t) \cr
&\leq\lim_{z \downarrow 0} {1 \over z}
Q^z_{-\beta_1} (\nu > t) = H_5(\nu > t) \leq
\int_t^\infty {1\over \sigma} (2 \pi s^3)^{-1/2} ds.&(3.6)}$$
A simple argument based on scaling
and the Markov property applied at times
$t= k\delta_1^2$, $k=1,2,\dots$,
shows that there exists
a constant $c_1>0$, such that
$$\wh Q^z_{-\beta_1} (\nu > t) \leq
\exp(- c_1 t \sigma^2/ \delta_1^2),\eqno(3.7)$$
for all $t> \delta_1^2/\sigma^2$ and $z \in(0,\delta_1)$.
Another standard estimate is
$$\wh Q^z_{-\beta_1} (\nu \geq \delta_1^2/\sigma^2) \leq z c_2.$$
This combined with the previous estimate gives
(with possibly new values for the constants),
$$\wh Q^z_{-\beta_1} (\nu > t) \leq z c_2
\exp(- c_1 t \sigma^2/ \delta_1^2),$$
for all $t> \delta_1^2/\sigma^2$ and $z \in(0,\delta_1)$.
We obtain from this an estimate analogous to (3.6)
but applicable for $t> \delta_1^2/\sigma^2$:
$$H_6(\nu>t) \leq c_2 \exp(- c_1 t \sigma^2/ \delta_1^2). \eqno(3.8)
$$
Since the excursion process is a Poisson point
process, we have from (3.6) and (3.8),
$$\eqalignno{
\var v_6&(s) = s \int_0^\infty t^2 H_6(\nu \in dt) \cr
& \leq s \int_0^{\delta_1^2/\sigma^2} 
{1\over \sigma} t^2 (2 \pi t^3)^{-1/2} dt
+ s \int_{\delta_1^2/\sigma^2}^\infty {1\over \sigma}
(\delta_1^2/\sigma^2)^2 (2 \pi t^3)^{-1/2} dt \cr
&\qquad + s \int_{\delta_1^2/\sigma^2}^\infty t^2 
c_2 {\delta_1^2 \over c_1 \sigma^2} 
\exp(- c_1 t \sigma^2/ \delta_1^2) dt \cr
&\leq c_3 s \delta_1^3/\sigma^4 
+ c_4 s \delta_1^3/\sigma^4 + c_5 s \delta_1^8/\sigma^8
\leq c_6 s \delta_1^3. &(3.9)
}$$

\medskip
\noindent{\it Step 3}. We will find a link between processes
reflected on sloped lines (in space-time) and within an interval.
We will need to define some more variables.
First of all, $s_0>0$ should be considered a small
constant whose value will be chosen later in the proof
and which does not change with $\delta$.
Recall $\eta$ defined in (3.5).
Let $u_0>0$ and $\delta_1 \in(0, \delta)$ be defined
by the following two equations
$u_0 = (\delta -\delta_1)/(\beta_1 - \beta_2)$,
and $s_0 = u_0/\eta(\beta_1, \delta,\sigma^2)$.

Recall that $B^6_t$ is a reflected Brownian motion
in $[0,\delta_1]$ and note that now $\delta_1$
is defined relative to $\delta$.
Let $B^7_t$ be the analogous reflected Brownian motion
in $[0,\delta]$.

Note that 
$\delta - (\beta_1 -\beta_2) t > \delta_1$
for $t \in(0,u_0)$. 
Hence, on the interval
$(0,u_0)$, the upper reflecting boundary for $B^6_t$
lies below that for $B^4_t$. This relationship
between the upper reflecting boundaries implies that
the excursion measure distribution of the lifetime 
of an excursion from $0$
of the process $B^4_t$ is stochastically larger
than that for an excursion of $B^6_t$, for
excursions within the interval $(0,u_0)$. It follows
that one can construct $B^4_t$ and $B^6_t$ on a common
probability space so that 
$v_6(s)\land t  \leq v_4(s) \land t $
for all $t \leq u_0$, where
$v_4(s) = \inf\{t: L^{4-}_t = s\}$. 
On the other hand, 
$\delta - (\beta_1 -\beta_2) t < \delta$
for $t >0$, so
the analogous relationship for $B^7$
goes in the opposite way, i.e.,
$v_7(s)\land t  \geq v_4(s) \land t $
for all $t \geq 0$. 

Although the process $B^4_t$ starts from $0$,
by construction, it will be necessary to consider
the case when it starts from some other value;
in other words, we will now consider a process with
the same transition probabilities but a different
starting point. The starting point $y$ will be
reflected in the notation by writing $\P ^y$ or $\E ^y$,
as usual.

Let $T^i_0$ be the hitting time of 0 
for the process $B^i_t$ for  $i=4,6,7$.
By the previous remarks,
we can construct versions of $B^4_t$ and $B^7_t$ on the
same probability space
so that they start from the same point $y$ and
$T^4_0 \land t \leq T^7_0 \land t$ for $t\leq u_0$.

By the strong Markov property applied
at $T^7_0$, we have 
$\E ^y v_7(s) = \E ^y T^7_0 + \E ^0 v_7(s)$.
It follows easily from (3.7), applied to $\delta$
rather than $\delta_1$, that
$$\E ^y T^7_0 \leq c_7 \delta^2 /\sigma^2.  \eqno(3.10)$$
Consider arbitrarily small $\eps\in(0,1/4)$.
We obtain using (3.5) (applied with $\delta_1$
replaced by $\delta$) and (3.10), 
$$ \E ^y (v_4(s_0)) \leq \E ^y (v_7(s_0 ))
\leq \E ^y T^7_0 + \E ^0 (v_7(s_0 ))
\leq c_7 \delta^2 /\sigma^2 + u_0. $$
For small $\delta>0$, 
(3.5) shows that $\eta(\delta)$ is approximately $2 \delta/\sigma^2$.
Hence, $u_0 = s_0\eta(\delta)$ is approximately equal to
$2 s_0\delta/\sigma^2$. This shows that for small
$\delta$, the last displayed inequality yields
$$ \E ^y (v_4(s_0 ))
\leq  u_0(1+\eps) = s_0 \eta(\delta)(1+\eps)
\leq 2 s_0 \delta (1+\eps)^2/\sigma^2. \eqno(3.11)$$

Next we will find a lower bound for the same
quantity.

By the strong Markov property applied
at $T^6_0$, we have 
$\E ^y v_6(s) = \E ^y T^6_0 + \E ^0 v_6(s)$
and
$\var( v_6(s)\mid B^6_0 = y) 
= \var( T^6_0\mid B^6_0 = y) + \var (v_6(s)\mid B^6_0 = 0)$.
We have an estimate analogous to (3.10):
$$\E ^y T^6_0 \leq c_7 \delta^2_1 /\sigma^2,  \eqno(3.12)$$
and another estimate following from (3.7):
$$\var( T^6_0\mid B^6_0 = y) \leq c_8 \delta^4_1 /\sigma^4, \eqno(3.13)$$
for any $y \in [0,\delta_1]$.

We obtain using (3.5) and (3.12), 
$$ \E ^y (v_6(s_0(1-\eps) ))
\leq \E ^y T^6_0 + \E ^0 (v_6(s_0 (1-\eps)))
\leq c_7 \delta^2_1 /\sigma^2 + s_0(1-\eps)\eta(\delta_1). $$
For small $\delta>0$, $\delta_1$ is also small
and (3.5) shows that $\eta(\delta_1)$ is about $2 \delta_1/\sigma^2$.
Hence, $s_0(1-\eps)\eta(\delta_1)$ is approximately equal to
$2 s_0(1-\eps)\delta_1/\sigma^2$. This shows that for small
$\delta$, the last displayed inequality yields
$$ \E ^y (v_6(s_0 (1-\eps)))
\leq  s_0(1-\eps/2)\eta(\delta_1)
\leq  s_0(1-\eps/2)\eta(\delta)
= u_0 (1-\eps/2). \eqno(3.14)$$
A similar estimate for the variance follows from (3.9)
and (3.13), for small $\delta$,
$$\var( v_6(s_0 (1-\eps))\mid B^6_0 = y)
\leq c_8 \delta^4_1 /\sigma^4 + 
c_6 s_0 (1-\eps) \delta_1^3/\sigma^4
\leq c_9 s_0(1-\eps) \delta_1^3/\sigma^4.\eqno(3.15)$$
This estimate, (3.14) and
the Chebyshev inequality yield,
$$\P ^y(v_6(s_0 (1-\eps)) \geq u_0)
\leq {c_9 s_0 (1-\eps) \delta_1^3/\sigma^4
\over (\eps u_0/2)^2}
\leq {c_{10}  \delta_1^3
\over \eta(\delta) \eps^2 u_0 \sigma^4}
= {c_{10}  \delta_1^3
\over \eta^2(\delta) \eps^2 s_0 \sigma^4}.
$$
For small $\delta$
we have 
$$\delta/\sigma^2 < \eta(\delta) < 4 \delta/\sigma^2. \eqno(3.16)$$
Hence,
$$\P ^y(v_6(s_0 (1-\eps)) \geq u_0)
\leq {c_{11}  \delta_1^3
\over (\delta/\sigma^2)^2 \eps^2 s_0 \sigma^4}
\leq {c_{11}  \delta_1
\over  \eps^2 s_0 }. \eqno(3.17)
$$
We have from (3.14)-(3.17), for small $\delta$,
$$\eqalign{
\E ^y (v_6(s_0(1-\eps)) )^2 
&= \var( v_6(s_0 (1-\eps))\mid B^6_0 = y)
+ (\E ^y v_6(s_0(1-\eps)) )^2 \cr
& \leq 
c_9 s_0(1-\eps) \delta_1^3/\sigma^4
+ (s_0(1-\eps/2)\eta(\delta))^2 \cr
&\leq 32 s_0^2 (1-\eps/2)^2 \delta^2 /\sigma^4 .
}$$
We use this estimate, (3.5) and (3.16)-(3.17) to obtain,
for sufficiently small $\delta$,
$$\eqalignno{
\E ^y(v_4(s_0)) &\geq
\E ^y (v_4(s_0) \land u_0)\cr
&\geq \E ^y (v_6(s_0) \land u_0) \cr
& \geq \E ^y (v_6(s_0(1-\eps)) \land u_0) \cr
&\geq \E ^y (v_6(s_0(1-\eps)) ) - \E ^y \left[v_6(s_0(1-\eps)) 
\bone _{\{ v_6(s_0(1-\eps)) \geq u_0\}} \right] \cr
& \geq s_0 (1-\eps)\eta(\delta_1) 
- \left( \E ^y (v_6(s_0(1-\eps)) )^2 
\E ^y (\bone _{\{ v_6(s_0(1-\eps)) \geq u_0\}} )^2 \right)^{1/2} \cr
&= s_0 (1-\eps)\eta(\delta_1) 
- \left( \E ^y (v_6(s_0(1-\eps)) )^2 
\P ^y ( v_6(s_0(1-\eps)) \geq u_0 ) \right)^{1/2} \cr
&\geq
2 s_0 (1-\eps)^2\delta_1/\sigma^2 
- \left( 32 s_0^2 (1-\eps/2)^2 \delta^2 /\sigma^4 \cdot
{c_{11}  \delta_1 \over  \eps^2 s_0 } \right)^{1/2}. &(3.18)
}$$
It follows from (3.16) and the definition of $\delta_1$ and $u_0$
that for small $\delta>0$,
$$\delta_1 
=\delta- s_0 \eta(\delta) (\beta_1 -\beta_2)
\geq \delta (1 - 4 s_0 (\beta_1 -\beta_2) /\sigma^2).$$
We will choose sufficiently small $s_0>0$ (relative to
$\sigma, \beta_1 ,\beta_2$ and $\eps$) so that
$\delta_1 > \delta(1-\eps/2)$.
Then the last inequality and (3.18) yield for small $\delta$,
$$
\E ^y(v_4(s_0)) \geq 2 s_0 (1-\eps)^3\delta/\sigma^2.\eqno(3.19)$$

\medskip
\noindent{\it Step 4}.
We will apply induction in order to obtain
estimates for $\E ^0 v_4(js_0)$ with integer $j\geq 1$.
At the time $v_4(s_0)$, the distance between
the reflecting barriers for $B^4_t$ is equal to
$\wt\delta=\delta - (\beta_1 -\beta_2) v_4(s_0)$, 
which is less than $\delta$, so we can use the 
estimates (3.11) and (3.19) with $\delta$ replaced by $\wt \delta$,
assuming that $\delta$ itself is sufficiently 
small for the estimates to hold.
By the strong Markov property,
$$\E ^0 (v_4(2s_0) - v_4(s_0) \mid v_4(s_0))
\leq 2 \wt\delta s_0 (1+\eps)^2/\sigma^2
= 2 [\delta - (\beta_1 -\beta_2) v_4(s_0)] s_0 (1+\eps)^2/\sigma^2,
$$
and so
$$\eqalign{\E ^0 (v_4(2s_0) - v_4(s_0) )
&\leq \E ^0 2 [\delta - (\beta_1 -\beta_2) v_4(s_0)] s_0 (1+\eps)^2/\sigma^2\cr
&\leq 2 [\delta - (\beta_1 -\beta_2) 2 \delta s_0 (1+\eps)^2/\sigma^2] 
s_0 (1+\eps)^2/\sigma^2 \cr
&= 2\delta (s_0/\sigma^2) [(1+\eps)^2 - 2(\beta_1 -\beta_2)  
(s_0/\sigma^2) (1+\eps)^4] .}
$$
It follows that for small $\delta>0$,
$$\eqalign{\E ^0 v_4(2s_0) &= \E ^0 v_4(s_0) + \E ^0 (v_4(2s_0) - v_4(s_0) )\cr
&\leq 2 \delta (s_0/\sigma^2) (1+\eps)^2
+ 2\delta (s_0/\sigma^2) [(1+\eps)^2 - 2(\beta_1 -\beta_2)  
(s_0/\sigma^2) (1+\eps)^4]\cr
& = 2\delta (s_0/\sigma^2) [2(1+\eps)^2 - 2(\beta_1 -\beta_2) 
(s_0/\sigma^2) (1+\eps)^4] .}
$$
More generally,
$$\eqalign{\E ^0 v_4((j+1)s_0) 
&= \E ^0 v_4(js_0) + \E ^0 (v_4((j+1)s_0) - v_4(js_0) )\cr
&\leq \E ^0 v_4(js_0)
+ 2 [ \delta - (\beta_1 -\beta_2) \E ^0 v_4(js_0)] s_0 (1+\eps)^2/\sigma^2\cr
&= \E ^0 v_4(js_0)[1 - 2  (\beta_1 -\beta_2) (s_0/\sigma^2) (1+\eps)^2]
+ 2\delta (s_0/\sigma^2) (1+\eps)^2
.}
$$
From this we obtain by induction,
$$\eqalign{\E ^0 v_4(js_0)
&\leq \E ^0 v_4(s_0) [1 - 2  (\beta_1 -\beta_2) (s_0/\sigma^2) (1+\eps)^2]^{j-1} \cr
&\qquad+ 2\delta (s_0/\sigma^2) (1+\eps)^2
\sum_{k=0}^{j-2} 
[1 - 2  (\beta_1 -\beta_2) (s_0/\sigma^2) (1+\eps)^2]^k \cr
&\leq 2\delta (s_0/\sigma^2) (1+\eps)^2
[1 - 2  (\beta_1 -\beta_2) (s_0/\sigma^2) (1+\eps)^2]^{j-1} \cr
&\qquad+ 2\delta (s_0/\sigma^2) (1+\eps)^2
\sum_{k=0}^{j-2} 
[1 - 2  (\beta_1 -\beta_2) (s_0/\sigma^2) (1+\eps)^2]^k \cr
&=2\delta (s_0/\sigma^2) (1+\eps)^2
{ 1 - [1 - 2  (\beta_1 -\beta_2) (s_0/\sigma^2) (1+\eps)^2]^j
\over 1 - [1 - 2  (\beta_1 -\beta_2) (s_0/\sigma^2) (1+\eps)^2]}\cr
&={ \delta \over \beta_1 -\beta_2}
(1- [1 - 2  (\beta_1 -\beta_2) (s_0/\sigma^2) (1+\eps)^2]^j )
.}$$

Now fix an arbitrary $a>0$, an arbitrarily small
$\eps >0$, and choose a sufficiently small
small $s_0>0$ so that $\delta_1 > \delta(1-\eps/2)$,
and such that
for some integer $j$ we have $js_0 = a$, and,
moreover, $j$ is sufficiently large to imply
the following:
$$\eqalign{
&{ \delta \over \beta_1 -\beta_2}
(1- [1 - 2  (\beta_1 -\beta_2) (s_0/\sigma^2) (1+\eps)^2]^j ) \cr
& \qquad
={ \delta \over \beta_1 -\beta_2}
(1- [1 - 2  (\beta_1 -\beta_2) (s_0/\sigma^2) (1+\eps)^2]^{a/s_0} ) \cr
&\qquad
\leq 
{\delta \over \beta_1 -\beta_2}
(1- \exp (-2 a (\beta_1 -\beta_2)  (1+\eps)^3/\sigma^2) )
.}$$
Then for sufficiently small $\delta>0$ we have,
$$\E ^0 v_4(a) = \E ^0 v_4(j s_0)
\leq 
{\delta \over \beta_1 -\beta_2}
(1- \exp ( -2 a (\beta_1 -\beta_2)  (1+\eps)^3/\sigma^2) ).\eqno(3.20)
$$
A completely analogous argument using (3.19)
in place of (3.11) yields
$$\E ^0 v_4(a) \geq {\delta \over \beta_1 -\beta_2}
(1- \exp (-2 a (\beta_1 -\beta_2)  (1-\eps)^4/\sigma^2) )
. \eqno(3.21)
$$
Since $\eps>0$ is arbitrarily small,
a standard argument based on (3.20)-(3.21) gives
for $\delta\downarrow 0$,
$$\E ^0 v_4(a) = {\delta \over \beta_1 -\beta_2}
(1- \exp (-2 a (\beta_1 -\beta_2)  /\sigma^2) ) + o(\delta)
. \eqno(3.22)
$$

\medskip
\noindent{\it Step 5}. The last part of the proof
exploits a relationship between local time
and certain stopping times.
Recall the local times $L^{4-}_t$ and $L^{4+}_t$,
introduced earlier in the proof. We have
for some standard Brownian motion $B^8_t$,
$$B^4_t = B^8_t  - \beta_1 t + L^{4-}_t - L^{4+}_t.$$
One has to check that the normalization
of the local time, defined relative to the
normalization of the excursion laws in (3.1),
is the correct one for the above ``L\'evy formula.''
This can be done, for example, by comparing
our normalizations with those in Theorems
3.6.17 and 6.2.23 in Karatzas and Shreve (1988).

Note that the $\sigma$-fields generated by $B^4_t$
and $B^8_t$ are identical so $v_4(a)$
is a stopping time for $B^8_t$. 
We have $B^4_{v_4(a)} = 0$ 
and $ L^{4-}_{v_4(a)} = a$, so
$$0 = B^8_{v_4(a)}  - \beta_1 {v_4(a)}+ a - L^{4+}_{v_4(a)}.\eqno(3.23)$$
Since $v_4(a) $ is bounded by $ \delta/ (\beta_1 - \beta_2)$
the optional stopping theorem yields
$\E  B^8_{v_4(a)} = 0$, and so, using (3.22),
$$\E  L^{4+}_{v_4(a)} = a - \E  \beta_1 {v_4(a)}
= a - \delta{\beta_1 \over \beta_1 -\beta_2}
(1- \exp (-2 a (\beta_1 -\beta_2)  /\sigma^2) ) + o(\delta)
. \eqno(3.24)$$
Now recall that $L^{4+}_{v_4(a)}$ has the same
distribution as $L^{x+\delta}_\infty$
given $\{L^x_\infty = a\}$. This observation
and the last formula
complete the proof of part (i) of the lemma.

(ii) The proof of part (ii) of the lemma
uses a formula analogous to (3.24), 
but requires some additional work.

Recall that $\delta$ was positive in part (i) of the proof;
it will be negative in the present part.

Recall the transformations of $B_t$ from the
proof of (i).
It is easy to see that analogous
transformations in the current case
do not lead to 
$B^4_t$ which is a Brownian motion starting from $0$,
with drift $-\beta_1$,
reflected on the horizontal axis and the line
$t\to \delta - (\beta_1 - \beta_2)t$ (with $\delta>0$),
but instead they give a Brownian motion $\wt B^4_t$ starting from $0$,
with drift $-\beta_2$,
reflected on the horizontal axis and the line
$t\to \delta + (\beta_1 - \beta_2)t$ (with $\delta<0$).

A subtle but significant difference from (i)
is that the infinite excursion of $B_t$ from
the graph of $X^x_t$ will go in the 
direction of the graph of $X^{x+\delta}_t$
and so it will generate some more local time.
By Lemma 3.1 (i) and Remark 3.2 (ii), the excursions with finite
lifetimes have the same intensities 
for Brownian motions with drifts $\beta_2$ and
$-\beta_2$, so we can use estimate (3.24)
for the portion of the local time generated
before the last, infinite excursion of 
$B_t$ from the graph of $X^x_t$. The estimate
has to be modified as $\beta_1$ has to be replaced
by $\beta_2$, and so we obtain
$$a - |\delta|{\beta_2 \over \beta_1 -\beta_2}
(1- \exp (-2 a (\beta_1 -\beta_2)  /\sigma^2) ) + o(\delta)
. \eqno(3.25)$$
To this we will have to add the local
time spent by $B_t$ on the graph of $X^{x+\delta}_t$
during its final, infinite excursion from
the graph of $X^x_t$. The rest of the proof is devoted
to that calculation.

Let $U$ be the first time when the
final, infinite excursion of $B_t$ from
the graph of $X^x_t$ hits the graph of $X^{x+\delta}_t$.
Let $\delta_1 = |X^x_U - X^{x+\delta}_U|$.
First, we will condition on $\delta_1$.
The process $\{B_t, t\geq U\}$ is a Brownian motion
conditioned not to hit the line 
$t \to B_U +\delta_1 + \beta_2 t$.
By subtracting the drift and flipping the process
to the other side of the horizontal axis,
we may consider a Brownian motion $B^9_t$
starting from $\delta_1$, with drift $\beta_2$,
conditioned not to hit $0$. We will estimate
the amount of the local time this process 
spends on the graph of a solution $Y_t$
to (1.1) with $\beta_2 $ replaced by $-(\beta_1-\beta_2)$,
$\beta_1$ replaced by $0$, and $B_t$ replaced by
$B^9_t$.

Let $H^{\delta_1}$ be the excursion law for excursions
above the level $\delta_1$ for Brownian motion
with drift $\beta_2$, conditioned not to hit $0$.
Let $F_\infty$ denote the set of excursions with infinite lifetime.
We will compute $H^{\delta_1}(F_\infty)$.
Let $Q^z_{\beta_2}$ be the distribution of Brownian motion
with drift ${\beta_2}$.
Then 
$$\eqalign{
H^{\delta_1}(F_\infty)
&= \lim_{z\downarrow 0} {1\over z}\cdot
Q^{{\delta_1}+z}_{\beta_2}( T_\infty < T_{\delta_1} \mid T_\infty < T_0)\cr
&= \lim_{z\downarrow 0} {1\over z}\cdot
{Q^{{\delta_1}+z}_{\beta_2}( T_\infty < T_{\delta_1} \hbox{   and   } T_\infty < T_0)
\over
Q^{{\delta_1}+z}_{\beta_2}( T_\infty < T_0)}\cr
&= \lim_{z\downarrow 0} {1\over z}\cdot
{Q^{{\delta_1}+z}_{\beta_2}( T_\infty < T_{\delta_1} )
\over
Q^{{\delta_1}+z}_{\beta_2}( T_\infty < T_0)}.
}$$
Recall that the scale function $S(y)$ for Brownian motion
with drift ${\beta_2}$ is equal to \break 
$\exp(-2{\beta_2} y/\sigma^2)$.
This gives
$$\eqalignno{
H^{\delta_1}(F_\infty)
&= \lim_{z\downarrow 0} {1\over z}\cdot
{S({\delta_1}+z) - S({\delta_1}) \over S(\infty) - S({\delta_1})}
\cdot
{S(\infty) - S(0) \over S({\delta_1}+z) - S(0)}
\cr
&= \lim_{z\downarrow 0} {1\over z}\cdot
{\exp(-2{\beta_2} ({\delta_1}+z)/\sigma^2)-\exp(-2{\beta_2} {\delta_1}/\sigma^2)
\over
0-\exp(-2{\beta_2} {\delta_1}/\sigma^2)}
\cdot
{0-1 \over \exp(-2{\beta_2} ({\delta_1}+z)/\sigma^2)-1}\cr
&= {2{\beta_2}  \over
\sigma^2 [ 1 - \exp(-2{\beta_2} {\delta_1}/\sigma^2)]}.
}$$
If we fix arbitrarily small $\eps>0$
then for sufficiently small ${\delta_1}>0$ we have
$${1\over{\delta_1}} \leq H^{\delta_1}(F_\infty)
\leq ( 1+\eps){1\over{\delta_1}}. \eqno (3.26)$$

We proceed to calculate the expected time
to hit ${\delta_1}$ for Brownian motion
with drift ${\beta_2}$,
starting from ${\delta_1}-z$ and conditioned
not to hit $0$, where $z\in(0,\delta_1)$.
If we take
$$s(z) = \exp\left [-\int_0^z 2 {\beta_2}/\sigma^2 dy\right]
= \exp(- 2 {\beta_2} z/\sigma^2), $$
and
$$S(z) = \int_0^z s(y) dy
= {\sigma^2 \over 2 {\beta_2}} [1 -\exp(- 2 {\beta_2} z/\sigma^2)], $$
then formula (9.9) on p. 264 of Karlin and Taylor (1981)
yields
$$\eqalign{
&\E ^{{\delta_1}-z} (T_{\delta_1} \mid T_{\delta_1} < T_0)\cr
&=
{2 [S({\delta_1}) - S({\delta_1}-z)] \over S({\delta_1})S({\delta_1}-z)}
\int_0^{{\delta_1}-z}
{S^2(y) \over \sigma^2 s(y)} dy
+ 2 \int_{{\delta_1}-z}^{\delta_1}
{S(y) [ S({\delta_1}) - S(y)] \over 
\sigma^2 s(y) S({\delta_1})} dy\cr
&= 
{2 \left[{\sigma^2 \over 2 {\beta_2}}
[1 -\exp(- 2 {\beta_2} {\delta_1}/\sigma^2)]
 - {\sigma^2 \over 2 {\beta_2}} 
[1 -\exp(- 2 {\beta_2} ({\delta_1} - z)/\sigma^2)]
\right] \over 
{\sigma^2 \over 2 {\beta_2}} [1 -\exp(- 2 {\beta_2} {\delta_1}/\sigma^2)]
{\sigma^2 \over 2 {\beta_2}} [1 -\exp(- 2 {\beta_2} ({\delta_1}-z)/\sigma^2)]
} \times \cr
&\qquad
\int_0^{{\delta_1}-z}
{{\sigma^4 \over 4 {\beta_2}^2} [1 -\exp(- 2 {\beta_2} y/\sigma^2)]^2
\over \sigma^2 \exp(- 2 {\beta_2} y/\sigma^2)} dy \cr
&\quad
+ 2 \int_{{\delta_1}-z}^{\delta_1} \Bigg (
{{\sigma^2 \over 2 {\beta_2}} [1 -\exp(- 2 {\beta_2} y/\sigma^2)]  
\over 
\sigma^2 \exp(- 2 {\beta_2} y/\sigma^2) 
{\sigma^2 \over 2 {\beta_2}} [1 -\exp(- 2 {\beta_2} {\delta_1}/\sigma^2)]} \times\cr
&\qquad
\left[{\sigma^2 \over 2 {\beta_2}} [1 -\exp(- 2 {\beta_2} {\delta_1}/\sigma^2)]
- {\sigma^2 \over 2 {\beta_2}} [1 -\exp(- 2 {\beta_2} y/\sigma^2)]\right]\Bigg) dy \cr
&= 
{2 [\exp(- 2 {\beta_2} ({\delta_1} - z)/\sigma^2)
-\exp(- 2 {\beta_2} {\delta_1}/\sigma^2)]
\over 
{\sigma^2 \over 2 {\beta_2}} [1 -\exp(- 2 {\beta_2} {\delta_1}/\sigma^2)]
[1 -\exp(- 2 {\beta_2} ({\delta_1}-z)/\sigma^2)]
} \times \cr
&\qquad
\int_0^{{\delta_1}-z}
{{\sigma^4 \over 4 {\beta_2}^2} [1 -\exp(- 2 {\beta_2} y/\sigma^2)]^2
\over \sigma^2 \exp(- 2 {\beta_2} y/\sigma^2)} dy \cr
&\quad
+ 2 \int_{{\delta_1}-z}^{\delta_1}
{{\sigma^2 \over 2 {\beta_2}} [1 -\exp(- 2 {\beta_2} y/\sigma^2)] 
[ \exp(- 2 {\beta_2} y/\sigma^2)-\exp(- 2 {\beta_2} {\delta_1}/\sigma^2)] 
\over 
\sigma^2 \exp(- 2 {\beta_2} y/\sigma^2) 
 [1 -\exp(- 2 {\beta_2} {\delta_1}/\sigma^2)]} dy.}$$
The expected lifetime of an excursion
below ${\delta_1}$ for the Brownian motion
with drift ${\beta_2}$,
starting from ${\delta_1}$ and conditioned
not to hit $0$ is therefore equal to
$$\eqalignno{
\lim_{z\downarrow 0} {1\over z}
&\E ^{{\delta_1}-z} (T_{\delta_1} \mid T_{\delta_1} < T_0) \cr
& =
{2 {2{\beta_2}\over\sigma^2}
\exp(- 2 {\beta_2} {\delta_1}/\sigma^2)
\over 
{\sigma^2 \over 2 {\beta_2}} [1 -\exp(- 2 {\beta_2} {\delta_1}/\sigma^2)]
[1 -\exp(- 2 {\beta_2} {\delta_1}/\sigma^2)] } 
\int_0^{{\delta_1}}
{{\sigma^4 \over 4 {\beta_2}^2} [1 -\exp(- 2 {\beta_2} y/\sigma^2)]^2
\over \sigma^2 \exp(- 2 {\beta_2} y/\sigma^2)} dy \cr
& \leq c_1 {\delta_1}, & (3.27)}$$
for small ${\delta_1}>0$ and some $c_1$ depending
on ${\beta_2}$ and $\sigma^2$ but not on ${\delta_1}$.

Fix arbitrarily small $\eps >0$.
We are ready to derive estimates for the total amount of local
time, say $L^{9+}_\infty$, that $B^9_t$ spends on the graph of $Y_t$.

On one hand, the estimate (3.26) shows that $L^{9+}_\infty$ is
stochastically bounded by an exponential random variable
with mean ${\delta_1}$, for sufficiently small
${\delta_1} >0$.

Let $v^+_9(a)$ be the time spent by $B^9_t$
between $Y_t$ and the horizontal axis before the time
when $L^{9+}_t$
hits $a$. Since $Y_t$ is non-increasing,
the estimate (3.27) can be used as an upper bound
for the expected duration of an excursion below $Y_t$,
for every $t\geq 0$.
Fix arbitrarily large $b<\infty$
and arbitrarily small $\eps>0$. 
We have from excursion theory,
$$\E  v^+_9(b{\delta_1}) \leq b{\delta_1} c_1 {\delta_1},$$
and so, for sufficiently small ${\delta_1}>0$,
$$\P (v^+_9(b{\delta_1}) \geq {\delta_1}\eps) 
\leq { b c_1 {\delta_1}^2 \over {\delta_1} \eps}
= {b c_1\over \eps} {\delta_1} < \eps.$$
We see that with probability greater than
$1-\eps$, the distance between $Y_t$ and the horizontal axis
remains greater than 
${\delta_1} - {\delta_1}\eps(\beta_1 - \beta_2)$,
at least until the time when $L^{9+}_t$ exceeds $b{\delta_1}$.
On this time interval and given this event,
the intensity for the arrival
process of the infinite excursion is bounded
above by $(1+\eps)/({\delta_1}(1-\eps))$, by (3.26).
Hence, $\E L^{9+}_\infty / \delta_1$
can be made arbitrarily close to $1$,
by choosing large $b$, then small $\eps$ and finally small $|\delta|>0$
(note that $\delta_1 \leq |\delta|$).

Finally, in order to obtain an unconditioned
estimate for $\E  L^{9+}_\infty$, we have to 
average over the possible values of ${\delta_1}$.
Let $\wt v_4(a)$ be the time when the local
time of $\wt B^4_t$ (defined earlier in the proof
of part (ii)) reaches $a$. The same argument which
gives (3.22) yields the following estimate,
$$\eqalign{
\E  L^{9+}_\infty & =
\E \delta_1 + o(\delta) \cr
&=
{|\delta|} -\E  (\beta_1-\beta_2) \wt v_4(a) + o(\delta)\cr
&= {|\delta|} -
{|\delta|}{\beta_1 -\beta_2\over \beta_1 -\beta_2}
(1- \exp (-2 a (\beta_1 -\beta_2)  /\sigma^2) ) + o(\delta) \cr
&= {|\delta|} \exp (-2 a (\beta_1 -\beta_2)  /\sigma^2) ) + o(\delta) .}$$
Adding this quantity to (3.25) gives the formula
in Lemma 3.4 (ii).
\qed
\bigskip

\noindent{\bf Lemma 3.5}. {\sl 
Fix $x,a,\beta_1,\beta_2>0$ and assume that $\beta_1-\beta_2 >0$. Then
$$\var (L^{x+\delta}_\infty \mid L^x_\infty = a)
= {\delta \over \beta_1 -\beta_2}
(1- \exp (-2 a (\beta_1 -\beta_2)  /\sigma^2) ) + o(\delta),$$
for $\delta \downarrow 0$. The same formula holds if
$x<0$ and $\delta \uparrow 0$.
}

\bigskip
\noindent{\bf Proof}. 
First suppose that $x,\delta >0$ and recall the notation
and definitions from the proof of Lemma 3.4 (i).
It follows from (3.23) that
$$  L^{4+}_{v_4(a)} - a =  B^8_{v_4(a)}  - \beta_1 v_4(a).$$
We have using (3.22),
$$\eqalignno{
\E ( &  L^{4+}_{v_4(a)} - a)^2
= \E  ( B^8_{v_4(a)}  - \beta_1 v_4(a))^2 & (3.28)\cr
&= \E   ( B^8_{v_4(a)})^2 
- 2 \beta_1 \E  \left[ B^8_{v_4(a)} v_4(a)\right]+ \E  (v_4(a))^2\cr
&= \E   v_4(a)
- 2 \beta_1 \E  \left[ B^8_{v_4(a)} v_4(a)\right]+ \E  (v_4(a))^2 \cr
& = {\delta \over \beta_1 -\beta_2}
(1- \exp (-2 a (\beta_1 -\beta_2)  /\sigma^2) ) + o(\delta)
- 2 \beta_1 \E  \left[ B^8_{v_4(a)} v_4(a)\right]+ \E  (v_4(a))^2.
}$$
Recall that $v_4(a) $ is bounded by $ \delta/ (\beta_1 - \beta_2)$.
Hence,
$$\E  (v_4(a))^2
\leq \delta^2/ (\beta_1 - \beta_2)^2,\eqno(3.29)$$
and
$$\eqalignno{\E   B^8_{v_4(a)} v_4(a)
&\leq \left( \E  (B^8_{v_4(a)})^2 \E (v_4(a))^2 \right)^{1/2} &(3.30)\cr
&= \left( \E  v_4(a) \E (v_4(a))^2 \right)^{1/2} \cr
& \leq \left( \left[{\delta \over \beta_1 -\beta_2}
(1- \exp (-2 a (\beta_1 -\beta_2)  /\sigma^2) ) + o(\delta)\right]
\cdot \delta^2/ (\beta_1 - \beta_2)^2 \right)^{1/2} \cr
&= o(\delta).
}$$
Combining (3.28)-(3.30) yields
$$
\E ( L^{4+}_{v_4(a)} - a)^2
= {\delta \over \beta_1 -\beta_2}
(1- \exp (-2 a (\beta_1 -\beta_2)  /\sigma^2) ) + o(\delta).$$
This implies
$$\eqalign{
\var L^{4+}_{v_4(a)}
& = \E  ( L^{4+}_{v_4(a)} - \E  L^{4+}_{v_4(a)})^2\cr
& = \E ( L^{4+}_{v_4(a)} - a)^2
- ( \E  L^{4+}_{v_4(a)} - a)^2 \cr
& = {\delta \over \beta_1 -\beta_2}
(1- \exp (-2 a (\beta_1 -\beta_2)  /\sigma^2) ) + o(\delta)\cr
&\qquad - 
\left[
\delta{\beta_1 \over \beta_1 -\beta_2}
(1- \exp (-2 a (\beta_1 -\beta_2)  /\sigma^2) ) + o(\delta)\right]^2 \cr
& = {\delta \over \beta_1 -\beta_2}
(1- \exp (-2 a (\beta_1 -\beta_2)  /\sigma^2) ) + o(\delta).
}$$
Similarly to the proof of Lemma 3.4 (i), we have
$\var(L^{x+\delta}_\infty \mid L^x_\infty = a) 
= \var L^{4+}_{v_4(a)}$, which combined with the last
formula proves the lemma in the case $x>0$.

Now consider the case $x,\delta<0$. We argue as in the proof
of Lemma 3.4 (ii) that we have to add a contribution
from the local time on $X^{x+\delta}_t$ generated by the
infinite excursion of $B_t$ below $X^x_t$. We have shown
in the proof of Lemma 3.4 (ii) that the local time
on $X^{x+\delta}_t$ generated by the
infinite excursion is stochastically bounded by an 
exponential random variable with mean $\delta$ so its variance 
is bounded by $2 \delta^2$, and, therefore, the contribution
to the variance from the infinite excursion is negligible. The same
formula holds in the case $x<0$ as in the case $x>0$.
\qed

\bigskip
Recall that $L^x_t$ denotes the local time of
$B_t - X^x_t$ at $0$.
The following result is analogous to Trotter's theorem
on the joint continuity of local times for Brownian motion
(see Karatzas and Shreve (1988) or Knight (1981)).

\bigskip
\noindent{\bf Theorem 3.6}. {\sl
Assume that $\beta_1, \beta _2 >0$ and $\beta_1 - \beta_2 >0$. 
There exists a version of the process $(x,t) \to L^x_t$
which is jointly continuous in both variables.
}

\bigskip

\proof 
Note that  $X_t^x$ and $X_t^y$ increase at the same rate when 
$B_t$ does not lie between $X_t^x$ and $X_t^y$, and by the assumptions on
$\beta_1$ and $\beta_2$,
they grow closer together
when $B_t$ does lie between them.
Therefore, for all $x,y$ and $t\geq 0$,
$$|X_t^x-X_t^y|\leq |x-y|. \eqno (3.31)$$

Define $G(x)=\E L_\infty^x$. 
The excursion law for Brownian motion below
the line $t \to \beta_2 t$ gives mass
$2 \beta_2 /\sigma^2$ to excursions with infinite lifetime,
by Lemma 3.1 (iii).
By excursion theory, the waiting time (in terms
of local time) for the first excursion with infinite
lifetime is exponential with mean $\sigma^2/(2 \beta_2)$.
This says that the distribution
of $L^0_\infty$ is exponential with mean $\sigma^2/(2 \beta_2)$.
This and the strong Markov property applied
at the first time when $B_t$ intersects
$X^x_t$ imply that for some $c_1<\infty $
and all $x$, we have $G(x)\leq c_1$.
An easy conditioning argument that combines this observation with Lemma 3.4
shows that for all $x$ and $y$,
$$|G(x)-G(y)|\leq c_2|x-y|. \eqno (3.32)$$

The process $(X_t^x, B_t)$ is strong Markov, and $L_t^x$ is an additive
functional. So by the Markov property,
$$\eqalign{\E[L_\infty^x-L_t^x\mid \F_t]&=\E[L_\infty^x\circ \theta_t
\mid \F_t]\cr
&=G(X_t^x-B_t).\cr}$$
Therefore
$$\E[L_\infty^x-L_t^x\mid \F_t]\leq c_1.$$
Also, using (3.31) and (3.32), 
$$\eqalign{|\E[(L_\infty^x-L_\infty^y)-(L_t^x-L_t^y)\mid \F_t]|&
=|G(X_t^x-B_t)-G(X_t^y-B_t)|\cr
&\leq c_2|X_t^x-X_t^y|\cr
&\leq c_2|x-y|.\cr}$$
By Bass (1995), Proposition I.6.14, 
$$\E[\sup_t|L_t^x-L_t^y|^4]\leq c_4|x-y|^2. \eqno (3.33)$$
By Kolmogorov's criterion and standard arguments (cf.\ the proof of
Proposition I.6.16 of Bass (1995)), we deduce that there exists a version
of $L_t^x$ that is jointly continuous in $x$ and $t$.
\qed

\bigskip

A classical Ray-Knight theorem (see Knight (1981),
Revuz and Yor (1991) or Yor (1997)) asserts, roughly speaking,
that if $L_t^x$ is the local time for the standard Brownian
motion then $x \to L_T^x$
is a diffusion for certain stopping
times $T$. As a part of that theorem, the infinitesimal parameters of the diffusion
are also given. We prove a similar
result for our family of local times, with $T\equiv \infty$.
Recall that we assume that $B_0=0$.

\bigskip
\noindent{\bf Theorem 3.7}. {\sl 
Suppose that $\beta_1, \beta _2 >0$ and $\beta_1 - \beta_2 >0$. 
The distribution
of $L^0_\infty$ is exponential with mean $\sigma^2/(2 \beta_2)$.
The process $\{L^x_\infty, x\geq 0\}$ is a diffusion
with the infinitesimal drift
$$\wt \mu(a) = - {\beta_1 \over \beta_1 -\beta_2}
(1- \exp (-2 a (\beta_1 -\beta_2)  /\sigma^2) ),$$
and infinitesimal variance
$$\wt \sigma^2(a) = {1 \over \beta_1 -\beta_2}
(1- \exp (-2 a (\beta_1 -\beta_2)  /\sigma^2) ).$$
The process $\{L^{-x}_\infty, x\geq 0\}$ is a diffusion
with the infinitesimal drift
$$\wh \mu(a) = -{\beta_2 \over \beta_1-\beta_2}
+ {\beta_1 \over \beta_1-\beta_2}
\exp (-2 a (\beta_1 -\beta_2)  /\sigma^2),$$
and the same infinitesimal variance
$$\wt \sigma^2(a) = {1 \over \beta_1 -\beta_2}
(1- \exp (-2 a (\beta_1 -\beta_2)  /\sigma^2) ).$$
}

\bigskip
\noindent{\bf Proof}.
We have already shown in the proof of Theorem 3.6 that the distribution
of $L^0_\infty$ is exponential with mean $\sigma^2/(2 \beta_2)$.

The Markovian character of the process 
$\{L^x_\infty, x\geq 0\}$ at any fixed ``time'' $x=y$ follows from the 
independence of the Poisson processes of excursions of $B_t$
below and above $X^y_t$. 
The same remark applies to $\{L^{-x}_\infty, x\geq 0\}$.
The infinitesimal
parameters of the processes were calculated in Lemmas 3.2 and 3.3.

The process $x \to L^x_\infty$ is continuous, by Theorem 3.6.
Since its infinitesimal drift is bounded and the infinitesimal variance is
nondegenerate, 
there is a unique (in law) Markov process with this
infinitesimal drift and variance (cf.~Bass (1997), Section IV.3), 
and this Markov process is in fact a 
strong Markov process.
\qed

\bigskip
\noindent{\bf Theorem 3.8}. {\sl
Suppose $\beta_1,\beta_2 >0$ and $\beta_1 - \beta_2 >0$.
For fixed $t>0$, we have a.s.,
for all $x,x_1, x_2\in \R$,
$$\left. {d \over dy} X^y_t \right|_{y=x}
=  \exp (-2 L^x_{t} (\beta_1 -\beta_2)  /\sigma^2),
$$
and
$$ X^{x_2}_t - X^{x_1}_t
= \int _{x_1} ^{x_2}  \exp (-2 L^x_{t} (\beta_1 -\beta_2)  /\sigma^2) dx. 
$$
}

\bigskip
\noindent{\bf Proof}. First we will prove an estimate
analogous to (3.22) except that it will hold
for $v_4(a)$ itself rather than its expectation.
Recall the notation and definitions from the proof of Lemma 3.4 (i).

The following estimate is completely analogous
to (3.17) except that we state it for the process
$B^7_t$ rather than $B^6_t$, so $\delta_1$
is replaced by $\delta$ in the bound.
$$\P ^y(v_7(s_0 (1-\eps)) \geq u_0)
\leq 
{c_{11}  \delta \over  \eps^2 s_0 \sigma^2}. 
$$
We can further modify the estimate by
replacing $s_0 (1-\eps)$ with $s_0$, so that
$$\P ^y(v_7(s_0) \geq u_0/(1-\eps))
\leq 
{c_{11}  \delta (1-\eps)\over  \eps^2 s_0 \sigma^2}. 
$$
This and (3.16) imply that
for small $\eps$ and $\delta$ we have
$$v_7(s_0) \leq u_0/(1-\eps) 
\leq 2 s_0 \delta (1+\eps)^2/\sigma^2\eqno(3.34)$$
with probability greater than or equal to
$1 - c_{11}  \delta (1-\eps)/(  \eps^2 s_0 \sigma^2)$.
The inequality (3.34) is analogous to (3.11)
and can be used in the same way as in the argument
between (3.19) and (3.20)
to prove a formula analogous to (3.20):
$$ v_4^\delta(a) = v_4(j s_0)
\leq 
{\delta \over \beta_1 -\beta_2}
(1- \exp (-2 a (\beta_1 -\beta_2)  (1+\eps)^3/\sigma^2) ),\eqno(3.35)
$$
where $\delta$ in $v^\delta_4(a)$ indicates the dependence
of $v^\delta_4(a)$ on $\delta$.
The above argument requires that we can use an estimate
analogous to (3.34) at every stage of the inductive 
procedure, i.e., at every stopping time
$v_4(m s_0)$ for $m=1,2, \dots , j-1$. All of these estimates
hold simultaneously with probability greater than
$$1 - (j-1)c_{11}  \delta (1-\eps)/(  \eps^2 s_0 \sigma^2).$$
This shows that the probability that (3.35) fails to hold
is smaller than
$$(j-1)c_{11}  \delta (1-\eps)/(  \eps^2 s_0 \sigma^2).$$
Now fix arbitrarily small $\eps_1>0$ 
and let $\delta_k = (1-\eps_1)^k$.
Let $A_k$ be the event in (3.35) with $\delta$ replaced by $\delta_k$, i.e.,
$$A_k = \left \{
v_4^{\delta_k}(a) 
\leq 
{\delta_k \over \beta_1 -\beta_2}
(1- \exp (-2 a (\beta_1 -\beta_2)  (1+\eps)^3/\sigma^2) )
\right\}.$$
We have
$$\sum_{k=0}^\infty
(j-1)c_{11}  \delta_k (1-\eps)/(  \eps^2 s_0 \sigma^2)
= \sum_{k=0}^\infty
(j-1)c_{11}  (1-\eps_1)^k (1-\eps)/(  \eps^2 s_0 \sigma^2)
<\infty,
$$
so only a finite number of events $A_k$ may fail to hold.
Consider an $\omega$ and $k_0$ such that all events
$A_k, k\geq k_0$, hold for this $\omega$. Suppose that
$\delta\in(0,\delta_{k_0})$. Then $\delta \in [\delta_{k_1-1}, \delta_{k_1}]$
for some $k_1\geq k_0$. Since $A_{k_1}$ holds,
we have
$$\eqalign{
v_4^{\delta}(a) &\leq v_4^{\delta_{k_1}}(a) 
\leq 
{\delta_{k_1} \over \beta_1 -\beta_2}
(1- \exp (-2 a (\beta_1 -\beta_2)  (1+\eps)^3/\sigma^2) )\cr
&\leq 
{\delta/(1-\eps_1) \over \beta_1 -\beta_2}
(1- \exp (-2 a (\beta_1 -\beta_2)  (1+\eps)^3/\sigma^2) ).
}$$
This inequality holds with probability one for all
sufficiently small $\delta>0$. Since $\eps>0$ and $\eps_1>0$ are
arbitrarily small, we see that a.s.,
$$\limsup_{\delta \to 0+} {v_4^{\delta}(a) \over \delta}
\leq {1 \over \beta_1 -\beta_2}
(1- \exp (-2 a (\beta_1 -\beta_2)  /\sigma^2) ).
$$
The same lower bound can be obtained for liminf in a completely
analogous way, so with probability one,
$$\lim_{\delta \to 0+} {v_4^{\delta}(a) \over \delta}
= 
{1- \exp (-2 a (\beta_1 -\beta_2)  /\sigma^2) \over \beta_1 -\beta_2}.
\eqno(3.36)
$$

Suppose $a>0$ and let $v(a) = \inf\{t\geq 0: L^x_t = a\}$. 
Fix some $x\in \R$ and consider $\delta>0$. We will first find
the right hand side derivative
$\left. {d^+\over dy} X^y_{v(a)} \right|_{y=x}$.
Let $T= \inf\{t: B_t = X^x_t\}$ and let
$U_1$ be the amount of time spent by $B_t$
between the graphs of $X^x_t$ and $X^{x+\delta}_t$
on the time interval $[0,T]$. We will write
$U_2$ to denote the amount of time spent by $B_t$
between the graphs of $X^x_t$ and $X^{x+\delta}_t$,
between times $T$ and $v(a)$.

If $x\geq 0$ then $U_1=0$. If $x<0$ then $U_1$
is not greater than the amount of time $U_3$ spent by $B_t$
between the lines $t\to x+ \beta_1 t$ and
$t \to x+\delta + \beta_1 t$, until the hitting time $T$.
Standard arguments show that for any arbitrarily small
$\eps >0$, we have $U_3/\delta^{2-\eps} \to 0$
as $\delta \to 0$, a.s.
Note that the distance between $X^{x+\delta}_t$ and $X^x_t$
decreases by $(\beta_1 -\beta_2)u$ on any interval
where the Brownian motion $B_t$ spends $u$ units between
these functions.
Hence,
$$X^{x+\delta}_{v(a)} - X^x_{v(a)} = \delta 
- (\beta_1 -\beta_2) (U_1 + U_2).$$
The random variable $U_2$ may be identified with $v^\delta_4(a)$, so
(3.36) gives for any fixed $a$, a.s.,
$$\eqalign{
\left. {d^+\over dy} X^y_{v(a)} \right|_{y=x}
&= \lim_{\delta \to 0} 
{X^{x+\delta}_{v(a)} - X^x_{v(a)} \over \delta} \cr
&= \lim_{\delta \to 0} 
{\delta 
- (\beta_1 -\beta_2) (U_1 + U_2) \over \delta} \cr
&= 1 - \lim_{\delta \to 0} 
{(\beta_1 -\beta_2) U_1  \over \delta} 
- \lim_{\delta \to 0} 
{(\beta_1 -\beta_2) v^\delta_4(a) \over \delta} \cr
&=1-0 -(1- \exp (-2 a (\beta_1 -\beta_2)  /\sigma^2)) \cr
&= \exp (-2 L^x_{v(a)} (\beta_1 -\beta_2)  /\sigma^2).}$$
The above holds simultaneously for all rational $a$, with
probability one. Since $t\to L^x_t$ 
and $t\to X^y_t - X^z_t$ are continuous monotone functions,
an elementary argument can be used to extend the last
formula to fixed times, i.e.,
$$\left. {d^+ \over dy} X^y_{t} \right|_{y=x}
= \exp (-2 L^x_{t} (\beta_1 -\beta_2)  /\sigma^2),\eqno(3.37)$$
simultaneously for all $t\geq 0$, a.s.

Fix some $t>0$. By Fubini's theorem, (3.37) holds
for almost all $x$, a.s.
We have $|X^y_t - X^z_t| \leq |y-z|$ for all $y$ and $z$.
Since the function $y \to X^y_t$ is Lipschitz, it has a derivative
almost everywhere and so for a fixed $t$, we may replace
the right hand derivative with the usual derivative in (3.37),
for almost all $x$. The function $x\to L^x_t$ is continuous,
so the derivative in (3.37) is equal almost everywhere
to a continuous function. This implies that the derivative is equal
to the function everywhere.
This proves the first assertion of the theorem.
The second one follows from the first one and from
the Lipschitz character of $y \to X^y_t$.
\qed

\bigskip
\noindent{\bf Remark 3.9}.
Suppose that $X^y_t$ are solutions to (1.1)
and assume that $\beta_1, \beta_2 >0$ and
$\beta_1 - \beta_2 >0$. Fix some $t>0$
and consider the function $y \to X^y_t$.
We will sketch an argument showing that $y \to X^y_t$ is $C^{1+\gamma}$
for every $\gamma < 1/2$, i.e., that the function
has a derivative which is H\"older continuous
with H\"older exponent $\gamma$.

Fix any $z \in \R$. With probability 1,
$B_t \ne X^z_t$, and with strictly positive probability,
there exists $\eps>0$ such that 
$B_s \ne X^y_s$ for all $y \in (z-\eps, z+\eps)$ and $s\geq t$.
It follows that if a local property holds
for the function $y \to L^y_\infty$ with probability 1,
it must hold for $y\to L^y_t$, with probability 1.
Since $y \to L^y_\infty$ is a diffusion, its paths
are H\"older continuous with exponent $\gamma $
for every $\gamma < 1/2$. It follows that the same
is true of $y\to L^y_t$. Theorem 3.8 now implies that
$y \to X^y_t$ is $C^{1+\gamma}$ for every $\gamma < 1/2$.
The same argument shows that $y \to X^y_t$ is not $C^{3/2}$.

\bigskip
\noindent{\bf 4. Time and direction of bifurcation}.
We will first address the question 
of the direction of bifurcation for the equation (1.3).
We will say that a positive bifurcation occurs
if for some $t_1$ we have $X_t > B_t$ for all
$t>t_1$. The definition of a negative bifurcation
is analogous. If $\beta_1$ and $\beta_2$ have the same
sign then it is easy to see that a bifurcation
will occur with probability one and its direction
will be the same as the sign of $\beta_k$'s.
If $\beta_1 >0$ and $\beta_2 < 0$ then there will
be no bifurcation. The next theorem deals with the only
remaining, non-trivial case.

\bigskip
\noindent{\bf Theorem 4.1}.
{\sl
Consider the equation (1.3) with $t_0 = x_0 = 0$.
Assume that $\beta_1 <0$ and $\beta_2 > 0$.
Let
$$\lambda_j =
{ 2 |\beta_j| ^{1/(\alpha_j + 1)} (\alpha_j + 1)^{\alpha_j/(\alpha_j + 1)}
\over
\sigma^{2/(\alpha_j + 1)}\Gamma(1/(\alpha_j +1)) },$$
for $j=1,2$.
The probability of a negative bifurcation is
equal to $\lambda_1/ (\lambda_1 + \lambda_2)$.
When $\alpha_1 = \alpha_2=0$, the formula simplifies to
$ |\beta_1| /(|\beta_1| + |\beta_2|)$.
}

\bigskip

Before proving Theorem 4.1 we present
a lemma which may have some interest of its own.

\bigskip
\noindent{\bf Lemma 4.2}.
{\sl Assume that $\beta_1 <0$ and $\beta_2 > 0$.
Consider a solution $X_t$ to (1.3) with
$t_0=x_0 =0$.
There exists $\gamma>0$, depending on $\alpha_1,\alpha_2,\beta_1,\beta_2$,
such that $X_t / t^{1/2+\gamma}$ converges in probability to $0$ as $t\to 0$.
}

\bigskip
\noindent{\bf Proof}.
Let us assume that $-1 < \alpha_1 \leq 0 \leq \alpha_2$.
The other cases may be treated in a similar way.
Let 
$U= \sup\{s \leq t: X_s - B_s \geq -t^{1/2} \}$.
For $s\in (U,t)$ we have $X_s - B_s < - t^{1/2}$ so
for such $s$, 
$|dX_s/ds| \leq \beta_1 t^{(1/2)\alpha_1}$.
It follows that 
$$|X_t - X_U | \leq -(t-U) \beta_1 t^{(1/2)\alpha_1}
\leq -\beta_1 t^{1+(1/2)\alpha_1},$$
and so
$$\eqalign{
B_t - X_t &\leq (B_t - B_U) + (B_U-X_U) + (X_U - X_t)\cr
&\leq \left( \max_{s\in(0,t)} B_s - \min_{s\in(0,t)} B_s\right)
+ t^{1/2} -\beta_1 t^{1+(1/2)\alpha_1}.}$$
This implies that
$$\eqalign{
\E &|B_t -X_t|^{\alpha_1}\bone_{\{ B_t - X_t >0 \}} \cr
&\leq \E 
\left|\left( \max_{s\in(0,t)} B_s - \min_{s\in(0,t)} B_s\right)
+ t^{1/2} -\beta_1 t^{1+(1/2)\alpha_1}\right|^{\alpha_1} \cr
&\leq 3^{\alpha_1}
\left[ \E  \left( \max_{s\in(0,t)} B_s - \min_{s\in(0,t)} B_s\right)^{\alpha_1}
+ t^{(1/2)\alpha_1} + |\beta_1|^{\alpha_1} t^{\alpha_1+(1/2)\alpha_1^2} \right ]\cr
&\leq c_1 t^{(1/2)\alpha_1} + c_2|\beta_1|^{\alpha_1} t^{\alpha_1+(1/2)\alpha_1^2}.
}$$
Recall from (2.2) that
$$X_t  = 
\int_{0}^t \left[
\beta_1  |X_s-B_s|^{\alpha_1} \bone_{\{X_s - B_s \leq 0\}} 
+ \beta_2  |X_s-B_s|^{\alpha_2}\bone_{\{X_s - B_s > 0\}}
\right] ds. $$
From this we have the following estimate
$$\eqalign{
\E  X_t &\geq 
\E  \int_{0}^t 
\beta_1  |X_s-B_s|^{\alpha_1} \bone_{\{X_s - B_s \leq 0\}} ds \cr
&= \int_{0}^t \E (
\beta_1  |X_s-B_s|^{\alpha_1} \bone_{\{X_s - B_s \leq 0\}}) ds \cr
&\geq \beta_1 \int_{0}^t (c_1 s^{(1/2)\alpha_1} 
+ c_2|\beta_1|^{\alpha_1} s^{\alpha_1+(1/2)\alpha_1^2}) ds \cr
& = \beta_1 \left(c_3 t^{1+ (1/2)\alpha_1} 
+ c_4 |\beta_1|^{\alpha_1} t ^{1+\alpha_1+(1/2)\alpha_1^2}\right).
}$$
Since $\alpha_1 > -1$, the exponents $1+ (1/2)\alpha_1$
and $1+\alpha_1+(1/2)\alpha_1^2$ are greater than $1/2$
and so for some $\gamma>0$ and every $c_5>0$,
$\liminf_{t\to 0} c_5 \E X_t/ t^{1/2+\gamma} \geq 0$.
It follows that 
$$\lim_{t\to 0} \P (X_t / t^{1/2+\gamma} < -c_6)=0, \eqno(4.1)$$ 
for every $c_6>0$.

Recall that $\alpha_2 \geq 0$. 
Since
$$X_t \leq 
\int_{0}^t  \beta_2  |X_s-B_s|^{\alpha_2}\bone_{\{X_s - B_s > 0\}}
 ds, $$
an elementary argument shows that for small $t$ we have
$X_t \leq  2\beta_2 t$ if $B_s \leq 1$ for all $s\in (0,t)$.
It is clear that $\P (\max_{s\in(0,t)} B_s >1)$
goes to $0$ as $t\to 0$ so
$$\lim_{t\to 0} \P (X_t / t > 2\beta_2)=0.$$ 
This and (4.1) prove the lemma.
\qed

\bigskip
\noindent{\bf Proof of Theorem 4.1}.
The assertion of the theorem
deals only with probabilities, so we can use
any solution to (1.3), as we have uniqueness in law
by Theorem 2.1. The same theorem shows that
a solution $X_t$ may be constructed so that
$(X_t,B_t)$ is a strong Markov process, and hence
we may apply excursion theory to it.
Recall the discussion at the beginning of Section 3.
The same analysis of excursion laws and the exit
system applies to the solutions of (1.3)
for arbitrary $\alpha_1, \alpha_2>-1$.
Let us briefly recall the facts that we will
need in our present argument.
Let $D= \{(b,x) \in \R^2: b = x \}$ and let
$(H^x, dL)$ be an exit system for the process of
excursions of $(B_t,X_t)$ from the set $D$.
The generic excursion may be denoted
$(e^1_t, e^2_t)$.
By the translation invariance
of the Brownian motion $B_t$ and the equation (1.3),
the distribution of $(e^1_t - x, e^2_t -x)$ under $H^x$
is the same for every $x \in \R$. Let this distribution
be called $H_1$. Let $H_{1+}$
denote the part of the measure $H_1$ which is supported
on excursions with $e^1_t > e^2_t$ and let $H_{1-}$
be the part supported on the set where $e^1_t < e^2_t$.
Let $H_{2+}$ be the distribution of 
$\{e^1_t - e^2_t, t \in (0,\nu)\}$ under $H_{1+}$
and let $H_{2-}$ have the same definition relative
to $H_{1-}$. We have,
up to a multiplicative constant,
$$H_{2+}(A) = \lim_{x \downarrow 0}
{1 \over |x|} Q^x_+(A),\eqno(4.2)$$
where $Q^x_+$ stands for the distribution
of the diffusion $Y_t$ with the same infinitesimal
variance as Brownian motion (i.e., $\sigma^2$)
but with drift $-\beta_1|Y_t|^{\alpha_1}$, killed
at the hitting time of $0$. We will normalize $H_{2+}$
as in (4.2). 
We will have to prove that the following formula
gives the correct normalization for $H_{2-}$,
$$H_{2-}(A) = \lim_{x \uparrow 0}
{1 \over |x|} Q^x_-(A).\eqno(4.3)$$
Here $Q^x_-$ denotes the distribution
of diffusion $Z_t$ with Brownian quadratic variation (name\-ly, $\sigma^2$)
and drift $-\beta_2|Z_t|^{\alpha_2}$, killed
at the hitting time of $0$.

The proof that (4.3) is the correct normalization for $H_{2-}$
can proceed exactly as the proof of Lemma 3.1 (iv),
thanks to Lemma 4.2. It only remains to find
and compare the formulae analogous to those for
$H_{2+}(F_h)$ and $H_{2-}(F_h)$.
Recall that $F_h$ is the event that the difference
between the maximum and the minimum
of an excursion exceeds $h$. 
The scale function for
a diffusion on $(0,\infty)$ with
infinitesimal drift $\mu(x) = -\beta_1 x ^{\alpha_1}$
and variance $\sigma^2$
is given by (see Karlin and Taylor (1981), p. 194),
$$\eqalignno{
S(x) &= \int_1^x \exp \left( 
- \int _0^ y { 2 \mu(z) \over \sigma^2} dz \right) dy 
= \int_1^x \exp \left( 
- \int _0^ y { -2 \beta_1 z ^{\alpha_1} \over \sigma^2} dz \right) dy \cr
&= \int_1^x \exp \left( 
 { 2 \beta_1 y ^{\alpha_1 +1} \over \sigma^2(\alpha_1 +1)}  \right) dy .&(4.4)
}$$
By (4.2),
$$\eqalign{H_{2+}(F_h) &= 
\lim_{x \downarrow 0}
{1\over x} Q^x_+ (T_h < T_0)
= \lim_{x \downarrow 0}
{1\over x} \cdot{ S(x) - S(0) \over S(h) - S(0) } \cr
&= \lim_{x \downarrow 0}
{1\over x} \cdot{\int_0^x \exp \left( 
{ 2 \beta_1 y ^{\alpha_1 +1} \over \sigma^2(\alpha_1 +1)}  \right) dy
\over
\int_0^h \exp \left( 
{ 2 \beta_1 y ^{\alpha_1 +1} \over \sigma^2(\alpha_1 +1)}  \right) dy
} 
= { 1 \over 
\int_0^h \exp \left( 
{ 2 \beta_1 y ^{\alpha_1 +1} \over \sigma^2(\alpha_1 +1)}  \right) dy
} .
}$$
If we use (4.3), we obtain in the same way
$$H_{2-}(F_h) =
{ 1 \over 
\int_0^h \exp \left( 
{ -2 \beta_2 y ^{\alpha_2 +1} \over \sigma^2(\alpha_2 +1)}  \right) dy
} ,
$$
which implies that 
$$\lim_{h\to 0} H_{2+}(F_h) / H_{2-}(F_h) = 1,$$
and this confirms that the normalization in (4.3) is correct.

The probability for the process $Y_t$ starting from
$\delta$ never to hit $0$ is equal to
$$\lim_{b\to \infty}
{S(\delta) - S(0) \over S(b) - S(0)}
=
{\int_0^\delta \exp \left( 
{ 2 \beta_1 y ^{\alpha_1 +1} \over \sigma^2(\alpha_1 +1)}  \right) dy
\over
\int_0^\infty \exp \left( 
{ 2 \beta_1 y ^{\alpha_1 +1} \over \sigma^2(\alpha_1 +1)}  \right) dy
}.
$$
It follows that $H_{2+}(F_\infty)$,
i.e., the measure given to positive excursions which do not
return to $0$ is given by
$$\eqalignno{
\lim_{\delta\to 0} &{1\over \delta} \cdot
{\int_0^\delta \exp \left( 
 { 2 \beta_1 y ^{\alpha_1 +1} \over \sigma^2(\alpha_1 +1)}  \right) dy
\over
\int_0^\infty \exp \left( 
 { 2 \beta_1 y ^{\alpha_1 +1} \over \sigma^2(\alpha_1 +1)}  \right) dy}
=
\left[ \int_0^\infty \exp \left( 
 { 2 \beta_1 y ^{\alpha_1 +1} \over \sigma^2(\alpha_1 +1)}  \right)
\right] ^{-1}\cr
&={ \left({- 2 \beta_1  \over \sigma^2(\alpha_1 +1)}\right)^{1/(\alpha_1 + 1)}
(\alpha_1 + 1) 
\over
\Gamma(1/(\alpha_1 +1)) } 
=
{ (-2 \beta_1) ^{1/(\alpha_1 + 1)} (\alpha_1 + 1)^{\alpha_1/(\alpha_1 + 1)}
\over
\sigma^{2/(\alpha_1 + 1)}\Gamma(1/(\alpha_1 +1)) } \df \lambda_1 .&(4.5)
}$$
An analogous formula holds for $\lambda_2 \df H_{2-}(F_\infty)$.
The processes of excursions on both sides are independent
so the probability of the negative bifurcation is the same
as the probability that the first arrival of an infinite excursion in the 
Poisson process on the negative side
comes before the analogous event on the other side.
The probability in question is the ratio of $\lambda_1$ and
$\lambda_1 + \lambda_2$. \qed

\bigskip
\noindent{\bf Remark 4.3}. Mike Harrison pointed out to us that
Theorem 4.1 may be proved without using excursion theory.
One can calculate the probability that 
the diffusion $X_t - B_t$
will go to infinity using an explicit formula
for the scale function of this diffusion.
The excursion theory approach has its advantages, though.
First, excursion theory seems to be the right tool
for the proof of Theorem 4.4 below. Second, the excursion
theory may be used to find the positive bifurcation 
probability when the vector process $(X_t, B_t)$
is Markov but $X_t-B_t$ is not. The solution of (1.2),
studied in Burdzy, Frankel and Pauzner (1998), is an example
of such a situation.

\bigskip
Let $T_*$ denote the bifurcation time, i.e.,
let $T_*$ be the supremum of $t$ with
$X_t =B_t$.

\bigskip
\noindent{\bf Theorem 4.4}.
{\sl
Consider the solution to (1.1) with
$t_0=x_0 =0$, $\beta_1 < 0$ and $\beta_2>0$. Then
$$\E  T_* = {\sigma^2 \over |2 \beta_1 \beta_2|}.$$
}

\bigskip
\noindent{\bf Proof}.
By Remark 3.2 (ii),
the distribution of the excursion law on excursions
with finite lifetime remains the same if we change $\beta$ to $-\beta$.
Hence, the formula (3.1) applies in the case
$\beta_1 < 0$ and $\beta_2>0$,
and we have
$$\E V_s = 
\left( {1 \over |\beta_1|} + {1 \over |\beta_2|}\right)s,$$
for $V_s = \inf \{t\geq 0: L_t \geq s\}$.
By excursion theory, the infinite excursion of $B_t$ from $X^0_t$
occurs independently from finite excursions in the Poisson point
process of excursions, so
the expected bifurcation time is equal to
$$\E  T_* = \int_0^\infty \lambda e^{-\lambda s} \E  V_s ds,$$
where $\lambda$ is the intensity of the Poisson process arrival for 
infinite excursions.
We have
$$\lambda = {2 (|\beta_1| +|\beta_2|) \over \sigma^2},$$
from (4.5), taking into account infinite excursions on both sides.
It follows that 
$$\eqalign{
\E  T_* & =
\int_0^\infty \lambda e^{-\lambda s} \E  V_s ds
=\int_0^\infty \lambda e^{-\lambda s} 
\left( {1 \over |\beta_1|} + {1 \over |\beta_2|}\right)s ds
={ 1 \over \lambda } \left( {1 \over |\beta_1|} + {1 \over |\beta_2|}\right)\cr
&={\sigma^2 \over 2 (|\beta_1| +|\beta_2|)}
\left( {1 \over |\beta_1|} + {1 \over |\beta_2|}\right)
={\sigma^2 \over |2 \beta_1 \beta_2|}. \squareinf}$$

\bigskip
\noindent{\bf Remark 4.5}. (i) It is also the case that
$$\E  T_* = {\sigma^2 (\beta_1 + \beta_2)\over 2 \beta_1 \beta_2^2}$$
if $\beta_1,\beta_2 >0$.
We leave the proof to the reader.

(ii) A similar result can be obtained for any values
of $\alpha_1, \alpha_2>-1$ but the formula
does not seem to have a compact form, so we only sketch
how it can be obtained. The proof of Theorem 4.4
needs two ingredients. One of them is the
expected amount of local time before the 
infinite excursion occurs. This is equal to the
expectation of the
minimum of two independent exponential 
random variables whose expected values
are inverses of the quantity in (4.5) (for
$(\alpha_1,\beta_1)$ and $(\alpha_2,\beta_2)$).

The second ingredient is the expectation of the
inverse local time at $s$,
for the process with finite excursions only.
This is equal to $s$ times the expected lifetime
of a finite excursion under the excursion law.
Here is how we can calculate this quantity.
For arbitrary $\alpha_1>-1$ we write as in (4.4),
$$\eqalignno{
s(x) & = \exp \left( -
 { 2 \beta_1 x ^{\alpha_1 +1} \over \sigma^2(\alpha_1 +1)}  \right), \cr
S(x) 
&= \int_1^x \exp \left( -
 { 2 \beta_1 y ^{\alpha_1 +1} \over \sigma^2(\alpha_1 +1)}  \right) dy .
}$$
For $0<x<y<\infty$, the Green function for Brownian
motion $Y_t$ with drift $-\beta_1$ (the negative sign
is due to restriction of the excursion law to finite
excursions) is given by
(see Remark 3.3 on p.~198 of Karlin and Taylor (1981)),
$$G(x,y)
= { 2 [S(x) - S(0) ] [S(\infty) - S(y)]
\over
\sigma^2 s(y) [S(\infty) - S(0)]}.$$
Hence, the expected lifetime of an excursion
is equal to 
$$\eqalign{
\lim_{x\downarrow 0} {1 \over x}
\int_0^\infty G(x,y)dy
& = \lim_{x\downarrow 0} {1 \over x}
\int_0^\infty { 2 [S(x) - S(0) ] [S(\infty) - S(y)]
\over
\sigma^2 s(y) [S(\infty) - S(0)]}dy \cr
&= { 2 \over
\sigma^2  [S(\infty) - S(0)]}
\int_0^\infty {[S(\infty) - S(y)]
\over s(y)}dy \cr
&= { 2 \over
\sigma^2  \int_0^\infty \exp \left( -
 { 2 \beta_1 z ^{\alpha_1 +1} \over \sigma^2(\alpha_1 +1)}  \right) dz}
\int_0^\infty {\int_y^\infty \exp \left( -
 { 2 \beta_1 z ^{\alpha_1 +1} \over \sigma^2(\alpha_1 +1)}  \right) dz
\over \exp \left( -
 { 2 \beta_1 y ^{\alpha_1 +1} \over \sigma^2(\alpha_1 +1)}  \right)}dy.
}$$
Adding this to the analogous quantity for $\alpha_2$
gives the expected lifetime
of a finite excursion under the excursion law.

\bigskip
\noindent{\bf 5. Lipschitz approximations}.
In this section we will address the question of how well
a Lipschitz function can approximate a Brownian path.
Our analysis will be based on the fact, proved in Lemma 5.2 below,
that a certain solution $X^*_t$ to (1.1) may be looked
upon as a Lipschitz approximation to $B_t$.

Our first lemma consists of two elementary observations
which are designed to help develop the mental
picture of the solutions $X^x_t$ of (1.1), in preparation
for Lemma 5.2.

\bigskip
\noindent{\bf Lemma 5.1} {\sl 
Let $X^x_t$ denote the solution of (1.1) with $X^x_0=x$.

\noindent (i) If $x < y$ then $X^x_t < X^y_t$ for all $t\in\R$, a.s.

\noindent (ii) For a fixed $t$, the function
$x \to X^x_t$ is continuous a.s.
}

\bigskip
\noindent{\bf Proof}.
(i) Suppose that we have $x < y$ and $X^x_s = X^y_s$
for some $s\in \R$. The two functions $X^x_t$ and $X^y_t$
are not identical
since $X^x_0 = x \ne y = X^y_0$, but they are both solutions
to (1.1) with $t_0=s$ and $x_0 = X^x_s$. This contradicts
the uniqueness of solutions to (1.1).

(ii) Consider any sequence $x_n$ converging 
monotonically to $x_\infty \in \R$.
By (i), the sequence $X^{x_n}_t$ is also monotone
in $n$, and by the Lipschitz property it must
converge to a limit $X^\infty_t$. The Lipschitz
property of the $X^{x_n}_t$'s implies that
of $X^\infty_t$. 
We can show that $X^\infty_t$ is a solution to (1.1)
using the same argument as in the proof of existence
in Theorem 2.2 for (1.1).
We must have $X^\infty_0 = X^{x_\infty}_0$,
so the uniqueness of the solutions to (1.1)
implies that $X^\infty_t = X^{x_\infty}_t$
for all $t$, a.s. We have shown that
$x_n \to x_\infty$ implies $X^{x_n}_t \to X^{x_\infty}_t$.
This completes the proof.
\qed

\bigskip
\noindent{\bf Lemma 5.2}. {\sl
Assume that $\beta_1 < 0 < \beta_2$.
For almost every $\omega$ there exists a  unique 
$\ol x=\ol x(\omega)$ such that the solution $X^{\ol x}_0$
to (1.1), that is, the solution satisfying $X^{\ol x}_0(\omega) = \ol x(\omega)$,
has the property that there exist arbitrarily large
$t$ with $X^{\ol x}_t(\omega) = B_t(\omega)$.
}
\bigskip

It is easy to see that if $\beta_1 < 0 < \beta_2$
then with probability 1, all
solutions $X^x_t$ have the property that there exist
arbitrarily small $t>-\infty$ such that $X^x_t = B_t$.
Lemma 2.13 shows that a result analogous
to Lemma 5.2 holds when $\beta_2 < 0 < \beta_1$,
and we require that the solution intersects
the Brownian path for arbitrarily small $t>-\infty$.

\bigskip
\noindent{\bf Proof}. We will first prove the existence.
The law of the iterated logarithm easily implies
that for some random $x>0$, the functions
$t \to x + \beta_2 t$ and $t \to -x + \beta_1 t$
stay above and below the trajectory of $B_t$, for $t\geq 0$, resp.
This shows that there exist both large
and small (random) $x$ such that $X^x_t$ does not intersect
the trajectory of $B_t$ for $t>0$.

Let $A$ be the set of all $x$ such that 
$X^x_t > B_t$ for all $t$ greater than some $t_1= t_1(x)$.
By Lemma 5.1 (i) and the above remarks, the set $A$ is 
a non-empty semi-infinite interval. We will show
that it is open.
Consider an $x$ such that 
$X^x_t > B_t$ for all $t$ greater than some $t_1$.
Then $X^x_t = X^x_{t_1} + \beta_2 (t-t_1)$
for $t > t_1$. Let $c_1 = \inf \{X^x_t - B_t: t > t_1+1\}$
and note that $c_1 >0$, by the continuity of $X^x_t - B_t$.
By Lemma 5.1 (ii),
the function $y \to X^y _{t_1 +1}$ is continuous
so we can find $\eps >0 $ such that 
$X^y_{t_1+1} > X^x _{t_1 +1} - c_1/2$
for all $y > x - \eps$. It follows easily
that for such $y$, we have
$X^y_t = X^x_{t_1+1} + \beta_2 (t-t_1-1)$
and so 
$X^y_t > B_t$ for $t > t_1 +1$.
This proves that $A$ is open.
The same is true of the set $A'$ of $x$'s with the property that
$X^x_t < B_t$ for all $t$ greater than some $t_1= t_1(x)$.
Hence, $(A\cup A')^c$ is non-empty and so
we must have at least one $x$ for which
$X^x_t =B_t$ for arbitrarily large $t$.

We turn to the proof of uniqueness.
Suppose that with positive probability there exist $x_1 <x_2$,
such that both trajectories $X^{x_1} _t$ and $X^{x_2}_t$
intersect $B_t$ for arbitrarily large times $t$.
Then we can find $\delta>0$ and $p>0$ such that with probability
greater than $p$, there exist $x_1$ and $x_2$ with
$x_2 > x_1 + \delta$ and such that $X^{x_1} _t$ and $X^{x_2}_t$
intersect $B_t$ for arbitrarily large times $t$.
We will show that this assumption leads to a contradiction.

Fix some $\gamma \in (1/2,1)$. By the law of the iterated logarithm, we can
find a large $t_1$ with the following property.
For every $x$, if $|X^x_{t_1}| \geq t_1^\gamma$
then $X^1_t \ne B_t$ for all $t>t_1$,
with probability greater than $1 - p/8$,.

Consider solutions $\wh X_t^{y_1} $ 
and $\wh X_t^{y_2} $ to (1.1)
with $\wh X_{t_1}^{y_1} = y_1 \df t_1^\gamma$
and $\wh X_{t_1}^{y_2} = y_2 \df -t_1^\gamma$. 
We enlarge $t_1$, if necessary,
so that the event 
$\{| B_{t_1/2}| \geq t_1^\gamma\} 
\cup \{| B_{t_1}| \geq t_1^\gamma\} $
has a probability smaller than $p/8$. If the event
$\{| B_{t_1/2}| \geq t_1^\gamma\} 
\cup \{| B_{t_1}| \geq t_1^\gamma\} $
does not occur and $t_1$ is sufficiently large,
then both processes $\wh X_t^{y_1} $ 
and $\wh X_t^{y_2} $ must intersect the trajectory of $\wt B_t$
between $t_1/2$ and $t_1$.
Let $T = \sup\{ t< t_1: \wh X_t^{y_1} = B_t\}$.
Then the process $\{Y_t = B_{T-t} - \wh X^{y_1}_{T-t}, t\geq 0\}$
is a Brownian motion with drift $\beta_2$ if $Y_t <0$
and $\beta_1$ if $Y_t >0$. Note that the distribution
of the process $Y_t$ does not depend on $t_1$.
We can apply Proposition 3.3 to the local time
$L^Y_t$ of $Y_t$ at $0$, to see that
$L^Y_t/t \to (1/|\beta_1| + 1/|\beta_2|)^{-1} \df \lambda$, 
as $t\to \infty$, a.s.
Enlarge $t_1$ again, if needed, so that 
$L^Y_{t_1/2} > \lambda t_1/4$ with probability exceeding
$1-p/8$. Let $\eta$ be the expected number
of positive excursions of $Y$ whose height
does not exceed $\delta$, whose duration exceeds $1$,
and which start at a time $t$ with $L^Y_t \leq 1$.
Then the total number of such excursions
which start at times $t$ with $L^Y_t < \lambda t_1/4$
has a Poisson distribution with mean
$\eta \lambda t_1/4$. We make $t_1$ large enough
so that with probability greater than
$1-p/8$, the total number of such excursions
which start at times $t$ with $L^Y_t < \lambda t_1/4$,
is greater than $\eta \lambda t_1/8$.
Collecting the above facts, we see
that with probability greater than $1-3p/8$,
we have all of the following:
(i) $T \in (t_1/2, t_1)$, (ii) the local time
for the process $\wh X^{y_1}_t - B_t$ at $0$
accumulated between times $0$ and $t_1$ exceeds
$\lambda t_1/4$, and (iii) the number of negative excursions
of $\wh X^{y_1}_t - B_t$ whose absolute height is less than $\delta$,
the duration is greater than $1$, and which lie
within interval $(0,t_1)$, is greater than $\eta \lambda t_1/8$.
Note that $dX^{y_2}_t /dt - dX^{y_1}_t /dt = \beta_2 - \beta_1$
for any $t$ within such an excursion provided
$\wh X^{y_1}_t > \wh X^{y_2}_t + \delta$.
The last observation shows that if 
$\wh X^{y_1}_0 > \wh X^{y_2}_0 + \delta$ then
$\wh X^{y_1}_{t_1} > \wh X^{y_2}_{t_1} + \delta
+ 2 (\beta_2-\beta_1) \eta \lambda t_1/8$. If $t_1$ is sufficiently
large the last inequality cannot hold because we would have
$t_1^\gamma > -t_1^\gamma + \delta
+ 2 (\beta_2-\beta_1) \eta \lambda t_1/8$. We conclude that
with probability greater than $1-3p/8$, we have
$\wh X^{y_1}_0 - \wh X^{y_2}_0 < \delta$.

We reformulate the last statement in terms of
$x_1$ and $x_2$. Using Lemma 5.1 (i), we see
that the probability that there exist
$x_1$ and $x_2$ with $x_2 >x_1 +\delta$,
$|X^{x_1}_{t_1}| \leq t_1^\gamma$ and
$|X^{x_2}_{t_1}| \leq t_1^\gamma$ is less than $3p/8$.
An earlier argument showed that the probability that
$x_2 >x_1 +\delta$ and $|X^{x_1}_{t_1}| \geq t_1^\gamma$
or $|X^{x_2}_{t_1}| \geq t_1^\gamma$ is less than $2p/8$.
We conclude that the probability of $x_2 >x_1 +\delta$
is bounded by $5p/8$, which contradicts our assumption.
\qed

\bigskip
Consider equation (1.1) with $-\beta_1=\beta_2=\beta>0$. Let $X_t^*$ denote
the solution of (1.1) constructed in Lemma 5.2. That is, $X_t^*=X_t^{\ol x}$.

\ms
\noindent{\bf Lemma 5.3}. {\sl We have with probability 1,
$$
\limsup _{t\to -\infty} { X^*_t - B_t \over \log t}
= \limsup _{t\to \infty} { X^*_t - B_t \over \log t}
= \limsup _{t\to -\infty} {B_t - X^*_t  \over \log t}
= \limsup _{t\to \infty} { B_t - X^*_t  \over \log t}
\geq {\sigma^2 \over  2\beta}. $$
}

\ms
Note the lim inf's as $t\to \infty$ are zero as $X_t^*$ crosses $B_t$ for arbitrarily
large $t$.
\bigskip
\noindent{\bf Proof}.
Let $\wt X_t$ be a solution to (1.1) 
with $-\beta_1 = \beta_2 = -\beta$ and 
let $Y_t = \wt B_t - \wt X_t$.
The process $Y_t$ is a diffusion which spends zero time on the real
axis, which behaves like Brownian motion with drift $\beta$
when $Y_t <0$, and it is a Brownian motion with drift
$-\beta$ when $Y_t >0$.  
By Karlin and Taylor (1981), Chapter 15.5, (5.34),
the process $Y_t$
has a stationary probability distribution with a density 
$$\psi(y) = {\beta \over  \sigma^2}
\exp \left ( - {2 \beta |y| \over \sigma^2 }\right ).$$
Let $\{\wh Y_t, t\in \R\}$ be the process which has density
$\psi(y)$ for every fixed $t$, and which has the transition
probabilities of $Y_t$. Let 
$$ \wh X_t = \int_0^t \sign(\wh Y_s) \beta ds$$
and
$$\wh B_t = \wh Y_t - \wh Y_0 + \wh X_t.$$
It is easy to check that $\wh B_t$ is a Brownian motion
with $\wh B_0 =0$, and that $\wh X_t$ solves (1.1)
with $\wt B_t$ replaced by $\wh B_t$
and $\beta_1 = -\beta_2 =\beta $. Moreover,
$\wh X_t$ has the property that $\wh X_t = \wh B_t$
for infinitely many arbitrarily large negative  and arbitrarily large positive $t$.
If we now time-reverse $\wh B_t$ and $\wh X_t$,
we will obtain a Brownian motion and a corresponding
solution to (1.1) which satisfies the defining properties of $X^*_t$.
Hence, we may construct $B_t$ and the corresponding
process $X^*_t$ by letting $B_t = \wh B_{-t}$ and
$X^*_t = \wh X_{-t}$.

The scale function $S(y)$ for Brownian motion 
with drift $-\beta$ is given by 
$S(y) = \exp( 2 \beta y/\sigma^2)$
(Karlin and Taylor (1981) Chapter 15.4).
Let $T_a$ be the hitting time of $a$ by the process $Y$.
The mass $H(F_h)$ given by the excursion law for the process
$\wh Y_t$ to positive excursions whose height
exceeds $h$ is equal to 
$$\eqalign{
\lim_{\eps \downarrow 0}
{1\over \eps} \P ^\eps (T_h < T_0)
&= \lim_{\eps \downarrow 0}
{1\over \eps} \cdot{ S(\eps) - S(0) \over S(h) - S(0) } 
= \lim_{\eps \downarrow 0}
{1\over \eps} \cdot{  \exp( 2 \beta \eps /\sigma^2) -1 \over
\exp( 2 \beta h/\sigma^2) -1 } \cr
& = {2 \beta \over \sigma^2} \cdot
{ 1 \over \exp( 2 \beta h/\sigma^2) -1 }.
}$$

Fix some small $\eps >0$ and let 
$h_k = k \log 2 \cdot(1- \eps) \sigma^2/(2 \beta)$.
Let $L_t$ denote the local time of $\wh Y_t$ at $0$
with $L_0=0$,
and let $A_k$ denote the event that there exists
a positive excursion of $\wh Y_t$ whose
height exceeds $h_k$, and which starts at a time
$t$ such that $2^k \leq L_t < 2^{k+1}$.
The probability
of $A_k$ is the probability that a Poisson random
variable with mean $\lambda_k =2^k H(F_{h_k})$ takes a non-zero value.
Thus, $\P (A_k^c) =  e^{-\lambda_k}.$
For large $k$,
$$ \eqalign{
\lambda_k &= 2^k \cdot{2 \beta \over \sigma^2} \cdot
{ 1 \over \exp( 2 \beta h_k/\sigma^2) -1 }
\geq 2^k \cdot{2 \beta \over \sigma^2} 
\exp(- 2 \beta h_k/\sigma^2) \cr
&= 2^k \cdot{2 \beta \over \sigma^2} 
\exp \left( - 2 {\beta \over \sigma^2} \cdot
{ k \log 2 \cdot (1- \eps) \sigma^2 \over 2 \beta } \right)
= 2^k \cdot{2 \beta \over \sigma^2} \cdot
 2^{- k (1-\eps)}
= {2 \beta \over \sigma^2} \cdot
 2^{ k \eps}.
}$$
This implies that $\sum_k \P (A_k^c) = \sum_k e^{-\lambda_k} < \infty$.
By the Borel-Cantelli Lemma, only a finite number
of the events $A^c_k$ occur.
Hence,
$$\eqalign{
\limsup_{t\to\infty} { \wh Y_t \over \log L_t}
&\geq \limsup_{k \to \infty} \ \sup  
\left\{{ \wh Y_t \over \log L_t} : L_t \in [2^{k},2^{k+1}] \right\}\cr
&\geq \limsup_{k \to \infty} \ \sup  
\left\{{ \wh Y_t \over \log 2^{k+1}} : L_t \in [2^{k},2^{k+1}] \right\}\cr
& \geq \limsup_{k \to \infty}
{  h_k  \over (k+1) \log 2} \cr
& = \limsup_{k \to \infty}
{ k \log 2 \cdot(1- \eps) \sigma^2 \over 2 \beta (k+1) \log 2 }  \cr
&= { (1-\eps) \sigma^2  \over  2\beta  } .
}$$
Since $\eps>0$ is arbitrarily small
and, by Proposition 3.3, $\lim_{t\to \infty} L_t /t = \beta/2$, a.s., 
we obtain, with probability 1,
$$\limsup_{t\to\infty}  {\wh Y_t \over\log t} = 
{  \sigma^2  \over  2\beta  }.$$
A similar argument yields,
$$-\liminf_{t\to\infty}  {\wh Y_t \over\log t} = 
\limsup_{t\to-\infty}  {\wh Y_t \over\log t} =
-\liminf_{t\to\infty}  {\wh Y_t \over\log t} =
{  \sigma^2  \over  2\beta  }.$$
Recall from the first part of the proof
that $\wh Y_t = B_{-t} - X^*_{-t} - \wh Y_0$.
This combined with the results for $\wh Y_t$ implies the proposition.
\qed
\bigskip

The function $t \to a + \beta |t|$ is Lipschitz
with constant $\beta$. For some random $a$, this
function is greater than $B_t$ for every $t$, by the law of the iterated logarithm.
Since the infimum of an arbitrary family of Lipschitz
functions with constant $\beta$ is again
a Lipschitz function with constant $\beta$,
there exists a smallest Lipschitz function $Z^+_t$
with constant $\beta$ with the property that
$Z^+_t \geq B_t$ for all $t$. Let $Z^-_t$
be the largest Lipschitz function 
with constant $\beta$ such that
$Z^-_t \leq B_t$ for all $t$. Note that
$Z^+_t$ and $Z^-_t$ are not measurable
with respect to $\sigma\{B_s, s\leq t\}$.

\bigskip
\noindent{\bf Lemma 5.4}. {\sl Assume that $B_0=0$. We have with probability 1,
$$
\limsup _{t\to -\infty} { Z^+_t - Z^-_t \over \log t}
= \limsup _{t\to \infty} { Z^+_t - Z^-_t \over \log t}
\leq {\sigma^2 \over  2\beta}. $$
}

\bigskip
\noindent{\bf Proof}.
Consider $a_1,a_2>0$. Let
$$\eqalign{
A_{++} &= \{ \exists t>0: B_t = a_1 + \beta t\}, \qquad
A_{+-} = \{ \exists t>0: B_t = -a_2 - \beta t\}, \cr
A_{-+} &= \{ \exists t<0: B_t = a_1 - \beta t\}, \qquad
A_{--} = \{ \exists t<0: B_t = -a_2 + \beta t\}.
}$$  
The probability that $B_t$ ever hits the line
$t\to a_1 + \beta t$ is equal to
$\exp(- 2 a_1 \beta/\sigma^2)$ (Karlin and Taylor (1975), p.~362).
The probability that $B_t$ crosses
the line $a_1 + \beta t$ at some $t_1>0$
and then crosses the line $-a_2 - \beta t$
for some $t>t_1$ is bounded by 
$\exp( - 2 a_1 \beta /\sigma^2)\exp( - 2 a_2 \beta /\sigma^2)$,
by the strong Markov property applied at $t_1$.
The probability of crossing
first $-a_2 - \beta t$ and then
$a_1 + \beta t$ is bounded by the same quantity.
Hence, 
$$\P (A_{++}\cap A_{+-}) \leq
2 \exp( - 2 (a_1 + a_2) \beta /\sigma^2).$$
The same estimate holds for $\P (A_{--}\cap A_{-+})$,
by symmetry. We obtain
$$\P (A_{++}\cap A_{--}) = \P (A_{-+}\cap A_{+-})
= \exp( - 2 (a_1 + a_2) \beta /\sigma^2),$$
from the independence of the processes $\{B_t, t\geq 0\}$
and $\{B_t, t\leq 0\}$. It follows that
$$\eqalign{
\P &(Z^+_0 -B_0 \geq a_1, B_0 - Z^-_0 \geq a_2) =
\P (Z^+_0 \geq a_1, Z^-_0 \leq -a_2) \cr
&\leq 
\P ([A_{++}\cap A_{+-}]\cup 
[A_{--}\cap A_{-+}]\cup
[A_{++}\cap A_{--}] \cup [A_{-+}\cap A_{+-}])\cr
&\leq 8 \exp( - 2 (a_1 + a_2) \beta /\sigma^2).}$$

Choose $\eps\in (0,1)$. Let $m>8$ be an integer 
large enough that $(m-1)/(m(1-\eps))>1$.
We have for any $y>0$,
$$\eqalign{
\P ( Z^+_{0} - Z^-_{0} \geq y)
&\leq \sum _{j=0}^m
\P (Z^+_0 - B_0 \geq j y/m , B_0 - Z^-_0 \geq (m-j-1) y /m)\cr
& \leq \sum _{j=0}^m
8 \exp\left( - 2 \left({j y\over m} + {(m-j-1) y \over m}\right) 
{\beta \over \sigma^2}\right) \cr
&\leq 9 m\exp\left( - 2 (m-1) y \beta \over m\sigma^2\right).
}$$

Fix some large $b<\infty$.
Consider
an integer $k>0$. Let $n$ be the integer part of
$${ 2 \beta (1-\eps) b 2^k \over \sigma^2 k \log 2},$$
and let $x_k =   b 2^k / n $, and $t^k_j = j  2^k /n$.
We have,
$$\eqalign{
\P (Z^+_{t^k_j} - Z^-_{t^k_j} \geq  x_k) 
&= \P (Z^+_0 - Z^-_0 \geq x_k)  \cr
& \leq 9m  \exp\left( {- 2(m-1) x_k \beta \over m\sigma^2}\right) \cr
&= 9m \exp \left ( - { (m-1)2 \beta b 2^k  
\over mn\sigma^2 }\right) \cr
& \leq 9m \exp \left ( - { (m-1)2 \beta b 2^k \sigma^2 k \log 2 
\over m2 \beta (1-\eps) b 2^k\sigma^2 }\right) \cr
& = 9m \exp \left( - k (m-1)\log 2 \over m (1-\eps)\right) 
= 9m \cdot 2^{-k(m-1)/(m(1-\eps))}.
}$$
For some $c_1< \infty$, using $(m-1)/(m(1-\eps)) >1$,
we obtain,
$$\eqalign{
\sum_{k=1}^\infty  & \sum_{0 \leq t^k_j \leq 2^k}
\P (Z^+_{t^k_j} - Z^-_{t^k_j} \geq  x_k) 
\leq c_1\sum_{k=1}^\infty 2 n \cdot 9m
\cdot 2^{-k(m-1)/(m(1-\eps))} \cr
&\leq c_1\sum_{k=1}^\infty
{ 2 \beta (1-\eps) b 2^k \over \sigma^2 k \log 2}\cdot
36 m \cdot 2^{-k(m-1)/(m(1-\eps))} < \infty.
}$$
By the Borel-Cantelli lemma, for all sufficiently
large $k$ and all $t^k_j \in[0, 2^k]$, we have
$Z^+_{t^k_j} - Z^-_{t^k_j} \leq  x_k$. 
If $Z^+_{t^k_j} - Z^-_{t^k_j} \leq   x_k$ then
for $t \in [ (t^k_j + t^k_{j-1})/2, (t^k_j + t^k_{j+1})/2]$,
$$Z^+_{t} - Z^-_{t} \leq   x_k + 2 \beta |t-t_j|
\leq  b 2^k / n + \beta  2^k /n
= (b 2^k /n) ( 1 + \beta/b)
=  x_k ( 1 + \beta/b).
$$
This implies that for large $k$, we have for all
$t\in[0, 2^k]$,
$$Z^+_{t} - Z^-_{t} \leq   x_k ( 1 + \beta/b).$$
We obtain
$$\eqalign{
\limsup_{t\to\infty} { Z^+_{t} - Z^-_{t} \over \log t}
&\leq \limsup_{k \to \infty} \sup _{t \in [2^{k-1},2^k]} 
{ Z^+_{t} - Z^-_{t} \over \log t} \cr
&\leq \limsup_{k \to \infty} \sup _{t \in [2^{k-1},2^k]} 
{ Z^+_{t} - Z^-_{t} \over \log 2^{k-1}} \cr
& \leq \limsup_{k \to \infty}
{  x_k ( 1 + \beta/b) \over (k-1) \log 2} \cr
& \leq \limsup_{k \to \infty}
{  (b 2^k/n) ( 1 + \beta/b) \over (k-1) \log 2} \cr
& \leq \limsup_{k \to \infty}
{  b 2^k \sigma^2 k \log 2 ( 1 + \beta/b) 
\over 2 \beta (1-\eps) b 2^k (k-1) \log 2} \cr
&= { \sigma^2  ( 1 + \beta/b) 
\over  2\beta (1-\eps) } .
}$$
Since $\eps$ may be chosen arbitrarily small
and $b$ may be chosen arbitrarily large,
with probability 1,
$$\limsup_{t\to\infty} { Z^+_{t} - Z^-_{t} \over \log t}
\leq { \sigma^2   \over  2\beta  }.$$
The result for $t\to -\infty$ follows by symmetry.
\qed
\bigskip

\noindent{\bf Theorem 5.5}.
{\sl (i) With probability 1,
$$\eqalign{
\limsup _{t\to -\infty} &{ X^*_t - B_t \over \log t}
= \limsup _{t\to \infty} { X^*_t - B_t \over \log t} \cr
&= \limsup _{t\to -\infty} {B_t - X^*_t  \over \log t}
= \limsup _{t\to \infty} { B_t - X^*_t  \over \log t} \cr
&= \limsup _{t\to -\infty} { Z^+_t - Z^-_t \over \log t}
= \limsup _{t\to \infty} { Z^+_t - Z^-_t \over \log t}
= {\sigma^2 \over  2\beta}. }$$

\noindent (ii) $\E |B_t - X^*_t| = {1 \over 2} \cdot {\sigma^2 / \beta}$,
for every $t\in\R$.

\noindent (iii) $\E  (Z^+_t - B_t) = \E  (B_t -Z^-_t) 
= {3 \over 4} \cdot\sigma^2 / \beta$, for every $t\in\R$.
}

\bigskip
Theorem 5.5 (i) shows, in a sense, that $Z^+$ and $Z^-$ are
as good Lipschitz approximations to $B_t$
as $X^*$. However, the comparison comes out differently when
we look at the averages presented in (ii) and (iii).

\bigskip
\noindent{\bf Proof}.
It is elementary to check that we always have
$Z^-_t \leq X^*_t \leq Z^+_t$. This and
Lemmas 5.3 and 5.4 yield (i).

Recall the stationary
density $\psi(y)$ for $Y_t$ from the proof
of Lemma 5.3. This is the same as
the density for the distribution of $B_t - X^*_t$.
Hence 
$$\E |B_t - X^*_t| = 
\int_{-\infty}^\infty |y|
{\beta \over  \sigma^2}
\exp \left ( - {2 \beta |y| \over \sigma^2 }\right ) dy =
{1 \over 2} \cdot {\sigma^2 \over \beta},$$
which proves (ii).

For $a>0$,
the probability that $B_t$ crosses
the line $a + \beta t$ for some $t>0$
is equal to $\exp( - 2 a \beta /\sigma^2)$
(Karlin and Taylor (1975) p. 362).
This is the same as the probability
of crossing the line $a - \beta t$ for some
$t<0$. The probability that none of these events
happen is $[1 - \exp( - 2 a \beta /\sigma^2)]^2$,
and so
$$\P (Z^+_0 < a) = [1 - \exp( - 2 a \beta /\sigma^2)]^2.$$
This yields 
$$\E  Z^+_0 = {3 \sigma^2 \over 4 \beta}.$$
We similarly have $\E  Z^-_0 = -3 \sigma^2 /( 4 \beta)$,
and by translation invariance, for every $t$,
$$\E  (Z^+_t - B_t) = \E  (B_t -Z^-_t) = {3 \sigma^2 \over 4 \beta}.
\squareinf$$

\bigskip
If we let $\alpha_1 = \alpha_2 = 1$
and choose suitable $\beta_1$ and $\beta_2$
in (1.3), then $Y_t = X_t - B_t$ is an
Ornstein-Uhlenbeck process. Results for such
a process, closely
related to Theorem 5.5 (i), can be found in the paper of
Darling and Erd\"os (1956).

\bigskip
\noindent{\bf Corollary 5.6}.
{\sl
For any random Lipschitz function $g(t)$ with constant
$\beta$ we have with probability one
$$\limsup _{t\to \infty} { g(t) - B_t \over \log t}
\geq  {\sigma^2 /(  4\beta)}.$$
}

\bigskip
\noindent{\bf Proof}.
Suppose that $X^*_t - B_t = a$
for some $t$ and $a>0$. Let $s$ be the largest
time less than $t$ such that $B_s = X^*_s$.
Then we see that the quantity 
$\sup_{s \leq u \leq t} |g(u) - B_u|$
cannot be smaller than $a/2$ for any Lipschitz function
$g(u)$ with constant $\beta$, by comparing $g(u)$ with
the function $u \to  B_s + a/2 + (u-s)\beta$.
Since
$\limsup _{t\to \infty} { (X^*_t - B_t )/ \log t}
= {\sigma^2 /(  2\beta)}$,
for any Lipschitz function $g(t)$ with constant
$\beta$ we must have
$$\limsup _{t\to \infty} { g(t) - B_t \over \log t}
\geq  {\sigma^2 /(  4\beta)}.\squareinf$$

\bigskip

Corollary 5.6 sheds some new light on
an old problem about strong approximations.
Let us assume that $\sigma^2 =1$, i.e., we will
consider now only standard Brownian motion.
Suppose that $\{V_k\}_{k\geq 1}$ are i.i.d.
random variables such that $|V_k| \leq \beta$, a.s.
Let $S_n = \sum_{k=1}^n V_k$ and extend the function
$n \to S_n$ to all positive real values 
by linear interpolation between $S_n$ and $S_{n+1}$.
Note that the random function $S_t$
is Lipschitz with constant $\beta$.

The following is an immediate consequence of Corollary 5.6.

\proclaim Theorem 5.7. Suppose that $V_k$ and $S_t$ are as above.
If $S_t$ and $B_t$ are constructed on the same probability
space (but not necessarily independent), then 
$$\limsup _{t\to \infty} { S_t - B_t \over \log t}
\geq  {1 /(  4\beta)}. \eqno(5.1)$$
\ms

Theorem 2.3.2 of Cs\"org\"o and R\'ev\'esz (1981)
says that if the $V_k$ have finite variance and
$$\limsup _{t\to \infty} { |S_t - B_t| \over \log t}
=0,$$
then the $V_k$ have a standard normal distribution.
Our result (5.1) may be interpreted as a quantitative
version of the same theorem, in the case when
$|V_k|$ are bounded. A remarkable
theorem of Koml\'os, Major and Tusn\'ady
(see Cs\"org\"o and R\'ev\'esz (1981) Theorem 2.6.1)
implies that if the $V_k$ are bounded, then one may
construct $S_t$ and $B_t$ on a common probability
space so that 
$$\limsup _{t\to \infty} { |S_t - B_t| \over \log t}
\leq C <\infty.\eqno(5.2)$$
It is striking that one can achieve the same logarithmic
order of approximation for a Lipschitz function
$S_t$ with independent increments $S_n- S_{n-1}$, as for
an arbitrary Lipschitz function $g(t)$ with constant $\beta$.
Rio (1991) proved that (5.2) holds with $C = 9/2$
if $V_k$ are centered Poisson variables
(the estimate had appeared in Section 5 of the preprint;
that section was not included in the final version
of the article, Rio (1994)). No other estimates
for $C$ seem to be known so (5.1) is our own modest
contribution to the field of strong approximations.

\bigskip
\noindent{\bf 6. Open problems}.
We list a few questions we were not able
to answer in this paper.

\item{(i)} Can one prove pathwise uniqueness
in Theorem 2.1 if one or both $\alpha_1$ and $\alpha_2$
belong to $(-1,0)$?

\item{(ii)} Does a result analogous to Theorem 3.7
hold for $\beta_1,\beta_2>0$ with $\beta_1-\beta_2 <0$?
A similar question can be asked about the case when
$\beta_1 < 0 < \beta_2$; in the last case a special
solution to (1.1), defined in Lemma 5.2, would have
to play an important role.
Can one generalize Theorem 3.7 to local times corresponding
to solutions of (1.3) with $\alpha_1$ and $\alpha_2$
not necessarily equal to $0$?

\item{(iii)} Find the best 
$\gamma = \gamma(\alpha_1,\alpha_2,\beta_1,\beta_2)>0$
in Lemma 4.2.

\item{(iv)} Find the best constants in (5.1) and (5.2).

\item{(v)} Does there exist a unique Lipschitz solution to (2.3)
if $B_t$ is a fractional Brownian motion of index $H\in (1/2,1)$?

\bigskip

\vskip1truein
\centerline{REFERENCES}
\bigskip

\item{[1]} R.J.~Adler (1981) {\it The Geometry of Random
Fields}. Wiley, New York.

\item{[2]} R.F.~Bass (1995) {\it Probabilistic Techniques in
Analysis}. Springer, New York.

\item{[3]} R.F.~Bass (1997) {\it Diffusions and Elliptic Operators}.
Springer, New York.

\item{[4]} J.~Bertoin (1996) {\it L\'evy Processes}.
Cambridge University Press, Cambridge.

\item{[5]} R.~Blumenthal, R. (1992)  {\it Excursions of
Markov Processes}. Birkh\"auser, Boston, Mass.

\item{[6]} K.~Burdzy (1987) {\it Multidimensional
Brownian Excursions and Potential Theory}. \break
Longman, Essex, England.

\item{[7]} K.~Burdzy, D.~Frankel and A.~Pauzner (1997)
``Fast equilibrium selection by rational
players living in a changing world''
(preprint)

\item{[8]} K.~Burdzy, D.~Frankel and A.~Pauzner (1998)
``On the time and direction of stochastic bifurcation''
In 
{\it Asymptotic Methods in Probability and Statistics.
A Volume in Honour of Mikl\'os Cs\"org\"o}.
Elsevier. (to appear)

\item{[9]} M.~Cs\"org\"o and P.~R\'ev\'esz (1981)
{\it Strong Approximations in Probability and Statistics}
Academic Press, New York.

\item{[10]} D.A.~Darling and P.~Erd\"os (1956)
``A limit theorem for the maximum of normalized sums
of independent random variables''
{\it Duke Math. J. \bf 23}, 143--155.

\item{[11]} L.~Decreusefond and A.S.~\"Ust\"unel (1997)
``Stochastic analysis of the fractional Brownian motion''
{\it Pot. Anal.}, to appear.

\item{[12]} E.~Fabes and C.E.~Kenig (1981) ``Examples of
singular parabolic measures and singular transition probability
densities''
{\it Duke Math.~J. \bf 48}, 845--856.

\item{[13]} H.~F\"ollmer, P.~Protter and A.~Shiryaev (1995)
``Quadratic covariation and an extension of It\^o's formula''
{\it Bernoulli \bf 1}, 149--169.

\item{[14]} J.M.~Harrison and L.A.~Shepp, L. A. (1981)
``On skew Brownian motion'' 
{\it Ann. Probab. \bf 9}, 309--313. 

\item{[15]} I.~Karatzas and S.E.~Shreve (1988)
{\it Brownian Motion and Stochastic Calculus}. \break
Springer, New York.

\item{[16]} S.~Karlin and H.M.~Taylor (1975)
{\it A First Course in Stochastic Processes}.
Academic Press, New York, 2-nd ed.

\item{[17]} S.~Karlin and H.M.~Taylor (1981)
{\it A Second Course in Stochastic Processes}.
Academic Press, New York.

\item{[18]} F.B.~Knight (1981)
{\it Essentials of Brownian Motion and Diffusion}.
Math. Surveys 18. American Mathematical Society,
Providence, RI.

\item{[19]} C.~Leuridan (1998) ``Le th\'eor\`eme de
Ray-Knight \`a temps fixe.''
{\it S\'eminaires de Probabilit\'es XXXII} (to appear).

\item{[20]} B.~Maisonneuve (1975)  ``Exit Systems.''
{\it Ann. Probab. \bf 3}, 399-411.

\item{[21]} A.~Mandelbaum, L.~Shepp and R.~Vanderbei (1990)
``Optimal switching between a pair of Brownian motions''
{\it Ann. Probab. \bf 18}, 1010--1033.

\item{[22]} J.R.~Norris, L.C.G.~Rogers and D.~Williams
(1987) ``Self-avoiding random walk: A Brownian motion
model with local time drift''
{\it Probab. Th. Rel. Fields \bf 74}, 271-287.

\item{[23]} D.~Revuz and M.~Yor (1991)
{\it Continuous Martingales and Brownian Motion}. \break
Springer, New York.

\item{[24]} E.~Rio (1991) ``Local invariance principles
and its applications to density estimation.''
{\it Pr\'epubl. Math. Univ. Paris-Sud \bf 91-71}.

\item{[25]} E.~Rio (1994) ``Local invariance principles
and their application to density estimation.''
{\it Prob. Th. Rel. Fields \bf 98}, 21--45.

\item{[26]} L.C.G.~Rogers (1997)
``Arbitrage with fractional Brownian motion''
{\it Math. Finance \bf 7}, 95--105.

\item{[27]} M.~Sharpe (1989)  {\it General Theory of Markov
Processes}. Academic Press, New York.

\item{[28]} D.W.~Stroock and S.R.S.~Varadhan (1979)
{\it Multidimensional Diffusions Processes}. Springer, New York.

\item{[29]} M.~Yor (1997) {\it Some Aspects of Brownian Motion.
Part II: Some Recent Martingale Problems}.
Birkh\"auser, Basel.

\vskip1truecm

\baselineskip=0.7\baselineskip
\parskip=0pt
\hbox{\vbox{
\obeylines
Department of Mathematics
University of Washington
Box 354350
Seattle, WA 98195-4350
bass@math.washington.edu
burdzy@math.washington.edu
}
}

\bye